\documentclass[12pt]{amsart}
\usepackage{amssymb}
\usepackage[all]{xy}
\usepackage[toc,page,title,titletoc,header]{appendix}
\usepackage{graphicx}
\usepackage{subfigure}

\makeindex

\begin{document}
\theoremstyle{plain}
\newtheorem{thm}{Theorem}[section]
\newtheorem*{thm1}{Theorem 1}
\newtheorem*{thm1.1}{Theorem 1.1}
\newtheorem*{thmM}{Main Theorem}
\newtheorem*{thmA}{Theorem A}
\newtheorem*{thm2}{Theorem 2}
\newtheorem{lemma}[thm]{Lemma}
\newtheorem{lem}[thm]{Lemma}
\newtheorem{cor}[thm]{Corollary}
\newtheorem{pro}[thm]{Proposition}
\newtheorem{propose}[thm]{Proposition}
\newtheorem{variant}[thm]{Variant}
\theoremstyle{definition}
\newtheorem{notations}[thm]{Notations}
\newtheorem{rem}[thm]{Remark}
\newtheorem{rmk}[thm]{Remark}
\newtheorem{rmks}[thm]{Remarks}
\newtheorem{defi}[thm]{Definition}
\newtheorem{exe}[thm]{Example}
\newtheorem{claim}[thm]{Claim}
\newtheorem{ass}[thm]{Assumption}
\newtheorem{prodefi}[thm]{Proposition-Definition}
\newtheorem{que}[thm]{Question}
\newtheorem{con}[thm]{Conjecture}

\newtheorem*{dmlcon}{Dynamical Mordell-Lang Conjecture}
\newtheorem*{condml}{Dynamical Mordell-Lang Conjecture}
\numberwithin{equation}{section}
\newcounter{elno}                
\def\points{\list
{\hss\llap{\upshape{(\roman{elno})}}}{\usecounter{elno}}}
\let\endpoints=\endlist
\newcommand{\res}{\rm res}
\newcommand{\Om}{\Omega}
\newcommand{\om}{\omega}
\newcommand{\la}{\lambda}
\newcommand{\mc}{\mathcal}
\newcommand{\mb}{\mathbb}
\newcommand{\surj}{\twoheadrightarrow}
\newcommand{\inj}{\hookrightarrow}
\newcommand{\zar}{{\rm zar}}
\newcommand{\Exc}{\rm Exc}
\newcommand{\an}{{\rm an}}
\newcommand{\red}{{\rm red}}
\newcommand{\codim}{{\rm codim}}
\newcommand{\Supp}{{\rm Supp}}
\newcommand{\rank}{{\rm rank}}
\newcommand{\Ker}{{\rm Ker \ }}
\newcommand{\Pic}{{\rm Pic}}
\newcommand{\Div}{{\rm Div}}
\newcommand{\Hom}{{\rm Hom}}
\newcommand{\im}{{\rm im}}
\newcommand{\Spec}{{\rm Spec \,}}
\newcommand{\Nef}{{\rm Nef \,}}
\newcommand{\Frac}{{\rm Frac \,}}
\newcommand{\Sing}{{\rm Sing}}
\newcommand{\sing}{{\rm sing}}
\newcommand{\reg}{{\rm reg}}
\newcommand{\Char}{{\rm char}}
\newcommand{\Tr}{{\rm Tr}}
\newcommand{\ord}{{\rm ord}}
\newcommand{\id}{{\rm id}}
\newcommand{\NE}{{\rm NE}}
\newcommand{\Gal}{{\rm Gal}}
\newcommand{\Min}{{\rm Min \ }}
\newcommand{\Max}{{\rm Max \ }}
\newcommand{\Alb}{{\rm Alb}\,}
\newcommand{\GL}{{\rm GL}\,}        
\newcommand{\PGL}{{\rm PGL}\,}
\newcommand{\Bir}{{\rm Bir}}
\newcommand{\Aut}{{\rm Aut}}
\newcommand{\End}{{\rm End}}
\newcommand{\Per}{{\rm Per}\,}
\newcommand{\ie}{{\it i.e.\/},\ }
\newcommand{\niso}{\not\cong}
\newcommand{\nin}{\not\in}
\newcommand{\soplus}[1]{\stackrel{#1}{\oplus}}
\newcommand{\by}[1]{\stackrel{#1}{\rightarrow}}
\newcommand{\longby}[1]{\stackrel{#1}{\longrightarrow}}
\newcommand{\vlongby}[1]{\stackrel{#1}{\mbox{\large{$\longrightarrow$}}}}
\newcommand{\ldownarrow}{\mbox{\Large{\Large{$\downarrow$}}}}
\newcommand{\lsearrow}{\mbox{\Large{$\searrow$}}}
\renewcommand{\d}{\stackrel{\mbox{\scriptsize{$\bullet$}}}{}}
\newcommand{\dlog}{{\rm dlog}\,}    
\newcommand{\longto}{\longrightarrow}
\newcommand{\vlongto}{\mbox{{\Large{$\longto$}}}}
\newcommand{\limdir}[1]{{\displaystyle{\mathop{\rm lim}_{\buildrel\longrightarrow\over{#1}}}}\,}
\newcommand{\liminv}[1]{{\displaystyle{\mathop{\rm lim}_{\buildrel\longleftarrow\over{#1}}}}\,}
\newcommand{\norm}[1]{\mbox{$\parallel{#1}\parallel$}}
\newcommand{\boxtensor}{{\Box\kern-9.03pt\raise1.42pt\hbox{$\times$}}}
\newcommand{\into}{\hookrightarrow}
\newcommand{\image}{{\rm image}\,}
\newcommand{\Lie}{{\rm Lie}\,}      
\newcommand{\CM}{\rm CM}
\newcommand{\sext}{\mbox{${\mathcal E}xt\,$}}  
\newcommand{\shom}{\mbox{${\mathcal H}om\,$}}  
\newcommand{\coker}{{\rm coker}\,}  
\newcommand{\sm}{{\rm sm}}
\newcommand{\pgcd}{\text{pgcd}}
\newcommand{\trd}{\text{tr.d.}}
\newcommand{\tensor}{\otimes}
\renewcommand{\iff}{\mbox{ $\Longleftrightarrow$ }}
\newcommand{\supp}{{\rm supp}\,}
\newcommand{\ext}[1]{\stackrel{#1}{\wedge}}
\newcommand{\onto}{\mbox{$\,\>>>\hspace{-.5cm}\to\hspace{.15cm}$}}
\newcommand{\propsubset}
{\mbox{$\textstyle{
\subseteq_{\kern-5pt\raise-1pt\hbox{\mbox{\tiny{$/$}}}}}$}}
\newcommand{\sA}{{\mathcal A}}
\newcommand{\sB}{{\mathcal B}}
\newcommand{\sC}{{\mathcal C}}
\newcommand{\sD}{{\mathcal D}}
\newcommand{\sE}{{\mathcal E}}
\newcommand{\sF}{{\mathcal F}}
\newcommand{\sG}{{\mathcal G}}
\newcommand{\sH}{{\mathcal H}}
\newcommand{\sI}{{\mathcal I}}
\newcommand{\sJ}{{\mathcal J}}
\newcommand{\sK}{{\mathcal K}}
\newcommand{\sL}{{\mathcal L}}
\newcommand{\sM}{{\mathcal M}}
\newcommand{\sN}{{\mathcal N}}
\newcommand{\sO}{{\mathcal O}}
\newcommand{\sP}{{\mathcal P}}
\newcommand{\sQ}{{\mathcal Q}}
\newcommand{\sR}{{\mathcal R}}
\newcommand{\sS}{{\mathcal S}}
\newcommand{\sT}{{\mathcal T}}
\newcommand{\sU}{{\mathcal U}}
\newcommand{\sV}{{\mathcal V}}
\newcommand{\sW}{{\mathcal W}}
\newcommand{\sX}{{\mathcal X}}
\newcommand{\sY}{{\mathcal Y}}
\newcommand{\sZ}{{\mathcal Z}}
\newcommand{\A}{{\mathbb A}}
\newcommand{\B}{{\mathbb B}}
\newcommand{\C}{{\mathbb C}}
\newcommand{\D}{{\mathbb D}}
\newcommand{\E}{{\mathbb E}}
\newcommand{\F}{{\mathbb F}}
\newcommand{\G}{{\mathbb G}}
\newcommand{\HH}{{\mathbb H}}
\newcommand{\I}{{\mathbb I}}
\newcommand{\J}{{\mathbb J}}
\newcommand{\M}{{\mathbb M}}
\newcommand{\N}{{\mathbb N}}
\renewcommand{\P}{{\mathbb P}}
\newcommand{\Q}{{\mathbb Q}}
\newcommand{\R}{{\mathbb R}}
\newcommand{\T}{{\mathbb T}}
\newcommand{\U}{{\mathbb U}}
\newcommand{\V}{{\mathbb V}}
\newcommand{\W}{{\mathbb W}}
\newcommand{\X}{{\mathbb X}}
\newcommand{\Y}{{\mathbb Y}}
\newcommand{\Z}{{\mathbb Z}}

\newcommand{\fix}{\mathrm{Fix}}

\newcommand{\SH}{\rm SH}
\newcommand{\Tan}{\rm Tan}

\title[The Dynamical Mordell-Lang Conjecture]{The Dynamical Mordell-Lang Conjecture for polynomial endomorphisms of the affine plane}

\author{Junyi Xie}
\address{Universit\'e de Rennes I
  Campus de Beaulieu,
  b\^atiment 22-23,
  35042 Rennes cedex
  France}
%
\email{junyi.xie@univ-rennes1.fr}
\date{\today}
\bibliographystyle{plain}
\thanks{The author is supported
by the ERC-starting grant project "Nonarcomp" no.307856.}

\maketitle

\begin{abstract}
In this paper we prove the Dynamical Mordell-Lang Conjecture for polynomial endomorphisms of the affine plane.
\end{abstract}

\tableofcontents

\newpage

%
%
%
%
%
%
%
%
%
%
%
%
%
%
%
%
%
%

\section*{Introduction}
\subsection{The dynamical Mordell Lang  conjecture}
This article is concerned with the so-called {\em dynamical Mordell-Lang conjecture} \index{dynamical Mordell-Lang conjecture}that was proposed by Ghioca and Tucker in~\cite{Ghioca2009}.
\begin{dmlcon}[\cite{Ghioca2009}]\label{dml}Let $X$ be a quasi-pro\-jective variety defined over $\mathbb{C}$, let
$f: X \rightarrow X$ be an endomorphism, and $V$ be any subvariety of $X$. For any point $p\in X(\mathbb{C})$  the set $\{n\in \mathbb{N}|\,\,f^n(p)\in V(\mathbb{C})\}$ is a union of at most finitely many arithmetic progressions\footnote{an arithmetic progression is a set of the form $\{an+b|\,\,n\in \mathbb{N}\}$ with $a,b\in \mathbb{N}$.}.
\end{dmlcon}

This conjecture is inspired by the
Mordell-Lang conjecture on subvarieties of semiabelian varieties (now a theorem of Faltings \cite{Faltings1994} and Vojta \cite{Vojta1996}), which says that if $V$ is a subvariety of a semiabelian variety $G$ defined over $\C$ and $\Gamma$ is a finitely generated subgroup  of $G(\C)$, then $V(\C)\bigcap \Gamma$ is a union of at most finitely many translates of
subgroups of $\Gamma$.

Observe that the dynamical Mordell-Lang conjecture implies the classical Mordell-Lang conjecture in the case $\Gamma\simeq (\mathbb{Z},+)$.

\medskip

It is also motivated by the  Skolem-Mahler-Lech Theorem \cite{Lech1953} on linear recurrence sequences.
More precisely, suppose $\{A_n\}_{n\geq 0}$ is any recurrence sequence satisfying $A_{n+l}=F(A_{n},\cdots,A_{n+l-1})$ for all $n\geq 0$, where $l\geq 1$ and $F(x_0,\cdots,x_l)=\sum_{i=0}^{l-1}a_{i}x_i$ is a linear form on $\C^{l}$. The  Skolem-Mahler-Lech Theorem asserts that the set $\{n\geq 0|\,\,A_n=0\}$ is a union of at most finitely many arithmetic progressions.

This statement is equivalent to the dynamical Mordell-Lang conjecture for the linear map $f:(x_0,\cdots,x_{l-1})\mapsto(x_1,\cdots,x_{l-1},F(x_0,\cdots,x_l))$
and the hyperplane $V=\{x_0=0\}.$

%

%
%
%
%
%
%
%
%

\smallskip

\subsection{The main results and comparison to previous results}
Our  goal is to prove this conjecture for {\em any} polynomial endomorphism on $\mathbb{A}^2_{\overline{\mathbb{Q}}}$.
\begin{thm}\label{thmdmlpoly}Let $f:\mathbb{A}^2_{\overline{\mathbb{Q}}}\rightarrow \mathbb{A}^2_{\overline{\mathbb{Q}}}$ be a polynomial endomorphism defined over $\overline{\mathbb{Q}}$. Let $C$ be an irreducible curve in $\mathbb{A}^2_{\overline{\mathbb{Q}}}$ and $p$ be a closed point in $\mathbb{A}^2_{\overline{\mathbb{Q}}}$. Then the set $\{n\in \mathbb{N}|\,\,f^n(p)\in C\}$ is a finite union of arithmetic progressions.
\end{thm}
Pick any polynomial $F(x,y)\in \overline{\mathbb{Q}}[x,y]$.
By applying this result to the map $f:\mathbb{A}^2_{\overline{\mathbb{Q}}}\rightarrow \mathbb{A}^2_{\overline{\mathbb{Q}}}$ defined by $(x,y)\mapsto (y,F(x,y))$ and $C=\{x=0\}$, we  obtain the following corollary about recurrence sequences.
\begin{cor}\label{correncurrenceaeq}Let $\{A_n\}_{n\geq 0}$ be a sequence of algebraic numbers satisfying $A_{n+2}=F(A_n,A_{n+1})$ for all $n\geq 0$, where $F(x,y)\in \overline{\mathbb{Q}}[x,y]$. Then the set $\{n\geq 0|\,\, A_n=0\}$ is a finite union of arithmetic progressions.
\end{cor}
A  direct induction on the dimension also yields the following
\begin{thm}\label{thmdmlmsv}
For any  non-constant polynomials $F_1, \ldots, F_m \in \overline{\mathbb{Q}}[T]$, let us consider the endomorphism $f:=(F_1(x_1),\cdots,F_m(x_m))$ on $\mathbb{A}^m_{\overline{\Q}}$.

For any irreducible curve $C \subset \mathbb{A}^m_{\overline{\Q}}$ defined over $\overline{\mathbb{Q}}$ and any point $ p \in \mathbb{A}^m(\overline{\mathbb{Q}})$, the set $\{n\geq 0|f^n(p)\in C\}$ is a finite union of arithmetic progressions.
\end{thm}

The dynamical Mordell-Lang conjecture has received quite a lot of attention  in the recent years and our theorems are closely related to several known results.

Bell, Ghioca and Tucker \cite{Bell2010} proved the Dynamical Mordell-Lang conjecture for \'etale maps of quasiprojective varieties of arbitrary dimension, thereby generalizing the Skolem-Mahler-Lech Theorem \cite{Lech1953} on linear recurrence sequences. The core of their argument is to
work in a $p$-adic field and to analyze the dynamics in a quasi-periodic region where they are
able to construct suitable invariant  curves.
Afterwards, the author \cite{Xie2014} proved the dynamical Mordell-Lang conjecture for birational endomorphisms of the affine plane. The techniques in \cite{Xie2014} are of a very different flavour. Particularly, in \cite{Xie2014}, we got a new proof of the dynamical Mordell-Lang conjecture for polynomial automorphisms of $\A^2$ which are {\em not} conjugated to an automorphism of some projective surface. However, we relied on Bell, Ghioca and Tucker's result in some cases, especially in the case of affine automorphisms of $\A^2$. In this paper, we develop the techniques used in \cite{Xie2014} in a more general situation and use them more systematically.


Our  Theorem~\ref{thmdmlmsv} also generalizes \cite[Theorem 1.5]{Benedetto2012} of Benedetto, Ghioca, Kurlberg, and Tucker (hence \cite[Theorem 1.4]{Ghioca2008} of  Ghioca, Tucker, and Zieve) which proved the Dynamical Mordell-Lang conjecture
in the case
$f=(F(x_1),\cdots, F(x_n)):\mathbb{A}^n_{\overline{\mathbb{Q}}}\rightarrow \mathbb{A}^n_{\overline{\mathbb{Q}}}$
where $F\in \overline{\mathbb{Q}}[t]$ is an indecomposable polynomial function defined over $\overline{\mathbb{Q}}$ which
has no periodic critical points other than the point at infinity and $V$ is a curve.

%

\bigskip

\subsection{Overview of the proof of the main theorem}
Since the proof of Theorem \ref{thmdmlpoly}  is quite long and involves many different cases, we provide in this introduction a detailed overview of  our strategy.
To simplify the discussion, we suppose that $f$ is a dominant polynomial map $f:=(F(x,y),G(x,y))$ defined over $\Z$ and  $p\in \Z^2$. We assume that the set $\{n\geq 0|\,\,f^n(p)\in C\}$ is infinite and $p$ is not preperiodic. We need to prove that the curve $C$ is periodic.



\smallskip

To do so we shall work in  suitable compactifications of $\A^2$ for which the induced map by $f$ at infinity is nice, in the sense that it does not contract any curve to a point of indeterminacy. These dynamically meaningful compactifications have been constructed and studied by Favre and Jonsson in~\cite{Favre2011}. To put it in broad terms, we shall use suitable height arguments to  focus what happens to the branches of $C$ at infinity under iteration, and conclude by applying  the construction of  polynomials in valuation subrings of $\overline{\mathbb{Q}}[x,y]$ that we have developped in a former paper \cite{Xieb}.

\smallskip

Let us now see in more details how our arguments work. We denote by $\la_2$ the topological degree \index{topological degree} of $f$ i.e. the number of preimages of a general point in $\mathbb{A}^2(\overline{\mathbb{Q}})$ and by $\la_1$ the dynamical degree \index{dynamical degree} of $f$ that is $\lim_{n\rightarrow \infty}(\deg f^n)^{\frac{1}{n}}$. These degrees are invariants of conjugacy and satisfy the inequality $\la_1^2\geq \la_2$.

\smallskip

\noindent {\em 1) The case $\la_1^2=\la_2$}. This case is quite special in the sense that the map $f$ exhibits some kind of dynamical rigidity. By
\cite[Theorem C]{Favre2011} either one can find a projective compactification of $\A^2$ in which $f$ induces an {\em endomorphism}, or $f$ is a skew product and there exists affine coordinates in which it can be written as $(f(x,y) = (P(x), Q(x,y))$.

Let us first explain how our main theorem is proved in this case.
There are two important ingredients: one is Siegel's theorem that give constraints on the geometry of the curve $C$ and its preimages by $f$; and the other is a local version of the dynamical Mordell-Lang conjecture for super-attracting germs. The latter statement was first used in \cite{Xie2014} to treat the special case of birational polynomial maps, and we use it here more systematically.

{\em 1a) The map $f$ is an endomorphism} on a projective compactification $X$ of $\A^2$, with boundary $D_{\infty}:=X\setminus \A^2$. We proceed as follows. Since $C$ contains infinitely many points in the orbit of $p$, it also contains infinitely many integral
points hence admits at most  two branches at infinity by Siegel's theorem. Arguing in the same way with the preimages of $C$ we  end up with a sequence of irreducible curves $\{C^{j}\}_{j\leq 0}$ with $C^0=C$, and  $f(C^{j})=C^{j+1}$ such that  $C^j$ has at most two branches at infinity and the set $\{n\geq 0|\,\,f^n(p)\in C^j\}$ is infinite for all $j$.

Then we look at the positions of $C^j$ at infinity.  One can show that two situations may appear: either  one branch of $C$ intersects the divisor at infinity $D_\infty$ at a superattracting periodic point; or $C^j$ have bounded intersection with $D_\infty$. In the former situation we apply
our local version of the dynamical Mordell-Lang conjecture to conclude. In the latter case,
either $C^j = C^{j'}$ for some $j >j'$ and $C$ is periodic, or the $C^j$'s belong to
a fibration that is preserved by $f$ in which case it is not difficult to conclude.

{\em 1b) The map $f$ is a skew product and  $\deg(f^n) \simeq n \lambda_1^n$}.
One may construct a dynamically nice compactification $X$ of $\A^2$
such that $X$ is isomorphic to a Hirzebruch surface, and $f$ preserves the unique rational fibration on $X$. One can then understand fairly well the dynamics of $f$ on
the divisor at infinity in $X$, and the proof goes in a very similar way as in the previous case 1a).

\smallskip

\noindent {\em 2) The case $\la_1^2<\la_2$}.
To analyze this case  the above two ingredients are no longer sufficient, and we need to get deeper in the action of the map $f$ at infinity in dynamically meaningful compactifications of $\A^2$. In other words we shall use extensively the analysis of the action of $f$ on the space of valuations at infinity initiated in \cite{Favre2011}.

\smallskip


As in \cite{Favre2011}, $V_\infty$ is defined as the set of valuations  $v: k[x,y] \to \R \cup \{ +\infty\}$
centered at infinity and normalized by $\min \{ v(x), v(y)\} = -1$.
This set becomes a compact space when endowed with the topology of the pointwise convergence.
It can be also endowed with a natural partial order relation given by $v \le v'$ if and only if $v(P) \le v'(P)$ for all $P \in k[x,y]$.
This partial order relation makes it to be an $\R$-tree. The unique minimal point for that order relation is the valuation $-\deg$.

Let $s$ be a formal branch of curve centered at infinity. We may associate to $s$ a valuation $v_s\in V_{\infty}$ defined by $P\mapsto -\min\{\ord_{\infty}(x|_s),\ord_{\infty}(y|_s)\}^{-1}\ord_{\infty}(P|_s)$. Such a valuation is called a curve valuation.

Pick any proper modification $\pi : X\to \mathbb{P}^2$ that is an isomorphism above the affine plane with $X$ a smooth projective surface.
Let $\{E_0, E_1,\cdots,E_m\}$ be the set of all irreducible components of $X\setminus \mathbb{A}^2_k$. For any irreducible component $E_i$, we can
define a valuation $v_{E_i}:=b_{E_i}^{-1}\ord_{E_i}$ where $b_{E_i}:=-\min\{\ord_{E_i}(x),\ord_{E_i}(y)\}.$ Observe that $v_{E_i}\in V_{\infty}$. Such a valuation is called a divisorial valuation. The set of divisorial valuations is dense in any segment in $V_{\infty}$.

%

To define the action of $f$ on $V_{\infty}$, we first define a function $d(f,\d)$ on $V_{\infty}$ by $v\mapsto -\min\{v(f^*L),0\}$ where $L$ is a general linear form in $\overline{\mathbb{Q}}[x,y].$ For simplicity, we suppose that $d(f,v)>0$ for all $v\in V_{\infty}.$
Then the action $f_{\d}$ on $V_{\infty}$ is defined by $f_{\d}(v): P\mapsto d(f,v)^{-1}v(f^*P)$ for all $P\in \overline{\mathbb{Q}}[x,y].$

\smallskip

In \cite[Appendix A]{Favre2011} and essentially in \cite{boucksomfavrejonsson}, Boucksom, Favre and Jonsson constructed an eigenvaluation $v_*$ in $V_{\infty}$ and a canonical closed subset $J(f)$ of $V_{\infty}$ (see Section \ref{sectiongreenf} for details).  The following Theorem is a key ingredient in our paper.

\begin{thm}\label{thmmostimportantstepsimple}Let $f$ be a dominant polynomial endomorphism on $\mathbb{A}^2$ defined over an algebraically closed field satisfying $\la_1^2>\la_2$ and $\#J(f)\geq 3$. Let $W$ be an open neighborhood of $v_*$ in $V_{\infty}$. There exists a finite set of polynomials $\{P_i\}_{1\leq i\leq s}$ and a positive integer $N$ such that for any set of valuations $\{v_1,v_2\}\subseteq V_{\infty}\setminus f_{\d}^{-N}(W)$, there exists an index $i\in \{1,\cdots,s\}$ such that $v_j(P_i)> 0$ for all $j=1,2$.
\end{thm}

%

%
%

\smallskip
%

{\em 2a) The case $\#J(f)\geq 3.$}
We first suppose that $v_*$ is nondivisorial.
As in case $1a)$, we use Siegel's theorem to constructs a sequence of irreducible curves $\{C^{j}\}_{j\leq 0}$ with $C^0=C$, and  $f(C^{j})=C^{j+1}$ such that  $C^j$ has at most two branches at infinity and the set $\{n\geq 0|\,\,f^n(p)\in C^j\}$ is infinite for all $j\leq 0$. There exists a neighborhood $W$ of $v_*$ such that the curve valuations defined by the branches of $C$ at infinity are not contained in $W$ and $f_{\d}(W)\subseteq W$. It follows that for any $N\geq 0$, the curve valuations defined by the branches of $C^j$, $j\leq -N$, at infinity are not contained in $f^{-N}_{\d}(W).$ For $N$ large enough, Theorem \ref{thmmostimportantstepsimple} allows us to construct a finite set of polynomials $\{P_i\}_{1\leq i\leq s}$ such that for any $j\leq -N$, there exists $i=1,\cdots, s$ satisfying $P_i|_{C^j}=0$; this implies that $C$ is periodic.

When $v_*$ is divisorial, we may find a smooth projective compactification $X$ of $\A^2$ such that there exists an irreducible component $E$ of $X\setminus \A^2$ such that $v_*=v_E$. Take $W$ a small enough neighborhood of $v_*$. Not like the former case,  in general we can not ask $W$ to be invariant under $f_{\d}$. In this case we need Theorem \ref{thmmostimportantstep} which is a stronger version of Theorem \ref{thmmostimportantstepsimple}. Relying on Theorem \ref{thmmostimportantstep}, we can show that there is always some branch $s^j$ of $C^j$ such that the valuation $v_{s^j}$ stays in $W$ for a long time.

 In the case $\deg f|_E=1$, we can show that the intersection number $(s^j\cdot l_{\infty})$ of $s^j$ and the line $l_{\infty}$ at infinity in $\P^2$  can not grow much when $v_{s^j}$ stays in $W$. Also we can use Theorem \ref{thmmostimportantstep} to show that if a branch satisfying $v_{s^j}\not\in W$, then $(s^j\cdot l_{\infty})/\deg C^j$ is bounded by $1-\varepsilon$ for some $\varepsilon>0$ and all $j$ negative enough. By some very careful analysis
, we can at the end bound the degree of  $C^j$.

In the case $\deg f|_E\geq 2$, the new ingredient is the Northcott property for number fields. More precisely, for any number field $K$ such that both $E,f$ are defined over $K$, for any point $x\in E(K)$, the set of inverse orbit of $x$ in $E(K)$ is finite. Using this fact, we can show that if the branch $s^j$ of $C^j$ stays in $W$ in a long time, then the center of $s^j$ is contained in a periodic point in $E$. Then we can conclude by some local argument.


%
%
%
%
%
%
%
%
%
%
%
%
%
%
{\em 2b)} The case $\#J(f)\leq 2$. The serious difficulty in this case is that we can not apply Theorem \ref{thmmostimportantstep} directly.
If all valuations in $J(f)$ are divisorial, we may prove that $f$ is either \'etale or preserves a fibration. We treat this case separately. Otherwise, we notice that all the nondivisorial valuations in $ J(f)$ are periodic and repelling under $f_{\d}$. This fact shows that the curve valuations associated to the branches of $f^n(C)$ at infinity
can not be too close to those nondivisorial valuations in $J(f)$. This fact allows us to modify $\theta^*$ a little
and  get a modified version of Theorem \ref{thmmostimportantstep}.
Then we can use a strategy  which parallels to the corresponding case in $2a)$ to conclude our theorem in this case.

\bigskip


\subsection{More remarks about our techniques}
In order to prove Theorem \ref{thmdmlpoly}, we have developed some new techniques in this paper based on the theory of Favre and Jonsson (\cite{Favre2004,Favre2007,Favre2011,Jonsson}). These techniques can
be applied to not only the dynamical Mordell-Lang conjecture but also many other
problems of polynomial endomorphisms of $\A^2$. In particular, in our recent work \cite{M.Jonssona}, Jonsson, Wulcan and I proved \cite[Conjecture 1]{Silverman2014} of Silverman for polynomial endomorphisms of $\A^2$ with $\la_1\geq \la_2$ and in another recent work \cite{Jonssona}, Jonsson, Wulcan and I classified all the polynomial endomorphisms $f$ on $\A^2$ that preserves a pencil $|P|$ up to changing affine coordinates and replacing $f$ by a suitable iteration. In the sequels to the papers \cite{Xied,Xiee}, these techniques will also
be used to study the orbits of point, the periodic points and the periodic curves of polynomial endomorphisms $f$ on $\A^2$.

\subsection{Further problems}In our proof of Theorem \ref{thmdmlpoly}, we have use the theorem of Siegel and the Northcott property for number fields. That's why we restrict our theorem for endomorphisms defined over $k=\overline{\Q}$. We suspect that Theorem \ref{thmdmlpoly} remains true when $k$ is an arbitrary algebraically closed field of characteristic $0$. It seems that we can prove it by induction on transcendence degree of $k$ over $\overline{\Q}$ and some reduction arguments; however, the step of the reduction seems not trivial.

It would be interesting to generalize Theorem \ref{thmdmlpoly} for endomorphisms on arbitrary affine surfaces. It might be possible to prove this by methods similar to those in this paper. However this seems to require substantial effort, since it needs to generalize the theory of valuative space at infinity for $\A^2$ developed in \cite{Favre2007,Favre2011,Xieb} and this paper to arbitrary affine surfaces.

\subsection{The plan of the paper}
In Part \ref{partpre}, we gather a number of results on the geometry and dynamics at infinity.
We first introduce the valuative tree at infinity in Section 1, and then turn our attention to the notion of subharmonic function in Section 2. We give an interpretation of this potential theory in terms of $b$-divisors
in Section 3. Finally we recall the main properties of the action of a polynomial
map on the valuation space in Section 4.

In Part 2, we collect some arguments of local nature that will be used in the proof of our main result.
We first recall the definition and basic properties of the
local valuation space in Section 5. Then we state and prove a local version of the
dynamical Mordell-Lang conjecture for superattracting analytic germs in Section 6.

In Part 3, we give some basic observations on the Dynamical Mordell-Lang
Conjucture. We first define the DML property in Section 7.  Then we get some constraints on the target curve by Siegel's theorem in Section 8 and we use these constraints and the local arguments in Part 2 to prove Theorem \ref{thmdmlmsv} in Section 9.

In Part 4, we prove our main theorem in the resonant case $\la_1^2=\la_2$. We first treat the case $\deg(f^n)\asymp n\la_1^n$ in Section 10 and then treat the case $\deg(f^n)\asymp \la_1^n$ in Section 11.

In Part 5, we study the valuative dynamics in the case $\la_2^2>\la_1$. We first describe some basic properties of the Green function of $f$ on $V_{\infty}$ in Section 12. Then use the Green function to study the valuative dynamics in Section 13. In Section 14 we get more information on the valuative dynamics in the case $J(f)$ is finite. Finally, in Section 15 we show that $f$ is \'etale or preserves a fibration when $J(f)$ is a finite set of divisorial valuations.

In Part 6, we prove our main theorem in the non-resonant case $\la_1^2>\la_2$ which completes the proof of Theorem \ref{thmdmlpoly}. We first treat the case $\#J(f)\geq 3$ in Section 16 and then treat the case $\#J(f)\leq 2$ in Section 17.

\section*{Acknowledgement}
I would like to thank C. Favre for his constant support and his direction during the
writing of this article, and M. Jonsson for useful discussions on the valuative tree. I would also thank S.Cantat for their comments on the first version of this article.

\newpage

\part{Preliminaries}\label{partpre}

In this part, we collect some basic informations and results on the principal tools of our article, namely the space of valuations on the ring of polynomials in two variables that are centered at infinity. We first describe its tree structure in Section \ref{sectiontvaluationtreeinfty}, and then turn our attention to the notion of subharmonic function in Section \ref{sectionpotential}.  We give an interpretation of this potential theory in terms of $b$-divisors in Section \ref{sectionriezari}. Finally we recall the main properties of  the action of a polynomial map on the valuation space in Section \ref{sectionbgdopm}.

This part does not contain any new material. Proofs will be omitted and we shall refer to
other sources.

\smallskip

In this part  $k$ is an algebraically closed field of characteristic zero. We shall also fix affine coordinates on $\A^2_k = \Spec k [x,y]$.

\section{The valuative tree at infinity}\label{sectiontvaluationtreeinfty}
 We refer to \cite[Section 2]{Xieb} for details, see also \cite{Favre2004,Favre2011}.
\subsection{The valuative tree centered at infinity}

In this article by a valuation on a unitary $k$-algebra $R$ we shall understand a function
$v : R \to \mathbb{R}\cup \{+\infty\}$ such that the restriction of $v$ to $k^* = k - \{ 0 \}$ is
constant equal to $0$, and $v$ satisfies $v(fg) = v(f) + v(g)$ and $v(f +g) \ge \min \{v(f), v(g) \}$.
It is usually referred to as a pseudo-valuation in the literature, see \cite{Favre2004}. We will however make a slight abuse of notation and call them valuations.

\smallskip

We denote by $V_{\infty}$ \index{$V_{\infty}$}the space of all normalized valuations centered at infinity \index{normalized valuation centered at infinity} i.e. the set of valuations $v:k[x,y]\rightarrow \mathbb{R}\cup \{+\infty\}$ normalized by $\min\{v(x),v(y)\}=-1$. The topology on $V_{\infty}$ is defined to be the weakest topology making the map $v\mapsto v(P)$ continuous for every $P\in k[x,y]$.

\smallskip

The set $V_{\infty}$ is equipped with a {\em partial ordering} defined by $v\leq w$ if and only if
$v(P)\leq w(P)$ for all $P\in k[x,y]$ for which  $-\deg: P \mapsto -\deg (P)$ is the unique minimal element.

\smallskip

Given any valuation $v\in V_\infty$, the set
$ \{ w \in V_\infty, \, - \deg \le w \le v \}$ is isomorphic as a poset to the real segment
$[0,1]$ endowed with the standard ordering. In other words, $(V_\infty, \le)$ is a rooted tree in the sense of \cite{Favre2004,Jonsson}.

It follows that given any two valuations $v_1,v_2\in V_{\infty}$,
there is a unique valuation in $V_{\infty}$ which is maximal in the set $\{v\in V_{\infty}|\,\,v\leq v_1 \text{ and } v\leq v_2\}.$ We denote it by  $v_1\wedge v_2$.

The segment $[v_1, v_2]$ is by definition the union of $\{w , \, v_1\wedge v_2 \le w \le v_1\}$
and $\{w , \, v_1\wedge v_2 \le w \le v_2\}$.

\smallskip

Pick any valuation  $v\in V_\infty$. We say that two points $v_1, v_2$
lie in the same direction at $v$ if the segment $[v_1, v_2]$ does not contain $v$.
A {\em direction} (or a tangent vector) \index{direction} at $v$ is an equivalence class for this relation.
We write $\Tan_v$\index{$\Tan_v$} for the set of directions at $v$.

\smallskip

When $\Tan_v$ is a singleton, then $v$ is called an endpoint\index{endpoint}. In $V_\infty$, the set of endpoints is exactly the set of all maximal valuations. 
When $\Tan_v$ contains exactly two directions, then $v$ is said to be regular\index{regular point}.
When $\Tan_v$ has more than three directions, then $v$ is a branched point\index{branched point}.

\smallskip

Pick any $v\in V_\infty$.
For any tangent vector $\vec{v}\in \Tan_v$, we denote by $U(\vec{v})$ the subset of those elements in $V_\infty$ that determine $\vec{v}$. This is an open set whose boundary is reduced to the singleton $\{v\}$. If $v\neq -\deg$, the complement of $\{w\in V_\infty, \, w \ge v\}$ is equal to $U(\vec{v}_0)$ where $\vec{v}_0$ is the tangent vector determined by $-\deg$.

It is a fact that finite intersections of open sets of the form $U(\vec{v})$\index{$U(\vec{v})$} form a basis for the topology of $V_\infty$.

\smallskip

The \emph{convex hull}\index{convex hull} of  any subset $S\subset V_\infty$ is defined as the set of valuations $v\in V_\infty$ such that
there exists a pair $v_1 , v_2 \in S$ with $v\in [v_1, v_2]$.

A \emph{finite subtree} of $V_\infty$ is, by definition, the convex hull of a finite collection of points in $V_\infty$. A point in a finite subtree  $T\subseteq V_{\infty}$ is said to be an end point if it is extremal in $T.$

\subsection{Compactifications of $\mathbb{A}^2_k$}

A {\em compactification}\index{compactification} of $\A^2_k$ is the data of a projective surface $X$ together with an open immersion $\A^2_k \to X$ with dense image.

A compactification $X$ dominates another one $X'$ if the canonical birational map $X \dashrightarrow X'$ induced by the inclusion of $\A^2_k$ in both surfaces is in fact a regular map.

The category $\mathcal{C}$\index{$\mathcal{C}$} of all compactifications of $\A^2_k$ forms an inductive system for the relation of domination.

\smallskip

Recall that we have fixed affine coordinates on $\A^2_k = \Spec k[x,y]$.
We let $\P^2_k$ be the standard compactification of $\A^2_k$ and denote by
$l_\infty := \P^2_k \setminus \A^2_k$\index{$l_\infty$} the line at infinity in the projective plane.

An {\em admissible compactification}\index{admissible compactification} of $\A^2_k$ is by definition a smooth projective surface $X$ endowed with a birational morphism $\pi_X : X \to \P^2_k$ such that $\pi_X$ is an isomorphism over $\A^2_k$ with the embedding $\pi^{-1}|_{\A_k^2}:\A_k^2\to X$. Recall that $\pi_X$ can then be decomposed as a finite sequence of point blow-ups.

We shall denote by $\mathcal{C}_0$\index{$\mathcal{C}_0$} the category of all admissible compactifications. It is also an inductive system for the relation of domination. Moreover $\mathcal{C}_0$ is a subcategory of $\mathcal{C}$ and for any compactification $X\in \mathcal{C}$, there exists $X'\in \mathcal{C}_0$ dominates $X$.

\subsection{Divisorial valuations}
Let $X\in \mathcal{C}$ be a compactification of $\A^2_k=\Spec k[x,y]$ and $E$ be an irreducible component of $X\setminus \mathbb{A}^2$. Denote by $b_E:=\min\{\ord_E(x),\ord_E(y)\}$ \index{$b_E$} and $v_E:=b_E^{-1}\ord_E$\index{$v_E$}. Then we have $v_E\in V_{\infty}$.

By Poincar\'e Duality there exists a unique
{\em dual divisor} \index{dual divisor} $\check{E}$ \index{$\check{E}$} of $E$ i.e. the unique divisor supported by $X\setminus \A^2$ such that $(\check{E}\cdot F)=\delta_{E,F}$ for all irreducible components $F$ of $X\setminus \A^2$.


\subsection{Classification of valuations}
There are four kinds of valuations in $V_{\infty}$. The first one corresponds to the {\em divisorial valuations}\index{divisorial valuation} which we have defined above. We now describe the three remaining types of valuations.
\subsubsection*{Irrational valuations}
Consider any two irreducible components $E$ and $E'$ of $X\setminus \A^2_k$ for some compactification $X\in \mathcal{C}$ of $\mathbb{A}^2_k$
intersecting at a point $p$. There exists a local coordinate $(z,w)$ at $p$ such that $E=\{z=0\}$
and $E'=\{w=0\}$. To any pair $(s,t)\in (\mathbb{R}^{+})^{2}$ satisfying  $sb_E+tb_{E'}=1$, we attach the valuation $v$ defined on the
ring $O_p$ of holomorphic germs at $p$ by
$v(\sum a_{ij}z^iw^j)=\min \{si+tj|\,\, a_{ij}\neq 0\}.$
Observe that it does not
depend on the choice of coordinates. By first extending $v$ to the common fraction
field $k(x,y)$ of $O_p$ and $k[x,y]$, then restricting it to $k[x,y]$, we obtain a valuation in $V_{\infty}$, called
{\em quasimonomial}\index{quasimonomial}.  It is divisorial if and only if either $t=0$ or the ratio $s/t$ is a
rational number. Any divisorial valuation is quasimonomial. An {\em irrational valuation}\index{irrational valuation} is
by definition a nondivisorial quasimonomial valuation.
\subsubsection*{Curve valuations}
Recall that $l_{\infty}$ is the line at infinity of $\mathbb{P}^2_k$.
For any formal curve $s$ centered at some point $q\in l_{\infty}$, denote by $v_s$\index{$v_s$} the valuation defined by $P\mapsto (s\cdot l_{q})\ord_{\infty}(P|_s)$. Then we have $v_s\in V_{\infty}$ and call it a {\em curve valuation}\index{curve valuation}.

\medskip

Let $C$ be an irreducible curve in $\P^2_k$. For any point $q\in C\cap l_{\infty}$, denote by $O_{q}$ the local ring at $q$, $m_q$ the maximal ideal of $O_q$ and $I_C$ the ideal of height 1 in $O_q$ defined by $C$. Denote by $\widehat{O}_{q}$ the completion of $O_{q}$ w.r.t. $m_q$, $\widehat{m}_q$ the completion of $m_q$ and $\widehat{I}_C$ the completion of $I_C$. For any prime ideal $\widehat{p}$ of height 1 containing $\widehat{I}_C$, the morphism $\Spec \widehat{O}_q/\widehat{p}\rightarrow\Spec \widehat{O}_{q}$ defines a formal curve centred at $q$. Such a formal curve is called a \emph{branch of $C$ at infinity}\index{branch of $C$ at infinity}.

\begin{figure}
\centering
\includegraphics[width=6.5cm]{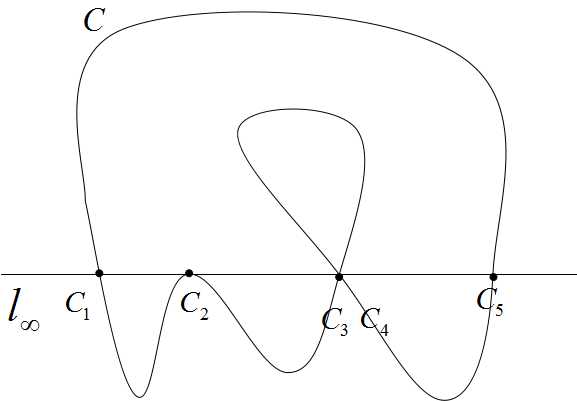}
\caption{}
\label{figc0}
\end{figure}

For example, in Figure \ref{figc0}, there are five branches at infinity of $C$.
Then for any branch $C_i$ $i=1,\cdots,5$, of $C$ at infinity, it corresponds to a curve valuation $v_{C_i}$ $i=1,\cdots,5$.

\subsubsection*{Infinitely singular valuations}
Let $h$ be a formal series of the form $h(z)=\sum_{k=0}^{\infty} a_kz^{\beta_k}$ with $a_k\in k^*$ and $\{\beta\}_k$ an increasing sequence of rational numbers
with unbounded denominators. Then $P\mapsto -\min\{\ord_{\infty}(x),\ord_{\infty}(h(x^{-1}))\}^{-1}\ord_{\infty}P(x,h(x^{-1}))$ defines
a valuation in $V_{\infty}$ namely an infinitely singular valuation\index{infinitely singular valuation}.

%
%
%
%

\medskip

A valuation $v\in V_{\infty}$ is a branch point in $V_{\infty}$ if and only if it is diviorial, it is a regular point in $V_{\infty}$ if and only if it is an irrational valuation, and it is an endpoint in $V_{\infty}$ if and only if it is a curve valuation or an infinitely singular valuation.
Moreover, given any smooth projective compactification $X$ in which $v=v_E$, one proves that the
map sending an element $V_\infty$ to its center in $X$ induces a  map $\Tan_v \to E$ that is a bijection.

\subsection{Parameterizations}
The {\em skewness} function\index{skewness function}
$\alpha:V_{\infty}\rightarrow [-\infty,1]$ \index{$\alpha$} is the unique  function on $V_{\infty}$
that is continuous on segments, and satisfies
$$\alpha(v_E)=\frac1{b_E^2}(\check{E}\cdot \check{E})$$ where $E$ is any irreducible component of $X\setminus \mathbb{A}^2_k$ of any compactification $X$ of $\mathbb{A}^2_k$ and $\check{E}$ is the dual divisor of $E$ as defined above.

The skewness function is strictly decreasing, and upper semicontinuous. Therefor it  induces a metric $d_{V_\infty}$ on $V_\infty$ by setting
$$d_{V_\infty}(v_1,v_2):=2\alpha(v_1\wedge v_2)-\alpha(v_1)-\alpha(v_2)$$
for all $v_1, v_2\in V_\infty.$
In particular, any segment in $V_\infty$ carries a canonical metric for which it becomes isometric to a real segment.

\medskip

In an analogous way, one defines the {\em thinness} function\index{thinness function} $A:V_{\infty}\rightarrow [-2,\infty]$\index{$A:V_{\infty}\rightarrow [-2,\infty]$} as the unique, increasing, lower semicontinuous function on $V_{\infty}$ such that for any irreducible exceptional divisor $E$ in some  compactification $X\in \mathcal{C}$, we have $$A(v_E)=\frac1{b_E}\, (1 + \ord_E(dx\wedge dy))~.$$ Here we extend the differential form $dx\wedge dy$ to a rational differential form on $X$.

These parameterizations behave in the following way:
\begin{points}
\item when $v$ is a divisorial valuation, then  $\alpha(v)$ and $A(v)$ are rational numbers;
\item when $v$ is an irrational valuation, then $\alpha(v), A(v)\in \mathbb{R}\setminus \mathbb{Q}$;
\item when $v$ is a curve valuation, then  $\alpha(v)=-\infty$, and $A(v)=+\infty$;
\item when $v$ is an infinitely singular valuation, then $\alpha(v)$ and $A(v)$ can be either finite or infinite.
\end{points}

\section{Potential theory on $V_{\infty}$}\label{sectionpotential}
We refer to \cite[Section 3]{Xieb} for details.

\subsection{Subharmonic functions on $V_{\infty}$}

To any $v\in V_\infty$ we attach its Green function\index{$g_v$}
$$
g_v(w) := \alpha(v \wedge w)~.
$$
This is a decreasing continuous function taking values in $[-\infty, 1]$, satisfying
$g_v(-\deg) =1$.

Given any positive Radon measure $\rho$ on $V_\infty$ we define
$$
g_\rho (w) := \int_{V_\infty} g_v(w) \, d\rho(v)~.
$$
Observe that $g_v(w)$ is always well-defined as an element in $[-\infty, 1]$ since
$g_v \le 1$ for all $v$.
%
Then we recall the following
\begin{thm}[\cite{Xieb}]\label{thmrhotogrhoinj}
The map $\rho \mapsto g_\rho$ is injective.
\end{thm}

One can thus make the following definition.
\begin{defi}
A function $\phi : V_\infty \to \R \cup \{-\infty\}$ is said to be \emph{subharmonic}\index{subharmonic} if there exists a positive Radon measure $\rho$ such that $\phi = g_\rho$. In this case, we write
$\rho = \Delta \phi$ \index{$\Delta \phi$}and call it the \emph{Laplacian}\index{Laplacian} of $\phi$.
\end{defi}
Denote by $\SH(V_{\infty})$\index{$\SH(V_{\infty})$} (resp. $\SH^{+}(V_{\infty})$)\index{$\SH^{+}(V_{\infty})$} the space of subharmonic functions on $V_\infty$ (resp. of non-negative subharmonic functions on $V_\infty$).

\smallskip

The next result collects some properties of subharmonic functions.

\begin{thm}[\cite{Xieb}]
Pick any subharmonic function $\phi$ on $V_\infty$. Then
\begin{points}
\item $\phi$ is decreasing and $\phi(-\deg)=\Delta\phi(V_{\infty})\ge0$;
\item $\phi$ is upper semicontinuous;
\item  for any valuation $v\in V_\infty$ the function
$t \mapsto \phi(v_t)$ is convex, where $v_t$ is the unique valuation in $[-\deg , v]$ of skewness $t$.
\end{points}
\end{thm}

\subsection{Subharmonic functions on finite trees}

Let $T$ be any finite subtree of $V_\infty$ containing $-\deg$. Denote by $r_T: V_\infty \to T$ \index{$r_T$}the canonical retraction defined by sending $v$ to the unique valuation $r_T(v)\in T$ such that
$[r_T(v), v] \cap T = \{ r_T(v)\}$.

For any function $\phi$, set $R_T \phi := \phi \circ r_T$\index{$R_T \phi$}. Observe that
$R_T\phi |_T = \phi |_T$ and that $R_T\phi$ is locally constant outside $T$.

Moreover we have the following
\begin{pro}
For any subharmonic function $\phi$, and any finite subtree $T$ containing $-\deg$, the function $R_T\phi$ is subharmonic. Moreover we have $R_T\phi \ge \phi$ and $\Delta (R_T\phi) = (r_T)_* \Delta \phi$.
\end{pro}
\subsection{Examples of subharmonic functions}
We refer to \cite{Xieb} for detail.
For any nonconstant polynomial $Q\in k[x,y]$, we define the function
$$\log |Q| (v) :=  -v(Q) \in  [-\infty,\infty)~.$$
\index{$\log|Q|(v)$}

\begin{pro}\label{prologqbranch}
The function $\log|Q|$ is subharmonic, and
$$
\Delta (\log |Q|) = \sum_i m_i \delta_{v_{s_i}}
$$
where $s_i$ are the branches of the curve $\{Q=0\}$ at infinity,
and $m_i$ is the intersection number of $s_i$ with the line at infinity in $\P^2_k$.
\end{pro}


\begin{pro}\label{proshpos}
The function $\log^+|Q| := \max \{ 0, \log|Q|\}$\index{$\log^+|Q|$} belongs to $\SH^+(V_{\infty})$.

Denote by  $s_1,\cdots,s_l$ the branches of $\{Q=0\}$ at infinity and by $T$ the convex hull of  $\{ - \deg, v_{s_1},\cdots, v_{s_l}\}$.
Then the support of $\Delta (\log^+|Q|)$ is the set of points $v\in T$ satisfying $v(Q)=0$ and $w(Q)<0$ for all $w\in (v,-\deg]$.
\end{pro}
In particular,  the support of $\Delta (\log^+|Q|)$ is finite.
\subsection{The Dirichlet pairing}
Let $\phi, \psi$ be any two subharmonic functions on $V_\infty$.
Since $\alpha$ is bounded from above one can define the Dirichlet pairing\index{Dirichlet pairing}\index{$\langle\phi,\psi\rangle$}
$$\langle\phi,\psi\rangle:=\int_{V_{\infty}^2}\alpha(v\wedge w)\, \Delta\phi(v)\Delta\psi(w)\in [-\infty,+\infty).$$
\begin{pro}The Dirichlet pairing induces a symmetric bilinear form on $\SH(V_{\infty})$ that satisfies
\begin{equation}\label{eq:dir-pair}
\langle \phi, \psi \rangle= \int_{V_\infty} \phi \,\Delta \psi
\end{equation}
\end{pro}

For any subharmonic function $\phi$ on $V_\infty$, we call $\langle\phi,\phi\rangle$ the \emph{energy}\index{energy} of $\phi$.

\medskip

Next, we recall the following useful result.
\begin{pro}[\cite{Xieb}]\label{prortgephipsi}\label{prophiphi}
Pick any two subharmonic functions $\phi, \psi \in \SH(V_\infty)$. For any finite subtree $T\subset V_\infty$ one has
$$
\langle R_T\phi, R_T\psi \rangle \ge \langle \phi, \psi \rangle~
.$$

In particular, we have $$\langle R_T\phi, R_T\phi \rangle \geq \langle \phi, \phi \rangle~$$  and  the equality holds if and only if $\Delta\phi$ is supported on $T$.
\end{pro}

Finally, we recall
a technical result that will play an important role in the rest of this paper.

%
%

For any set $S\subset V_\infty$ we define $B(S) := \cup_{v\in S} \{ w, \, w\ge v\}$\index{$B(S)$}.

\begin{pro}[\cite{Xieb}]\label{proapoxishplus}Let $\phi$ be a function in $\SH^+(V_{\infty})$ such that $\langle\phi,\phi\rangle=0$ and  the support of the positive measure $\Delta\phi$ is finite, equal to $\{v_1,\cdots,v_s\}$ for some positive integer $s$.

Then for any finite set $S\subseteq B(\{v_1,\cdots,v_s\})$ satisfying $\{v_1,\cdots,v_s\}\not\subseteq S$, there exists a function $\psi\in \SH^+(V_{\infty})$ such that
\begin{itemize}
\item $\psi(v)=0$ for all $v\in B(S)$;
\item
 $\langle\psi,\psi\rangle>0.$
\end{itemize}
\end{pro}
\rem Let $Q$ be a nonconstant polynomial in $k[x,y]$ and set $\phi:=\log^+|Q|$, then $\phi\in \SH^+(V_{\infty})$, $\langle\phi,\phi\rangle=0$ and  the support of the positive measure $\Delta\phi$ is finite.
\endrem

%
%
%
%
%
%
%
%
\subsection{The class of $\mathbb{L}^2$ functions\index{$\mathbb{L}^2$ function}}See \cite[Section 3.7]{Xieb}

We define $\mathbb{L}^2(V_{\infty})$\index{$\mathbb{L}^2(V_{\infty})$} to be the set of functions $$\phi:\{v\in V_{\infty}|\,\,\alpha(v)>-\infty\}\to \mathbb{R}$$ such that
there exist $\phi_1, \phi_2\in \SH (V_\infty)$ with $\langle\phi_1,\phi_1\rangle>-\infty$, $\langle\phi_2,\phi_2\rangle>-\infty$, and
$\phi(v)=\phi_1(v)-\phi_2(v)$
for all valuations with $\alpha(v)>-\infty$.

Observe that $\mathbb{L}^2(V_{\infty})$ is an infinite dimensional vector space
\begin{pro}
The restriction map
$ g \mapsto g|_{\{\alpha > -\infty\}}$ is injective from $\SH(V_\infty) \cap \{ \langle \phi, \phi \rangle > - \infty\}$ into $\mathbb{L}^2(V_{\infty})$.
\end{pro}
We shall thus always identify a subharmonic function with finite energy with its image in $\mathbb{L}^2(V_{\infty})$ so that we have in particular the inclusion $\SH^+(V_{\infty})\subseteq \mathbb{L}^2(V_{\infty})$.

\smallskip

It follows from the Hodge index theorem, see \cite[Theorem 3.18]{Xieb} that
\begin{pro}
For any two subharmonic functions $\phi_1, \phi_2$  with finite energy,
we have $\langle \phi_1, \phi_2 \rangle > -\infty$.
\end{pro}

This result allows one to extend the Dirichlet pairing to the space $\mathbb{L}^2(V_{\infty})$ as a symmetric bilinear form.

%
%
%
%
%
%

%
%



\subsection{Polynomials taking nonnegative values on valuations}

The results in this section are proven in \cite{Xieb}. They will play a crucial role in the proof of our main result.

\smallskip

Given any   finite subset $S$ of  $V_{\infty}$, we define the $k$-algebra $$R_S:=\{P\in k[x,y]|\,\,v(P)\geq 0 \text{ for all }v\in S \}~.$$\index{$R_S$}
When the transcendence degree of the fraction field of $R_S$ over $k$ is equal to $2$, then we say that $S$ is {\em rich}\index{rich}.

One of the main result of \cite{Xieb} is a characterization of rich subsets of $V_\infty$ in terms of the existence of suitable subharmonic functions. To state this result we need to introduce some more notation.

We set
\begin{itemize}
\item
$S^{\min}:=\{v\in S|\,\,v \text{ is minimal in }S\}$\index{$S^{\min}$};
\item
$B(S)^{\circ}$ to be the interior of $B(S)$\index{$B(S)^{\circ}$}.
\end{itemize}
It is easy to check that $R_{S'}\subseteq R_{S}$ if $S\subseteq B(S')$ and then we have $R_S=R_{S^{\min}}$.


The following result is \cite[Theorem 4.7]{Xieb}.
\begin{thm}\label{prodimtwothengeqz}\label{thmsuvsectionprorefshr}\label{thmsuvsectionprorefshr}Let $S$ be a finite set of valuations in $V_{\infty}$. Then the following statements are equivalent.
\points
\item The subset $S$ is rich.
\item There exists a nonzero polynomial $P\in R_{S}$ such that $v(P)>0$ for all $v\in S.$
\item For every valuation $v\in S^{\min}$, there exists a nonzero polynomial $P\in R_{S}$ such that $v(P)>0.$
\item There exists a function $\phi\in \SH^+(V_{\infty})$ such that $\phi(v)=0$ for all $v\in B(S)$ and $\langle\phi,\phi\rangle>0.$
\item There exists a function $\phi\in \mathbb{L}^2(V_{\infty})$, satisfying $\phi(v)=0$ for all $v\in B(S)$ and $\langle\phi,\phi\rangle>0$.
\item There exist a finite set $S'\subseteq V_{\infty}$ satisfying $S\subseteq B(S')^{\circ}$ and $S'$ is rich.
\endpoints
\end{thm}

\begin{rem}In (v) of Theorem \ref{prodimtwothengeqz}, for any $v\in V_{\infty}$ satisfying $\alpha(v) = - \infty$, we say $\phi(v)=0$ if $0\in [\liminf_{w<v,w\to v}\phi(w),\limsup_{w<v,w\to v}\phi(w)]$.
\end{rem}

The next result is \cite[Theorem 4.12]{Xieb}.
\begin{thm}\label{lemgeqm}
Let $S$ be a finite set of valuations in $V_\infty$.
Suppose that there exists a function $\phi \in \SH(V_\infty)$ such that
$\langle\phi,\phi\rangle>0$ and $\phi(v)=0$ for all $v\in B(S)$.

For any integer $l \ge 0$, there exists a real number $M_l \le 1$ such that for any set
$S¡ä$ of valuations such that
\begin{points}\item[(1)]$S'\setminus B(S)$ has at most $l$ elements;
\item[(2)] $S'\subseteq B(S)\cup \{v\in V_{\infty}|\,\, \alpha(v)\leq M_l\}$;
\end{points}
then there exists a function $\phi'\in \mathbb{L}^2(V_{\infty})$
satisfying $\phi'(v)=0$ for all $v\in B(S')$ and $\langle\phi',\phi'\rangle>0.$
\end{thm}

\section{The Riemann-Zariski space at infinity}\label{sectionriezari}

\subsection{Weil and Cartier classes}
See \cite[Appendix A]{Favre2011} or \cite{boucksomfavrejonsson,Cantat2007,Manin1986} for details.


Formally, the Riemann-Zariski space\index{Riemann-Zariski space} of $\mathbb{P}^2_k$ at infinity is
defined as $\mathfrak{X}:=\lim\limits_{\overleftarrow{X\in \mathcal{C}}}X.$ In our paper,
we concern ourself with its classes rather than itself.

For each compactification $X\in \mathcal{C}$, we denote by $N^1(X)_{\mathbb{R}}$ the
$\R$- linear space consisting of $\R$-divisors supported on $X\setminus \A^2_{k}$. For any morphism $\pi:X'\rightarrow X$ between compactifications, we have the pushforward $\pi_*:N^1(X')_{\mathbb{R}}\rightarrow N^1(X)_{\mathbb{R}}$ and the pullback $\pi^*:N^1(X)_{\mathbb{R}}\rightarrow N^1(X')_{\mathbb{R}}$, see \cite{debarre,Lazarsfeld} for details.

The space of \index{Weil classes} Weil classes of $\mathfrak{X}$ is defined to be the projective limit \index{$W(\mathfrak{X})$}$$W(\mathfrak{X}):=\lim\limits_{\overleftarrow{X\in \mathcal{C}}}N^1(X)_{\mathbb{R}}$$ with respect to pushforward arrows.
Concretely, a Weil class $\alpha\in W(\mathfrak{X})$ is given by its $incarnations$ $\alpha_{X}\in N^1(X)_{\mathbb{R}}$, compatible with pushforwards; that is, $\pi_*\alpha_{X}=\alpha_{X'}$ as soon as $\pi:X\rightarrow X'$. Observe that we may define a Weil class by its incarnations.

If $\alpha_X\in N^1(X)_{\mathbb{R}}$ is a class in some compactification $X\in \mathcal{C}$, then $\alpha_X$ defines a
Weil class $\alpha$, whose incarnation $\alpha_{X'}=\mu_*\pi_*\alpha_X$ where $\pi:X_1\rightarrow X$ and $\mu:X\rightarrow X'$ are morphisms between compactifications. We say that $\alpha$ is determined in $X$. A {\em Cartier class} \index{Cartier class} is a Weil class determined in a certain compactification. Denote by \index{$C(\mathfrak{X})$}$C(\mathfrak{X})$ the space of Cartier classes.




For each $X$, the intersection pairing $N^1(X)_{\mathbb{R}}\times N^1(X)_{\mathbb{R}}\rightarrow \mathbb{R}$ is denoted by $(\alpha\cdot \beta )_{X}$. By the pull-back formula, it induces a perfect pairing $W(\mathfrak{X})\times C(\mathfrak{X})\rightarrow \mathbb{R}$ which is denoted by $( \alpha \cdot \beta ).$ It induce an inner product on $C(\mathfrak{X})$.
The space \index{$\mathbb{L}^2(\mathfrak{X})$}$$\mathbb{L}^2(\mathfrak{X}):=\{\alpha\in W(\mathfrak{X})|\,\inf_X{(\alpha_X\cdot \alpha_X)}>-\infty\}$$ is the completion of $C(\mathfrak{X})$ under inner product. It is an infinite dimensional subspace of $W(\mathfrak{X})$ that contains $C(\mathfrak{X})$. It is endowed with a natural intersection product extending the one on Cartier classes and that is of Minkowski's type, see \cite{Cantat2007} or \cite{favreBourbaki}.
\subsection{Nef $b$-divisors and subharmonic functions}\label{subsectionclassesandva}In this section, we summarize the relations between the classes of  the Riemann-Zariski space at infinity and the potential theory of $V_{\infty}$.

Let $\mathcal{E}$ be the set of all irreducible component of $X\setminus \mathbb{A}^2_k$ for all compactifications $X$ of $\mathbb{A}^2_k$, modulo the following equivalence relations:
two divisors $E$, $E'$ in $(X,\iota)$ and $(X',\iota')$ are equivalent
if there exists a birational morphism $\pi:X\dashrightarrow X'$ such that $\pi\circ \iota_1=\iota_2$ sends $E$ onto $E'$.
As in \cite[Section 1.3]{boucksomfavrejonsson}, we may identify $W(\mathfrak{X})$ to $\mathbb{R}^{\mathcal{E}}$ and $C(\mathfrak{X})$ to $\oplus_{\mathcal{E}}\mathbb{R}$. The pairing is given by $(\alpha\cdot \beta)=\sum_{E\in \mathcal{E}}c_Ed_E$ where $\alpha=(c_E)_{E\in \mathcal{E}}$ and $\beta=\oplus_{E\in \mathcal{E}} d_E$ are Weil and Cartier divisors respectively. We first describe these identifications.

Given a compactification $X\in \mathcal{C}$ and let $E_1,\cdots, E_m$ be all irreducible exceptional divisors of $X$, the incarnation of $\alpha=(c_E)_{E\in \mathcal{E}}$ is $\alpha_X=\sum_{i=1,\cdots,m}c_{E_i}E_i$.

For any $E\in \mathcal{E}$, pick a compactification $X\in \mathcal{C}$ such that $E$ is an irreducible component of $X\setminus \mathbb{A}^2_k$. We denote by $\check{E}$ the unique class in $N^1(X)_{\mathbb{R}}\subseteq C(\mathfrak{X})$ satisfying $(\check{E}\cdot F)_X=0$ when $F$ is an irreducible component different from $E$ and $(\check{E}\cdot E)_X=1$. As a Cartier class, $\check{E}$ does not depend on the choice of the compactification $X$. The identification of $\oplus_{\mathcal{E}}\mathbb{R}$ to $C(\mathfrak{X})$ is given by $\oplus_{E\in \mathcal{E}} d_E\mapsto \sum_{E\in \mathcal{E}}d_E\check{E}$.

We define a map $i_C:C(\mathfrak{X})\to \mathcal{C}^0(V_{\infty},\mathbb{R})$ where $\mathcal{C}^0(V_{\infty},\mathbb{R})$ is the set of continuous functions on $V_{\infty}$ by $\sum_{E\in \mathcal{E}}d_E\check{E}\mapsto \sum_{E\in \mathcal{E}}b_Ed_Eg_{v_E}.$ Observe that $i_C$ is an embedding.


\medskip

We denote by $\Nef_{\infty}(\mathfrak{X})$\index{$\Nef_{\infty}(\mathfrak{X})$} the set of all Weil classes $\alpha\in W(\mathfrak{X})$ such that for any compactification $X\in \mathcal{C}$, the incarnation $\alpha_X$ is nef at infinity i.e. for any irreducible component $E$ of $X\setminus \mathbb{A}^2$, we have $(\alpha_X\cdot E)\geq 0.$

Let $g$ be a continuous function on $V_{\infty}$, by \cite[Lemma 3.5]{Xieb}, we can prove that there exists a sequence Cartier classes $\beta_n\in C(\mathfrak{X})$ satisfying $i_C(\beta_n)\to g$ uniformly as $n\rightarrow\infty$.

\begin{lem}\label{lembntogabn}The limit $\lim_{n\to\infty}(\beta_n\cdot \alpha)$ exists and does not depend on the choice of the sequence $\beta_n.$
\end{lem}
\proof
We only have the show that given a real number $\varepsilon>0$, for any  Cartier class $\beta$ in $C(\mathfrak{X})$ satisfying $|i_C(\beta)|\leq \varepsilon$ on $V_{\infty},$ there exists a constant $C>0$ such that $|(\alpha_X\cdot\beta_X)|\leq C\varepsilon.$

There exists an admissible compactifiction $X\in \mathcal{C}$ such that $\beta$ is determined in a $X$ and then $(\alpha\cdot\beta)=(\alpha_X\cdot\beta_X).$
Observe that $\beta_X=\sum b_Ei_C(\beta)E$ where the sum is over all irreducible components of $X\setminus \mathbb{A}^2.$
Denote by $\pi:X\rightarrow \mathbb{P}^2$ the dominant morphism between compactifications and $L_{\infty}$ the line at infinity of $\mathbb{P}^2.$
Observe that $\pi^*L_{\infty}=\sum b_EE$ and then $\varepsilon\pi^*L_{\infty}\pm \beta_X$ are effective. It follows that
$$|(\alpha\cdot\beta)=(\alpha_X\cdot\beta_X)|\leq \varepsilon(\alpha_X\cdot\pi^*L_{\infty})=\varepsilon(\alpha_{\mathbb{P}^2}\cdot L_{\infty}).$$
\endproof

The same argument in the proof of Lemma \ref{lembntogabn} also shows that the map $g\mapsto \lim_{n\to\infty}(\beta_n\cdot \alpha)$ is continuous on $\mathcal{C}^0(V_{\infty},\mathbb{R})$. This map defines a real Radon measure $\rho_{\alpha}.$ Observe that if $\beta$ is effective, $(\alpha,\beta)\geq 0$.
By \cite[Lemma 3.5]{Xieb}, we can prove that $\int_{V_{\infty}}fd\rho_{\alpha}\geq 0$ when $f$ is nonnegative. Then we get
\begin{lem}\label{lemrhoalpisposi}The real Radon measure $\rho_{\alpha}$ is positive.
\end{lem}

We define a map $i_N:\Nef_{\infty}(\mathfrak{X})\rightarrow \SH(V_{\infty})$ by $\alpha\mapsto g_{\rho_{\alpha}}$ and we have
\begin{pro}The map $i_N$ is bijective.
\end{pro}
\proof We define a map $j_N:\SH(V_{\infty})\rightarrow W(\mathfrak{X})$ be $\phi\mapsto (b_E\phi(v_E))_{E\in \mathcal{E}}.$ We only have to show that $j_N$
is the inverse of $i_N.$

We first claim that $j_N(\SH(V_{\infty}))\subseteq \Nef_{\infty}(\mathfrak{X}).$ Set $\alpha:=j_N(\phi)=(b_E\phi(v_E))_{E\in \mathcal{E}}\in W(\mathfrak{X}).$
For any compactifiction $X\in \mathcal{C}$, denote by $E_1,\cdots,E_m$ all the irreducible components of $X\setminus \mathbb{A}^2.$
We have $\alpha_X=\sum_{i=1}^m b_{E_i}\phi(v_{E_i})E_i$. Observe that $i_C(\alpha_X)$ is the unique function on $V_{\infty}$ satisfying
\begin{points}
\item $i_C(\alpha_X)(v_{E_i})=\phi(v_{E_i})$ for all $i=1,\cdots,m;$
\item $i_C(\alpha_X)$ is linear outside $\{v_{E_1},\cdots,v_{E_m}\}.$
\end{points}
It follows that $i_C(\alpha_X)$ takes form $\sum_{i=1}^ma_ig_{v_{E_i}}$ where $a_i\geq 0$ for all $i=1,\cdots,m.$ Then we have $\alpha_X=\sum a_ib_{E_i}^{-1}\check{E}_i$ and thus it is nef at infinity. It follows that $\alpha\in \Nef_{\infty}(V_{\infty}).$

For any $\phi\in \SH(V_{\infty})$, we have $j_N(\phi)=(b_E\phi(v_E))_{E\in \mathcal{E}}.$ By the definition of $i_N$, for all divisorial valuation $v_E\in V_{\infty}$,  we have
$$i_N(j_N(\phi))(v_E)=\int_{V_{\infty}}g_{v_E}d\rho_{j_N(\phi)}=(j_N(\phi)\cdot b_E^{-1}\check{E})=\phi(v_E).$$
It follows that $i_N\circ j_N=\id$ and then we conclude our proof.
\endproof

At last, we define an embedding $j_{\mathbb{L}^2}:\mathbb{L}^2(V_{\infty})\rightarrow \mathbb{L}^2(\mathfrak{X})$ as follows:
Let $\phi$ be any function in $\mathbb{L}^2(V_{\infty})$. Write $\phi=\phi_1-\phi_2$ on $\{v\in V_{\infty}|\,\,\alpha(v)>-\infty\}$ where $\phi_1,\phi_2$ are functions in $\SH(V_{\infty})$ satisfying $\langle\phi_i,\phi_i\rangle>-\infty$ for $i=1,2.$ Then $j_{\mathbb{L}^2}(\phi)$ is defined to be $i_N^{-1}(\phi_1)-i_N^{-1}(\phi_2).$
This definition does not depend on the choice of $\phi_1,\phi_2$ and satisfying $\langle\phi,\psi\rangle=\langle j_{\mathbb{L}^2}(\phi),j_{\mathbb{L}^2}(\psi)\rangle$ for all $\phi,\psi\in\mathbb{L}^2(V_{\infty}).$

\smallskip

For any $v\in V_{\infty}$, set $Z_v:=i_N^{-1}(g_v)$ the Weil class in  $\Nef_{\infty}(\mathfrak{X}).$ Observe that $Z_v\in \mathbb{L}^2(\mathfrak{X})$ when $\alpha(v)>-\infty$ and $Z_v\in C(\mathfrak{X})$ when $v$ is divisorial. If $v=v_E$ where $E$ is an irreducible component of $X\setminus \mathbb{A}^2_k$ for compactification $X\in \mathcal{C}.$ Denote by $\check{E}$ the duality of $E$ in $N^1(X)$ w.r.t. the intersection pairing. View $\check{E}$ as a Cartier class of $\mathfrak{X}$, then
we have $Z_{v_E}=b_E^{-1}\check{E}.$
Finally we recall the following
\begin{pro}\label{prozvzwealp}\cite[Lemma A.2]{Favre2011}For any two valuations $v,w\in V_{\infty}$ one has $(Z_v\cdot Z_w)=\alpha(v\wedge w).$
\end{pro}

\section{Background on dynamics of polynomial maps}\label{sectionbgdopm}

In this section we assume that $k$ is an algebraically closed field of characteristic zero.
Recall that the affine coordinates have been fixed, $\A^2_k = \Spec k[x,y]$.

\subsection{Dynamical invariants of polynomial mappings}
The (algebraic) degree of  a  dominant polynomial endomorphism   $f=(F(x,y),G(x,y))$ defined on $\mathbb{A}^2_k$ is  defined by
$$\deg(f):=\max\{\deg(F),\deg(G)\}~.$$\index{$\deg(f)$}
It is not difficult to show that the sequence $\deg(f^n)$ is sub-multiplicative, so that the limit
$\la_1(f):=\lim_{n\rightarrow \infty}(\deg(f^n))^{\frac{1}{n}}$\index{$\la_1(f)$} exists.
It is referred to as the {\em dynamical degree} of $f$, and it is a Theorem of Favre and Jonsson that $\la_1(f)$ is always a quadratic integer, see \cite{Favre2011}.

\smallskip

The (topological) degree  $\la_2(f)$\index{$\la_2(f)$} of $f$ is defined to be the number of preimages of a general closed point in $\A^2(k)$; one has $\la_2(fg)=\la_2(f)\la_2(g).$

It follows from B\'ezout's theorem that $\la_2(f) \le \deg(f)^2$ hence
\begin{equation}\label{eq:ineq}
\la_1(f)^2\geq \la_2(f)~.
\end{equation}

The resonant case $\la_1(f)^2 = \la_2(f)$ is quite special and the following structure theorem for these maps is proven in \cite{Favre2011}.

\begin{thm}
Any polynomial endomorphism $f$ of $\A^2_k$ such that $\la_1(f)^2 = \la_2(f)$ is proper\footnote{We say a polynomial endomorphism $f$ of $\A^2_k$ is proper if it is a proper morphism between schemes. When $k=\C$, it means that the preimage of any compact set of $\C^2$ is compact.}, and we are in one of the following two exclusive cases.
\begin{enumerate}
\item
$\deg(f^n) \asymp \la_1(f)^n$: there exists a compactification $X$ of $\A^2_k$ to which $f$ extends as a regular map $f : X \to X$.
\item
$\deg(f^n) \asymp n\la_1(f)^n$: there exist affine coordinates $x,y$
in which $f$ takes the form
$$f(x,y) = (x^l + a_1 x^{l-1}+ \ldots + a_l , A_0(x) y^l + \ldots + A_l(x))
$$
where $a_i \in k$ and $A_i \in k[x]$ with $\deg A_0 \ge1$, and $l = \la_1(f)$.
 \end{enumerate}
\end{thm}
\rem Regular endomorphisms as in (i) have been classified in \cite{Favre2011}.
\endrem
\subsection{Valuative dynamics}\label{subsectionvaldy}
Any dominant polynomial endomorphism $f$ as in the previous section induces a natural map on the space of valuations at infinity in the following way.

\smallskip

For any $v\in V_\infty$ we may set\index{$d(f,v)$}
$$d(f,v):=-\min\{v(F),v(G),0\}\geq 0~.$$
In this way, we get a non-negative continuous  decreasing function on $V_\infty$
such that  $d(f,v)\geq \deg(f)\alpha(v)$. Observe also that $d(f,-\deg)=\deg(f)$. It is a fact that $f$ is proper if and only if $d(f,v)>0$ for all $v\in V_{\infty}.$

\smallskip

We now set\index{$f_* v$}
\begin{itemize}
\item
$f_* v :=0$ if $d(f,v) = 0$;
\item
$f_*v(P) =v(f^*P)$ if $d(f,v)> 0$.
\end{itemize}
In this way one obtains a valuation on $k[x,y]$ (that may be trivial); and
we then get a continuous map\index{$f_{\d}$}
$$f_{\d}:\{v\in V_{\infty}|\,\,d(f,v)>0\}\rightarrow V_{\infty}$$ by $$f_{\d}(v):=d(f,v)^{-1}f_*v ~.$$

For any subset $S$ of $V_{\infty},$ set $f_{\d}^{-1}(S):=\{v\in V_{\infty}|\,\, d(f,v)>0 \text{ and }f_{\d}(v)\in S\}.$ If $f$ is an open set, then $f_{\d}^{-1}(S)$ is also open.

This map $f_{\d}$ can extend to a continuous map $f_{\d}:\overline{\{v\in V_{\infty}|\,\,d(f,v)>0\}}\rightarrow V_{\infty}$. The image of any $v\in \partial{\{v\in V_{\infty}|\,\,d(f,v)>0\}}$ is a curve valuation defined by a rational curve with one place at infinity.

\begin{lem}\label{lembrafact}
Let $C$, $D$ be two branches of curves at infinity satisfying $f(C)=D$. Then we have
$m_Cd(f,v_{C})=\deg(f|_{C})m_D$ where
$m_C=(C\cdot l_{\infty})$ and $m_D=(D\cdot l_{\infty})$.
\end{lem}
\proof
 Let $L$ be a general linear form in $k[x,y]$, we have
$$m_Cv_{C}(f^*L)=\deg(f|_{C})m_Dv_{D}(L)=\deg(f|_{C})m_D.$$ On the other hand, we have $v_{C}(f^*L)=d(f,v_{C})$. It follows that $$m_Cd(f,v_{C})=\deg(f|_{C})m_D.$$
\endproof

We now recall the following key result, \cite[Proposition 2.3\,,Theorem 2.4,\,Proposition 5.3.]{Favre2011}.
\begin{thm}\label{thm:eigenval}
There exists a unique valuation $v_*$
such that $\alpha(v_*)\geq 0\geq A(v_*)$, and
$f_*v_*=\la_1 v_*$.

If $\la_1(f)^2>\la_2(f)$, this valuation is unique .

If $\la_1(f)^2=\la_2(f)$, the set of such valuations is a closed segment.
\end{thm}

This valuation $v_*$ is called the {\em eigenvaluation} \index{eigenvaluation}of $f$ when $\la_1(f)^2>\la_2(f)$.

%
%
%

%
%
%
%
%

\subsection{Functoriality of classes of the Riemann-Zariski space}\cite[Appendix A]{Favre2011}
Let $f$ be a dominant polynomial endomorphism on $\mathbb{A}^2$ defined over $k$.

we have natural actions $f^*:C(\mathfrak{X})\rightarrow C(\mathfrak{X})$ induced by the pullback between the N\'eron-Severi groups and $f_*:W(\mathfrak{X})\rightarrow W(\mathfrak{X})$ induced by the pushforward between the N\'eron-Severi groups. Further, we have the projection formula
$$(f_*\beta\cdot\gamma)=(\beta\cdot f_*\gamma)$$ for $\beta\in C(\mathfrak{X})$ and $\gamma\in W(\mathfrak{X}).$

The pushforward (resp. pullback) preserves (resp. extends to) $\mathbb{L}^2$
-classes.
We obtain bounded operators $f^*,f_*:\mathbb{L}^2(\mathfrak{X})\rightarrow \mathbb{L}^2(\mathfrak{X})$
and $(f_*\beta\cdot \gamma)=(\beta\cdot f^*\gamma)$ for $\beta,\gamma\in \mathbb{L}^2(\mathfrak{X})$. We have $f_*f^*=\la_2(f)$ on $\mathbb{L}^2(\mathfrak{X})$.

\begin{lem}\label{lemfajcas}\cite[Lemma A.6]{Favre2011}We have $f_*Z_v = d(f,v)Z_{f_{\d}(v)}$ for all $v\in \widehat{V}_{\infty}.$
\end{lem}

\newpage

\part{Local arguments}
In this part we collect some arguments of local nature that will play an
important role in the proof of our main result.  We first recall the
definition of the local valuation space as in \cite{Favre2007}.
Then we state and prove a local version of the dynamical Mordell-Lang
conjecture for superattracting analytic germs (Theorem \ref{thmlodml}).

\section{The local valuative tree and the local Riemann-Zariski space}
Let $(X,q)$ be a smooth germ of surface at a closed point $q$ defined over an algebraically closed field $k.$
Pick a local coordinate $(z,w)$ at $q$ and set $\mathfrak{m}:=(z,w)$.

\subsection{The local valuative tree\index{local valuative tree}}See \cite{Jonsson} for details.
We first introduce the local avatar of the valuative tree at infinity defined in  \cite{Favre2004}.

We define the space $V_q$ \index{$V_q$}of  valuations that are trivial on $k^*$
and centered at $q$, and normalized by the condition
$$v(\mathfrak{m}) =\min\{v(z),v(w)\}=1.$$
The order of vanishing $\ord_q$ at the point $q$ is a valuation in $V_q$.


%

The space $V_{q}$  is
equipped with a partial ordering defined by $v\leq w$ if and only if $v(f)\leq w(f)$ for all $f\in k[[z,w]]$ for which is again a real tree (see \cite{Favre2004,Favre2007,Jonsson}).
The valuation $\ord_q$ is the minimal element of $V_{q}$.

Let $\pi:Y\rightarrow X$ be a morphism between compactifications in $\mathcal{C}$ such that $\pi$ is an isomorphism above $X\setminus \{q\}$.
Let $F$ be an irreducible component of $\pi^{-1}(q)$.
Set \index{$b^q_F$}$b^q_F:=\ord_F\pi^*\mathfrak{m}\in \mathbb{N}^+$, then \index{$v_{F}^q$}$v_{F}^q:=b^q_F\ord_F$ is contained in $V_{q}$.
Let $\check{F_q}$\index{$\check{F_q}$} be the unique divisor supported on $\pi^{-1}(q)$ such that $(\check{F_q},F')=\delta_{F,F'}$. The quantity
$(\check{F_q}\cdot \check{F_q})$ is independence on the choice of $Y$.

There exists a unique increasing and lower semicontinuous function\index{$\alpha^q$} $\alpha^q:V_{q}\rightarrow [1,+\infty]$ on $V_{q}$
satisfying $\alpha^q(v^q_F)=-(b^q_F)^{-2}(\check{F_q}\cdot \check{F_q}).$

\medskip

\bigskip


At last we talk about the connection between the local valuative tree and the global one.
Now we suppose that $X$ is a compactification of $\mathbb{A}^2_k$ in $\mathcal{C}$ defined over an algebraically closed field $k$ and
$q$ be a $k$-point in $X\setminus \mathbb{A}^2_k$.
Let $\{E_1,\cdots,E_s\}$ be the set of irreducible exceptional divisors containing $q$. We have $s=1$ or $2$.

\begin{exe}For $i=1,\cdots,s$,
there exists a valuation $v^q_{E_i}$ defined by $P\mapsto \ord_q(P|_{E_{i}})$ for $P\in k[[z,w]].$
\end{exe}

Denote by $U(q)$ \index{$U(q)$}that set of
valuations in $V_{\infty}$ whose centres in $X$ are $q$ and set $\overline{U(q)}:=U(q)\cup \{v_{E_1},\cdots, v_{E_s}\}$. For any $v\in U(q)$, there exists $r_q(v)\in \mathbb{R}^+$ such that $r_q(v)v\in V_q.$ Set $v^q:=r_q(v)v$ when $v\in U_q$ and $v^q:=v_{E_i}^q$ when $v\in \{v_{E_1},\cdots, v_{E_s}\}.$
The map $\overline{U(q)}\rightarrow V_q$ defined by $v\mapsto v^q$ is a homeomorphism.
When $v^q\in V_{q}\setminus \{v_{E_1}^q,\cdots,v_{E_s}^q\}$, the type of $v^q$ is the same as the type of $v$ as a valuation in $\overline{U(q)}$; if $v^q=v_{E_i}^q$, $v^q$ is a curve valuation.


\medskip

\subsection{The local Riemann-Zariski space\index{local Riemann-Zariski space}}Analogue to the Riemann-Zariski space at infinity, we can also define the Riemann-Zariski space at a point.

Let $(X,q)$ be a smooth germ of surface at a closed point $q$ defined over an algebraically closed field $k.$
Pick a local coordinate $(z,w)$ at $q$ and set $\mathfrak{m}:=(z,w)$.
We define $\mathcal{C}^q$ be the category of biratonal model $\pi:X_{\pi}\rightarrow X$ such that $\pi$ is an isomorphism above $X\setminus \{q\}.$ We denote by $N^1_q(X_{\pi})_{\mathbb{R}}$ the kernel of $\pi_*:N^1(X_{\pi})_{\mathbb{R}}\rightarrow N^1(X)_{\mathbb{R}}$.

As in Section \ref{subsectionclassesandva},
formally, the Riemann-Zariski space of $X$ at $q$ is
defined as $\mathfrak{X}^q:=\lim\limits_{\overleftarrow{X_{\pi}\in \mathcal{C}^q}}X_{\pi}.$
The space of Weil classes of $\mathfrak{X}^q$ is defined to be the projective limit\index{$W(\mathfrak{X}^q)$} $$W(\mathfrak{X}^q):=\lim\limits_{\overleftarrow{X_{\pi}\in \mathcal{C}^q}}N^1_q(X_{\pi})_{\mathbb{R}}$$ with respect to pushforward arrows. The space of Cartier classes on $\mathfrak{X}^q$ is defined to be the inductive limit\index{$C(\mathfrak{X}^q)$} $$C(\mathfrak{X}^q):=\lim\limits_{\overrightarrow{X_{\pi}\in \mathcal{C}^q}}N^1_q(X_{\pi})_{\mathbb{R}}$$ with respect to pullback arrows. As in Section \ref{subsectionclassesandva}, we embed $C(\mathfrak{X}^q)$ in $W(\mathfrak{X}^q)$.
The intersection pairing in $N^1_q(X_{\pi})$ induced an intersection pairing $W(\mathfrak{X}^q)\times C(\mathfrak{X}^q)\rightarrow \mathbb{R}.$ This pairing is perfect.

We identify $W(\mathfrak{X}^q)$ to $\mathbb{R}^{\mathcal{E}^q}$ and $C(\mathfrak{X}^q)$ to $\oplus_{\mathcal{E}}^q\mathbb{R}$ where $\mathcal{E}^q$ is the set of equivalence classes of irreducible exceptional divisor above $q$.

There exists a continuous embedding $V_q\rightarrow W(\mathfrak{X}^q)$ defined by\index{$Z^q_v$} $v\mapsto Z^q_v:=(b^q_E\alpha^q(v_E^q\wedge v))_{E\in \mathcal{E}^q}.$
If $v$ is divisorial, we see that $Z^q_v$ is a Cartier class.  By continuality, we can define the pairing  $(Z^q_{v_1}\cdot Z^q_{v_2}):=\alpha^q(v_1\wedge v_2)$.

\subsection{Dynamics on the local valuative tree}\label{sectionlocaldynamics}
In this section, we recall some background on dynamics on the local valuative tree.

\smallskip
Let $(X,q)$ be a smooth germ of surface at a closed point $q$ defined over an algebraically closed field $k.$
Pick a local coordinate $(z,w)$ at $q$ and set $\mathfrak{m}:=(z,w)$.
Let $f:(X,q)\rightarrow (X,q)$ be a germ of dominant endomorphism of $(X,q)$.

For any valuation $v\in V_q$, we define a valuation $f_*(v)$ by $f_*(v)(\phi):=v(f^*\phi)=v(\phi\circ f).$
In general, $f^*(v)$ is not normalized and may be identically $+\infty$ in $\mathfrak{m}$. The latter situation appears
exactly when $v=v_C$ is a $contracted$ $curve$ $valuation$ i.e. $C$ is a branch of curve at $q$ contracted by $f$.
Denote by $\mathfrak{C}_f$\index{$\mathfrak{C}_f$} the set of contracted curve valuation.
Observe that $\mathfrak{C}_f$ is finite.
For $v\in V_q$, set\index{$c(f,v)$} $c(f,v):=\min\{v(f^*x),v(f^*y)\}\in [0,\infty]$. Observe that $c(f,v)=0$ if and only if $v\in \mathfrak{C}_f$.

If $v\in V_q$ is not a contracted curve valuation, $f_{\d}(v)$ is defined to be $c(f,v)^{-1}f_*v$ and we have $\min\{f_{\d}(v)(x),f_{\d}(v)(y)\}=1.$ Then we have $f_{\d}(v)\in V_q.$

Set $c(f):=c(f,\ord_q).$ Observe that $c(f,v)$ is increasing in $V_q$ and by \cite[Proposition 7.14]{Jonsson} we have $c(f,v)\leq \alpha^q(v)c(f)$.

See \cite[Theorem 3.1]{Favre2007}, we can extend the map $f_{\d}:V_q\setminus \mathfrak{C}_f\rightarrow V_{q}$ to a unique continue endomorphism $V_{q}\rightarrow V_q$.

Then we recall \cite[Proposition 3.4]{Favre2007} as follows:
\begin{pro}\label{prolocalcfvlocalcon}
The subset $T_f^q$ of $V_{q}$ where $c(f,\d)$ is not locally constant is a finite
closed subtree of $V_{q}$. Its maximal elements are exactly the maximal
elements in the finite set $\mathfrak{E}^q_f$ consisting of divisorial valuations $v$ with $f_{\d}v=\ord_q$ and
of contracted curve valuations
\end{pro}

\medskip

Next, we recall

\begin{defi}\cite[Definition 4.1]{Favre2007}The {\em asymptotic attraction rate}\index{asymptotic attraction rate} of $f$ is  \index{$c_{\infty}(f)$}$c_{\infty}(f):=\lim_{n\rightarrow \infty}c(f^{n})^{\frac{1}{n}}$.
\end{defi}
This limit exists and does not dependent on the choice of coordinate.

Observe that,
if $df(q)=0$, then we have $c_{\infty}(f)>1.$
\begin{defi}We say a valuation $v_*\in V_q$ is an {\em eigenvaluation}\index{eigenvaluation} if it satisfies the following conditions:
\begin{points}
\item $f_{\d}v_*=v_*$;
\item $d(f,v_*)=c_{\infty};$
\item either $v_*$ is divisorial
or there exists an arbitrary small neighborhood $U$ of $v_*$ taking form
$U=\{v, v>v_1\}$ or $U=\{v,v_2>v\wedge v_2>v_1\}$  such that $f_{\d}U\subseteq U.$
\end{points}
\end{defi}

We recall
\begin{thm}\label{thmfavre2007eigencaluation}\cite[Theorem 4.2, Proposition 5.2]{Favre2007}
If $f$ is dominant with $c_{\infty}(f)>1$ at $q$, then there exists an eigenvaluation $v_*$ in $V_{q}.$
\end{thm}

An eigenvaluation $v_*$ is said to be attracting if it has a neighborhood $U$ in $V_q$ such that for any valuation $v\in U$ we have $f^n_{\d}(v)\to v_*$ as $n\to \infty.$

Recall that a fixed point germ is called {\em rigid}\index{rigid} if its critical set is contained in a
totally invariant set with normal crossings.

\begin{thm}\label{thmfavrejnormalform}\cite[Theorem 5.1]{Favre2007}
Suppose that $f$ is dominant with $c_{\infty}(f)>1$ and $v_*$ is an eigenvaluation in $V_q$.  Then one
can find a modification $\pi:(\widetilde{X},p)\rightarrow (X,q)$
such that the lift $\widetilde{f}$ of $f$ is regular at $p$,
$\widetilde{f}(p) = p$ and $\widetilde{f} :(\widetilde{X},p)\rightarrow(\widetilde{X},p)$ is
rigid. Moreover, if $v_*$ is nondivisorial and attracting, we may ask $p$ to be the center of $v_*$ in $\widetilde{X}$ and $d\widetilde{f}(p)^2=0.$
\end{thm}

%

Finally we prove a technical lemma which is useful in the rest of the paper.

\begin{lem}\label{lemlocalallessinftytendtocu}Let $C$ be an irreducible formal curve in $X$ containing $q$ such that $f^*C=dC$ and $f_*C=mC$ locally. If $d>m$, then there exists $w_1<v_C$ arbitrary close to $v^q_C$ such that for any $v\in W:=\{v\in V_1|\,\, v\wedge v^q_C>w_1\}$, we have $f_{\d}^n(v)\rightarrow v^q_C$ as $n\to \infty$ and $f_{\d}(W)\subseteq W.$ Moreover, for any $M>0$, there exists $N\geq 0$ such that $\{v\in V_q|\,\,\alpha^q(v)\leq M\}\subseteq f^{-N}_{\d}(W).$

In particular,
for any $v\in V_{q}$ satisfying $\alpha(v)<+\infty$ we have $f_{\d}^n(v)\rightarrow v^q_C$ as $n\to \infty.$
\end{lem}
\proof
Since $f^*C=dC$, we have $f^*Z^{q}_{v^q_C}=dZ^{q}_{v^q_C}.$
For all $v\in V_{q}$ and $n\geq 0$, we have $$\alpha^{q}(v_C\wedge v)=1/d^n(f^{*n}Z^{q}_{v^q_C}\cdot Z^{q}_v)$$$$=c(f^n,v)/d^n(Z^{q}_{v^q_C}\cdot Z^{q}_{f^n_{\d}(v)})=c(f^n,v)/d^n\alpha^{q}(v^q_C\wedge f^n_{\d}(v)).$$

Let $P$ be any polynomial in $k[x,y]$, we have $v^q_C(f^*P)=mv^q_C(P).$ It follows that $d(f,v^q_C)=m.$
Since the function $c(f,\cdot)$ is continuous on $V_q$ and $v_C$ is a curve valuation, there exists $w_1<v^q_C$ such that $c(f,v)< c(f,v^q_C)+1/2<d$ for all $v> w_1$. Set $W:=\{v\in V_{q}|\,\,v>w_1\}$, it follows that $f_{\d}(W)\subseteq W$ and for any $v\in W$, $f_{\d}^n(v)\rightarrow v^q_C$ as $n\rightarrow \infty$.

For all $n\geq 0$, we have $c(f^n)\leq c(f^n,v^q_C)=m^n$. It follows that $c(f^n,v)\leq c(f)\alpha^q(v)\leq m^n\alpha^q(v)$ for all $n\geq 0$ and all $v\in V_{\infty}.$

For any $M\geq 0$, there exists $N\geq 0$ such that $\frac{d^n}{Mm^n}>\alpha^q(w_1).$ It follows that for all $v\in V_q$ satisfying $\alpha(v)\leq M$, we have
$$\alpha^{q}(v^q_C\wedge f^N_{\d}(v))=d^nc(f^n,v)^{-1}\alpha^{q}(v^q_C\wedge v)\geq d^nc(f^n,v)^{-1}\geq \frac{d^n}{Mm^n}>\alpha^q(w_1).$$ It follows that $f^N_{\d}(v)\in W$ and then $f_{\d}^n(v)\rightarrow v^q_C$ as $n\to \infty.$
\endproof

\subsection{Compute local intersection of curves at infinity}
Let $C_1$, $C_2$ be two formal curves at infinity. The aim of this section is to compute the local intersection of them.

Denote by $l_{\infty}$ the line at infinity in $\mathbb{P}^2_k.$
Set $m_i:=(C_i\cdot l_{\infty})$ for $i=1,2.$  Pick a compactification $\pi:X\rightarrow \mathbb{P}^1_k$ which dominates $\mathbb{P}^2_k$ such that centers $q_i$ of the strict transform $\pi_i^{\#}C_i$'s of $C_i$'s are distinct, each $q_i$ lies in a unique irreducible exceptional divisor $E_i$ and $C_i$ is smooth at $q_i$ for $i=1,2$. We may suppose that $E_i\neq l_{\infty}$ for $i=1,2.$

Write $\pi^*C_i=\pi_i^{\#}C_i+Z_i$ where $Z_i\in N^1(X)_{\mathbb{R}}$ for $i=1,2.$ It follows that we have
$$(Z_i\cdot E_i)=\left((\pi^*C_i-\pi^{\#}C_i)\cdot E_i\right)=(C_i\cdot \pi_*E_i)-(\pi^{\#}C_i\cdot E_i)=-1;$$
$$(Z_i\cdot \pi^{\#}l_{\infty})=\left((\pi^*C_i-\pi^{\#}C_i)\cdot \pi^{\#}l_{\infty}\right)=m_i;$$
and
$$(Z_i\cdot E)=\left((\pi^*C_i-\pi^{\#}C_i)\cdot E\right)=0$$ for irreducible exceptional divisor $E$ different from $E_i$ and $\pi^{\#}l_{\infty}.$
Observe that we have $$m_i=(C_i\cdot l_{\infty})=(\pi^{\#}C_i\cdot \pi^*l_{\infty})=b_{E_i}$$ for $i=1,2$.

It follows that $Z_i=m_iZ_{-\deg}-b_{E_i}Z_{v_{E_i}}=m_i(Z_{-\deg}-Z_{v_{E_i}}).$ Then the coefficient of $Z_i$ of an irreducible exceptional divisor $E$ of $X$ is
$b_Em_i(1-\alpha(v_E\wedge v_{E_i})).$

Then we have $$(C_1\cdot C_2)=(\pi^*C_1\cdot \pi^*C_2)=\left((\pi^{\#}C_1+Z_1)\cdot (\pi^{\#}C_2+Z_2)\right)$$
$$=(Z_1\cdot Z_2)+(\pi^{\#}C_1\cdot Z_2)+(\pi^{\#}C_2\cdot Z_1)$$$$=m_1m_2(-1+\alpha(v_E\wedge v_{E_i}))+2m_1m_2(1-\alpha(v_{E_1}\wedge v_{E_2}))$$
$$=m_1m_2(1-\alpha(v_{E_1}\wedge v_{E_2}))=m_1m_2(1-\alpha(v_{C_1}\wedge v_{C_2})).$$
Then we have the following
\begin{pro}\label{prointersectionformalcurveinfi}If $C_1$, $C_2$ are two formal curves at infinity, then we have
$$(C_1\cdot C_2)=(C_1\cdot l_{\infty})(C_2\cdot l_{\infty})(1-\alpha(v_{C_1}\wedge v_{C_2})).$$
\end{pro}

\section{The dynamical Mordell-Lang Theorem near a superattracting point}
In this section, we study the dynamical Mordell-Lang Theorem when $C$ passing through a superattracting point.

We begin with the following simple property.
\begin{pro}\label{prolocaldml}
Let $X$ be a smooth projective variety defined of a valued field $(K,|\cdot|).$ Let $f:X\dashrightarrow X$ be a rational endomorphism on $X$ defined over $K$. Endow $X(K)$ the topology induced by $|\cdot|.$
Let $q$ be a $K$-point in $X$ satisfying
\begin{points}
\item $f(q)=q$;
\item $q\not\in I(f)$;
\item $df(q)=0$.
\end{points}
Let $C$ be a curve in $X$ satisfying $q\not\in C.$ Let $p$ be a $K$-point in $X$ satisfying $f^n(p)\not\in I(f)$ for all $n\geq 0.$
If there exists a sequence $n_j$ such that $f^{n_i}(p)\rightarrow q$ as $i\rightarrow \infty$, then the set $\{n|\,\,f^n(p)\in C\}$ is finite.
\end{pro}
\proof
Since $df(q)=0$ and $q\not\in C$, there exists a neighborhood $U$ of $q$ satisfying $U\cap I(f)=\emptyset$, $U\cap C(K)=\emptyset$ and $f(U)\subseteq U.$ Observe that $f^n(p)$ is defined over $K$ for all $n\geq 0.$
Since
$f^{n_i}(p)\rightarrow q$ as $i\rightarrow \infty$, there exists $m\geq 0$ such that $f^m(p)\in U$. It follows that $f^n(p)\in U$ for all $n\geq m$. Then we have $f^n(p)\not\in C$ for all $n\geq m$ which conclude our proof.
\endproof

Let $f:\mathbb{A}^2\rightarrow \mathbb{A}^2$ be a dominant polynomial morphism defined over $\overline{\mathbb{Q}}$.
Let $X$ be a compactificaiton defined over $\overline{\mathbb{Q}}$. Then $f$ extends to a rational endomorphism on $X$.
Let $q$ be a closed point in $X_{\infty}$ satisfying
\begin{points}
\item $f(q)=q$;
\item $q\not\in I(f)$;
\item $df(q)=0$.
\end{points}
Theorem \ref{thmfavre2007eigencaluation} implies that there exists an eigenvaluation $v_*\in V_q$ for $f$.
Then we have the following

\begin{thm}\label{thmlodml}
Let $C$ be an irreducible curve in $X$ containing $q$. Let $C_1$ be a branch of $C$ at $q$ such that the valuation $v_{C_1}\in V_q$ satisfies $f^n_{\d}(v_{C_1})\rightarrow v_*$ as $n\rightarrow \infty$. If $v_*$ is attracting and nondivisorial, and $v_{c_1}\neq v_*$, then $C$ is not preperiodic and for any point $p\in \mathbb{A}^2_{\overline{\mathbb{Q}}}$ which is not preperiodic under $f$, the set $\{n\in \mathbb{N}|\,\,f^n(p)\in C\}$ is finite.
\end{thm}

\proof[Proof of Theorem \ref{thmlodml}]
By contradiction, we suppose that there exists an infinite sequence $\{n_1<\cdots <n_i<n_{i+1}<\cdots \}$ such that $f^{n_i}(p)\in C$ for all $i\geq 1.$
By Theorem \ref{thmfavrejnormalform}, we may suppose that $f^n_{\d}(v)\rightarrow v_*$ for all $v\in V_q$ as $n\rightarrow \infty$.
By  Theorem \ref{thmfavrejnormalform} again, there exists a birational morphism $\pi:\widetilde{X}\to X$ which is an isomorphism above $X\setminus \{q\}$ such that the center $Q$ of $v_*$ is not contained in the strict transform $\pi^{\#}(C)$ of $C.$ Lift $f$ to a rational map $\widetilde{f}$ on $\widetilde{X}.$ We may suppose that $Q\not\in I(\widetilde{f})$, $\widetilde{f}(Q)=Q$ and $d\widetilde{f}(Q)=0$.
Set $\widetilde{C}:=\pi^{\#}C.$
Observe that $Q\in \widetilde{f}^{N_1}(\widetilde{C})$ for some $N_1\geq 1$ since $f^n_{\d}(v_{C_1})\rightarrow v_*$.  Set $\widetilde{p}:=\pi^{-1}(p).$


\medskip

Let $K$ be a number field such that $\widetilde{p}, \widetilde{X},\widetilde{f},Q$ and $\widetilde{C}$ are all defined over $K.$
For any place $v$ of $K$,
endow $X(K)$ with a metric $d_v$ induced by $v$.

We have $\widetilde{f}^{n_i}(\widetilde{p})\in \widetilde{C}$ and then $\widetilde{f}^{n_i+N_1}(\widetilde{p})\in f^{N_1}(\widetilde{C}).$
Since $Q$ is supperattacting, by \cite[Proposition 6.2]{Xieb}, there exists one place $v$ of $K$ such that $\widetilde{f}^n(\widetilde{p})\rightarrow Q$ with respect to the topology on $X$ with respect to $|\cdot|_v$. We conclude our proof by Proposition \ref{prolocaldml}.
\endproof

The following corollary comes from Theorem \ref{thmlodml} immediately.
\begin{cor}\label{cormonomial}
Let $f:\mathbb{A}^2\rightarrow \mathbb{A}^2$ be a dominant polynomial morphism defined over $\overline{\mathbb{Q}}$.
 Let $X$ be a compactificaiton defined over $\overline{\mathbb{Q}}$. Then $f$ extends to a rational endomorphism on $X$.
 Let $q$ be a closed point in $X\setminus \A^2$ such that in some local coordinates at $q$, $f$ takes form $(x^s,y^d)$ for $2\leq s\leq d-1.$
 Let $C$ be an irreducible curve in $X$ containing $q$. Let $p$ be a closed point in $\mathbb{A}^2(\overline{\mathbb{Q}})$. If $C$ is not fixed and $p$ is not preperiodic, then the set $\{n\in \mathbb{N}|\,\,f^n(p)\in C\}$ is finite.
 \end{cor}

\newpage

\part{The Dynamical Mordell-Lang Conjecture}
In this section, we give some basic observations on the Dynamical Mordell-Lang Conjecture and prove Theorem \ref{thmdmlmsv} as an application of these observations.

We first notice the following
\begin{pro}\label{proreducetofdom}Let $f:\mathbb{A}^2_{\overline{\mathbb{Q}}}\rightarrow \mathbb{A}^2_{\overline{\mathbb{Q}}}$ be a polynomial endomorphism defined over $\overline{\mathbb{Q}}$. Let $C$ be an irreducible curve in $\mathbb{A}^2_{\overline{\mathbb{Q}}}$ and $p$ be a closed point in $\mathbb{A}^2_{\overline{\mathbb{Q}}}$.
If $f$ is not dominant, then the set $\{n\in \mathbb{N}|\,\,f^n(p)\in C\}$ is a finite union of arithmetic progressions.
\end{pro}
\proof[Proof of Proposition \ref{proreducetofdom}]We suppose that the set $\{n\in \mathbb{N}|\,\,f^n(p)\in C\}$ is infinite.
If $f$ in the Main Theorem is not dominant, $f(\mathbb{A}^2_{\overline{\mathbb{Q}}})$ is an irreducible subvariety in $\overline{\mathbb{Q}}$ of dimension at most one.

If $\dim f(\mathbb{A}^2_{\overline{\mathbb{Q}}})=0$, then $f^{n}(p)=f(p)$ for all $n\geq 1$. Proposition \ref{proreducetofdom} holds in this case.

If $\dim f(\mathbb{A}^2_{\overline{\mathbb{Q}}})=1$ and $C=f(\mathbb{A}^2_{\overline{\mathbb{Q}}})$, then $f^n(p)\in C$ for $n\geq 1$ which conclude our proposition.

If $\dim f(\mathbb{A}^2_{\overline{\mathbb{Q}}})=1$ and $C\neq f(\mathbb{A}^2_{\overline{\mathbb{Q}}})$, then $C\cap f(\mathbb{A}^2_{\overline{\mathbb{Q}}})$ is finite. It follows that $p$ is preperiodic which concludes our proposition.
\endproof
In the rest of our paper, we suppose that $f$ is dominant.
\section{The DML property}\label{sedml}
As in \cite{Xie2014}, we introduce the following
\begin{defi}Let $X$ be a smooth surface defined over an algebraically closed field, and $f:X\dashrightarrow X$ be a rational endomorphism. We say that the pair $(X,f)$ satisfies the DML property for a curve $C$\index{DML property for a curve $C$} if
for any closed point $p\in X$ such that $f^n(p)\not\in I(f)$ for all $n\geq 0$, the set $\{n\in \mathbb{N}|\,\, f^n(p)\in C\}$ is a union of at most finitely many arithmetic progressions.

We say that the pair $(X,f)$ satisfies the DML property\index{DML property} if it satisfies the DML property for all curve $C$ in $X$.
\end{defi}

%
%

The DML property is equivalent to the following property.
\begin{pro}\cite[Proposition 4.2]{Xie2014}\label{prodefidml}Let $X$ be a smooth surface defined over an algebraically closed field, and $f:X\dashrightarrow X$ be a rational transformation.
The following statements are equivalent.
\begin{itemize}
\item[(1)] The pair $(X,f)$ satisfies the DML property.
\item[(2)] For any irreducible curve $C$ on $X$ and any closed point $p\in X$ such that $f^n(p)\not\in I(f)$ for all $n\geq 0$ and the set $\{n\in \mathbb{N}| f^n(p)\in C\}$ is infinite, then $p$ is preperiodic or $C$ is periodic.
\item[(3)] For any irreducible curve $C$ on $X$ and any closed point $p\in X$ such that $f^n(p)\not\in I(f)$ for all $n\geq 0$ and the set $\{n\in \mathbb{N}| f^n(p)\in C\}$ is infinite, then $p$ is preperiodic or $C$ is preperiodic.
\end{itemize}
\end{pro}
\proof[Proof of Proposition \ref{prodefidml}]
We first prove the equivalence of (1) and (2).

Suppose (1) holds. Let $C$ be any curve in $X$ and $p$ be a closed point in $X$ such that $f^n(p)\not\in I(f)$ for all $n\geq 0$. Assume that the set $\{n\in \mathbb{N}|\,\, f^n(p)\in C\}$ is infinite. The DML property of $(X,f)$ implies that there are integers $a>0$ and $b\geq 0$ such that $f^{an+b}(p)\in C$ for all $n\geq 0.$ If $p$ is not preperiodic, the set $O_{a,b}:=\{f^{an+b}(p)|\,\,n\geq 0\}$ is Zariski dense in $C$ and $f^{a}(O_{a,b})\subseteq O_{a,b}$. It follows that $f^a(C)\subseteq C$, hence $C$ is periodic.

Suppose (2) holds. If the set $S:=\{n\in \mathbb{N}|\,\, f^n(p)\in C\}$ is finite or $p$ is preperiodic, then there is nothing to prove. We may assume that $S$ is infinite and $p$ is not preperiodic. The property (2) implies that $C$ is periodic. There exists an integer $a>0$ such that $f^a(C)\subseteq C$. We may suppose that $f^{i}(C)\not\subseteq C$ for $1\leq i\leq a-1$. Since $p$ is not preperiodic, there exists $N\geq 0$, such that $f^n(p)\not\in (\cup_{1\leq i\leq a-1}f^{i}(C))\cap C$ for all $n\geq N.$
So $S\setminus\{1,\cdots,N-1\}$ takes form $\{an+b|\,\,n\geq 0\}$ where $b\geq 0$ is an integer, and it follows that $(X,f)$ satisfies the DML property.

So we only need to show that (3) implies (2).

We suppose that there exists a closed point $p\in X$ such that $f^n(p)\not\in I(f)$ for all $n\geq 0$ and the set $\{n\in \mathbb{N}| f^n(p)\in C\}$ is infinite. Moreover we may suppose that $C$ is preperiodic.

If $C$ is not periodic, there exist $m>0$ such that $f^m(C)$ is periodic. Then $\cup_{i=m}^{\infty}f^i(C)$ is a union of finitely many irreducible curves and  $f^n(p)\in \cup_{i=m}^{\infty}f^i(C)$ for $n\geq m.$ Since $C$ is not periodic, $C\cap \cup_{i=m}^{\infty}f^i(C)$ is finite. It follows that $p$ is preperiodic, which is a contradiction.
\endproof

\begin{thm}\label{thmbasicdmlra}Let $X$ be a smooth surface defined over an algebraically closed field, and $f:X\dashrightarrow X$ be a rational endomorphism, then the following properties hold.
\begin{points}
\item For any $m\geq 1$, $(X,f)$ satisfies the DML property if and only if $(X,f^m)$ satisfies the DML property.
\item Suppose $U$ is an open subset of $X$ such that the restriction $f_{|U}:U\rightarrow U$ is a morphism. Then $(X,f)$ satisfies the DML property, if and only if $(U,f_{|U})$ satisfies the DML property.
\item Suppose $\pi:X'\rightarrow X$ is a generic finite morphism between smooth projective surfaces, and $f:X\dashrightarrow X$, $f':X'\dashrightarrow X'$ are two rational maps satisfying $\pi\circ f'=f\circ\pi$. For any curve $C$ in $X$, if the pair $(X',f')$ satisfies the DML property for $\pi^{-1}(C)$, then $(X,f)$ satisfies the DML property for $C$.
\end{points}
\end{thm}

\proof[Proof of Theorem \ref{thmbasicdmlra}]
\proof[$( \rm i)$]
The "only if" part is trivial, so that we only have to deal with the "if" part. We assume that $(X,f^m)$ satisfies the DML property. Let $C$ be a curve in $X$ and $p$ be a point in $X$ such that $f^n(p)\not\in I(f)$ for all $n\geq 0$. Suppose that the set $\{n\in \mathbb{N}|\,\, f^n(p)\in C\}$ is infinite. Since $$\{n\in \mathbb{N}|\,\, f^n(p)\in C\}=\bigcup_{i=0}^{m-1}\{n\in \mathbb{N}|\,\, f^{nm}(f^i(p))\in C\},$$ then for some $i$, the set $\{n\in \mathbb{N}| \,\,f^{nm}(f^i(p))\in C\}$ is also infinite. Since $(X,f^m)$ satisfies the DML property, $C$ is periodic or $f^i(p)$ is preperiodic. It follows that $C$ is periodic or $p$ is preperiodic.
\proof[$(\rm ii)$] If $(X,f)$ satisfies the DML property, since $f_{|U}:U\rightarrow U$ is a morphism, $(U,f_{|U})$ satisfies the DML property.

Conversely suppose that $(U,f_{|U})$ satisfies the DML property. Let $C$ be an irreducible curve in $X$, $p$ be a closed point in $X$ such that $f^n(p)\not\in I(f)$ for all $n\geq 0$ and the set $\{n\in \mathbb{N}| f^n(p)\in C\}$ is infinite. The set $E=X-U$ is a proper closed subvariety of $X$. If $p\in U$, then we have that $C\not\subseteq E$. Since $(U,f_{|U})$ satisfies the DML property, we have either $p$ is preperiodic or $C$ is periodic. Otherwise, we may assume that for all $n\geq 0,$ $f^n(p)\in E,$ then the Zariski closure $D$ of $\{f^n(p)|\,\,n\geq 0\}$, is contained in $E$. We assume that $p$ is not preperiodic, then $C\subseteq D$. Since $D$ is fixed, we have that $C$ is periodic.
\proof[$(\rm iii)$]Let $p\in X$ be a nonpreperiodic point satisfying $f^n(p)\not\in I(f)$ of all $n\geq 0$ and $C$ be an irreducible curve in $X$. Suppose that the set $\{n\in \mathbb{N}|\,\,f^n(p)\in C\}$ is infinite. The set $I(f')$ is finite, so its image $\pi(I(f'))$ is finite. Let $S$ be the set of point $x$ in $X$ satisfying $\pi^{-1}(x)$ is infinite. Then $S$ is finite.
 Since $p$ is not preperiodic, there exists $N\geq 0$ such that $f^n(p)\not\in \pi(I(f'))\cup S$ for all $n\geq N.$ By replacing $p$ by $f^N(p)$, we may suppose that $N=0.$ Let $q$ be a point in $\pi^{-1}(p).$ We have $f^{'n}(q)\not\in I(f')$ and the set $\{n\in \mathbb{N}|\,\,f^{'n}(q)\in \pi^{-1}(C)\setminus \pi^{-1}(S)\}$ is infinite. Then there exists an irreducible component $C'$ of $\pi^{-1}(C)$ satisfying $\pi(C')=C$ and the set  $\{n\in \mathbb{N}|\,\,f^{'n}(q)\in C'\}$ is infinite. We see that $q$ is not preperiodic, so $C'$ is periodic. It follows that $C$ is periodic, which concludes our proof.
\endproof
%
%
%
%
%
%
%
%
\section{Constraints on the geometry of the target curve}
In this section the situation is as follows:
$f$ is a dominant polynomial map of $\A^2$ defined over $\bar{\Q}$, $C$ is an irreducible curve in $\mathbb{A}^2_{\overline{\mathbb{Q}}}$  containing infinitely many iterate of a non-preperiodic point $p\in \mathbb{A}^2(\overline{\mathbb{Q}})$. The follows theorem gives us some constraints on the geometry of $C$.

\begin{thm}\label{thmasequencecurves}
Let $f$ be a dominant endomorphism on $\mathbb{A}^2$ defined over $\overline{\mathbb{Q}}$, $C$ an irreducible curve in $\mathbb{A}^2_{\overline{\mathbb{Q}}}$ and $p$ be a closed point in $\mathbb{A}^2(\overline{\mathbb{Q}})$.

If the set $\{n\geq 0|f^n(p)\in C\}$ is infinite and $p$ is not $f$-preperiodic, then there exists a sequence of rational curves $\{C_i\}_{i\in \mathbb{Z}}$ with at most two branches at infinity such that
\begin{points}
\item $C^0=C$;
\item $f(C^i)=C^{i+1}$;
\item for all $i\in \mathbb{Z}$, the set $\{n\geq 0|f^n(p)\in C^i\}$ is infinite.
\end{points}
\end{thm}
Since $f$ is polynomial, the number branches of $C^i$ is increasing as $i\rightarrow -\infty$ but bounded by two. So there exists $N\leq 0$, such that the number of branches of $C^i$ is stable when $i\leq N$. So we have the following

\rem By replacing $C$ by $C^j$ for some $j\leq 0$, we may suppose that for all $j\leq 0$, number of places of $C^j$ at infinity are the same number $s\in\{1,2\}$. Let $C^j_i$'s be branches of $C^j$, we may suppose that $f(C^j_i)=C_i^{j+1}$ for $j\leq -1$ and $1\leq i\leq s.$
\endrem

The following theorem shows how to apply this sequence of curves to the Dynamical Mordell-Lang Conjecture.

\begin{thm}\label{thmseqcurvebounddegdml}
Let $f$ be a dominant endomorphism on $\mathbb{A}^2$ defined over $\overline{\mathbb{Q}}$ that is not birational.

Pick any smooth projective compactification $X$ of $\mathbb{A}^2$ and suppose that
there exists a sequence of irreducible curves in $X$ satisfying $f(C^i)=C^{i+1}$ for $i\leq -1$
and such that $\sup_{i\in \Z_-}(C^i\cdot L)$ is bounded for some ample line bundle $L\to X$.
Then the pair $(X,f)$ satisfies the DML property for the curve $C^i$ for some $i\leq 0.$
\end{thm}
%

%
%
%
%
%
%
%
%
%
%
%

\subsection{Proof of Theorem \ref{thmasequencecurves}}
We first fix some notations:
\begin{points}
\item[$\d$]$K$ is a number field;
\item[$\d$]$\mathcal{M}_K$\index{$\mathcal{M}_K$} is the set of places on $K$;
\item[$\d$]$\mathcal{M}_K^{\infty}$\index{$\mathcal{M}_K^{\infty}$} is the set of archimedean places on $K$;
\item[$\d$]$S$ is a finite set of places
of $K$ containing all the archimedean places;
\item[$\d$]$O_{K,S}$\index{$O_{K,S}$} is the ring of $S$-integers.
\end{points}
Theorem \ref{thmasequencecurves} is a corollary of the Siegel's Theorem (see \cite{Hindry2000} for details).

\begin{thm}[Siegel's Theorem]\label{thmsiegel}
Let $C$ be a curve over a number field $K$ and $g\in K(C)$ be
a nonconstant rational function on $C$. If either $C$ is not rational or
$g$ has at least three distinct poles, then the set $\{p\in C(K)|g(p)\in O_{K,S}\}$
is finite.
\end{thm}

Next we recall two obvious facts.
\begin{points}
\item[$\bullet$] If $C\in \mathbb{A}^2(K)$ is a plane curve which has at least $3$ branches at the infinity, by taking $g=ax+by$ where $x,y$ are the coordinate functions and $a,b$ are two general integers, Siegel's Theorem shows that the set of $S$-integral points of $C$ i.e. $\{(x,y)\in C(K)|x,y\in O_S\}$ is finite.
\item[$\bullet$] If $f:\mathbb{A}^2\rightarrow \mathbb{A}^2$ is a polynomial endomorphism of $\mathbb{A}^2$ whose coefficients are all contained in $O_{K,S}$, and $p\in \mathbb{A}^2(K)$ is a $S$-integer point. For any $n\geq 0,$ $f^n(p)$ is a $S$-integer point.
\end{points}

Then we have the following

%
%
\proof[Proof of Theorem \ref{thmasequencecurves}]
We may suppose that there exists a number field $K$, such that $f$ and $p$ are all defined over $K$. Further we may suppose that there exists a finite set $S$ of $\mathcal{M}_K$ containing $\mathcal{M}_K^{\infty}$ such that all coordinates of $p$ and all coefficients of $f$ are contained in $O_{K,S}$. It follows that all points in the orbit of $p$ are $S$-integral points.

For $i\geq 0$, we just set $C^i:=f^i(C)$. For $j\leq -1$, we
construct this sequence by induction. If we have $C^i$ for some $i\leq 0$ such that the set $\{n\geq 0|\,\,f^n(p)\in C^i\}$ is infinite. Then the set $\{n\geq 0|\,\,f^n(p)\in f^{-1}(C^i)\}$ is infinite. There exists an irreducible component $C^{i-1}$ of $f^{-1}(C^i)$ such that the set $\{n\geq 0|\,\,f^n(p)\in C^{i-1}\}$ is infinite. By Theorem \ref{thmsiegel}, $C^{i-1}$ is rational,  has at most two branches at infinite and satisfies $f(C^{i-1})=C^i$.

\endproof

\subsection{Proof of Theorem \ref{thmseqcurvebounddegdml}}Theorem \ref{thmseqcurvebounddegdml} is the corollary  of the following more general result.

\begin{pro}\label{procurveboundeddegfibration}Let $X$ be a smooth rational surface defined over an algebraically closed field and $f:X\dashrightarrow X$ be a dominant rational endomorphism on $X$ with $\la_2\geq 2$. Let $L$ be an ample line bundle and let $\{C^i\}_{i\leq 0}$ be a sequence of distinct curves in $X$ satisfying $f(C^i)=C^{i+1}$ for $i\leq -1$.

If that there exists $M>0$ such that $(C^i\cdot L)\leq M$, then up to a positive iterate there exists a generic finite cover $g:X'\rightarrow X$ with a rational endomorphism $f':X'\dashrightarrow X'$ satisfying  $f\circ g=g\circ f'$ such that
 we have that $f'$ preserves a rational fibration $\pi$ and
 for some $i\leq 0$, every component of $g^{-1}(C^i)$ contains in a fiber of $\pi$.
%
\end{pro}

\proof[Proof of Theorem \ref{thmseqcurvebounddegdml}]Let $p$ be a closed point in $X$ such that $f^n(p)\not\in I(f)$ for all $n\geq 0$ and the set $\{n\geq 0|\,\,f^n(p)\in C\}$ is infinite. By Proposition \ref{prodefidml}, we may suppose that $C$ is not periodic and $p$ is not preperiodic. Then the curves $C^{i}$'s, $i\leq -1$ are distinct. Since there exists $M>0$ such that $(C^i\cdot L)\leq M$, by Proposition \ref{procurveboundeddegfibration}, up to a positive iterate there exists a generic finite cover $g:X'\rightarrow X$ with a rational endomorphism $f':X'\dashrightarrow X'$ satisfying  $f\circ g=g\circ f'$ such that
 we have that $f'$ preserves a rational fibration $\pi$ and
 for some $i\leq 0$, every component of $g^{-1}(C^i)$ is contained in a fiber of $\pi$. Pick any point $q\in g^{-1}(p)$. By replacing $p$ by $f^n(p)$ for some $n\geq 0$, we may suppose that $(f')^n(q)\not\in I(f')$ for all $n\geq 0$. Then set of $n\geq 0$ such that $(f')^n(p)\in g^{-1}(C)$ is infinite. Pick $C'$ an irreducible component of  $g^{-1}(C)$ for which the set $\{n\geq 0|\,\, (f')^{n}\in C'\}$ is infinite. Then $\pi(C')$ is a periodic points. It follows that $C'$ is periodic and then $C$ is periodic.
\endproof

\proof[Proof of Proposition \ref{procurveboundeddegfibration}]
There exist a smooth projective surface $\Gamma$, a birational morphism $\pi_1:\Gamma\rightarrow X$ and morphism $\pi_2:\Gamma\rightarrow X$ satisfying $f=\pi_2\circ\pi_1^{-1}.$ We denote by $f_*$ the map $\pi_{2*}\circ\pi^*_1:\Div X\rightarrow \Div X$. Let $E_{\pi_1}$ be the union of exceptional irreducible divisors of $\pi_1$ and $\mathfrak{E}$ be the set of effective divisors in $X$ supported by $\pi_2(E_{\pi_1}).$ It follows that for any curve $C$ in $X$, there exists $D\in \mathfrak{E}$ such that $f_*C=\deg(f|_C)f(C)+D$.

For any effective line bundle $K\in \Pic(X)$, the projective space $H_K:=\mathbb{P}(H^0(K))$ parameterizing the curves $C$ in the linear system $|K|$.
Since $\Pic^0(X)=0$, for any $l\geq 0$, there are only finitely many effective line bundle satisfying $(K\cdot L)\leq l.$

Then $H^l:=\coprod_{(K\cdot L)\leq l}H_K$ is a finite union of projective spaces and
it parameterizing the curves $C$ in $X$ satisfying $(C\cdot L)\leq l$.

There exists $d\geq 1$ such that $dL-f^*L$ is nef. Then for any curve $C$ in $X$, we have $(f_*C\cdot L)=(C\cdot f^*L)\leq d(C\cdot L).$
It follows that $f_*$ induce a morphism $F:H^l\rightarrow H^{dl}$ by $C\rightarrow f_*C$ for all $l\geq 1.$ For all $l\geq 1$, $a\in \mathbb{Z}^+$ and $D\in \mathfrak{E}$, there exists an embedding $i_{a,D}:H_l\rightarrow H_{al+(D\cdot L)}$ by $C\mapsto aC+D$. Let $Z_1,\cdots,Z_m$ be all irreducible components of the Zariski closure of $\{C^j\}_{j\leq -1}$ in $H^M$ whose dimensions are maximal. For any $i\in\{1,\cdots,m\}$, there exists $l\leq M$ such that $(C\cdot L)=l$ for all $C\in Z_i$. Let $S$ be the finite set of pairs $(a,D)$ where $a\in \mathbb{Z}^+$, $D\in \mathfrak{E}$ satisfying $al+(D\cdot L)\leq dM$. Then we have $F(Z_i)\subseteq \cup_{j=1,\cdots,m}\cup_{(a,D)\in S}i_{a,D}(Z_j)$. It follows that there exists a unique $j_i\in \{1,\cdots,m\}$, and a unique $(a,D)\in S$ such that $F(Z_i)=i_{a,D}(Z_{j_i}).$ Observe that, the map $i\mapsto j_i$ is an one to one map of $\{1,\cdots,m\}$. By replacing $f$ by a positive iterate, we may suppose that $j_i=i$ and $F(Z_i)\subseteq i_{a_{Z_i},D_{Z_i}}Z_i$ for all $i=1,\cdots,m.$ Set $Y:=Z_1$, $a=a_{Z_1}$, $D=D_{Z_i}$ and $T=i_{a_{Z_1},D_{Z_1}}^{-1}\circ F|_Y$.

Observe that $Y$ is a projective variety and $T$ is an endomorphism on $Y.$
Let $K$ be the line bundle such that $\mathbb{P}(H^0(X,K))$ contains $Y$, $H$ the hyperplane line bundle on $\mathbb{P}(H^0(K))$ and $H'$ be the hyperplane line bundle on $\mathbb{P}(H^0(f_*K)).$ Observer that $i_{a,D}^*H'=H^{\otimes a}$ and $F^*(H)=H^{\otimes \la_2}$ where $\la_2$ is the topological degree of $f$. It follows that $T^*(H|_Y^{\otimes a})=H|_Y^{\otimes \la_2}$. It follows that the topological degree of $T$ is $(\la_2/a)^{\dim Y}$. Then $\la_2/a$ is a positive integer.

Let $S$ be the subvariety of $Y\times X$ whose set of closed point is $\{(C,q)|\,\, q\in C\}.$
Denote by $p_1:S\rightarrow Y$ and $p_2:S\rightarrow X$ the projections to the first and the second coordinates. For any $i\geq N$, $C_i$ is a fiber of $\pi_1$ and it is irreducible. Set $R$ be the infinite set of $j\leq 0$ such that $C^j\in Y$. Since $\{C^j\}_{j\in R}$ is dense in $Y$, we have the following properties.
\begin{points}\item[(1)] The generic fiber of $p_1$ is irreducible.
\item[(2)] Every fiber of $p_1$ is dimensional $1.$
\item[(3)] The restriction of $p_2$ on a fiber of $p_1$ is an embedding.
\item[(4)] The images of two different $p_1$-fibers by $p_2$ are different.
\end{points}
 Observe that $S$ is invariant by the rational endomorphism $T\times f: Y\times X\dashrightarrow Y\times X$ and then denote by $f_S$ the restriction of $T\times f$ to $S.$

Then the diagram
$$\xymatrix{
X \ar[rr]^{f} &   & X \\
S \ar[rr]^{f_S}\ar[d]_{p_1}\ar[u]^{p_2} &   & S \ar[d]_{p_1}\ar[u]^{p_2} \\
Y \ar[rr]^{T}        &   & Y
}$$
commutes.

For a general point $C\in Y$, set $T^{-1}(C)=\{C_1,\cdots, C_{(\la_2/a)^{\dim Y}}\}$. If we view them as curves in $X$, we have $f_*(C_i)=aC+D$ for $i=1,\cdots,(\la_2/a)^{\dim Y}.$ For a general points $p$ in $C$, the number of its preimages by $f|_{C_i}$ is $a$. So we have $\la_2=\#f^{-1}(q)\geq a(\la_2/a)^{\dim Y}$.

If $\la_2/a\geq 2$, we have $\dim Y=1$. Then $S$ is a surface which concludes our Proposition.

Otherwise, we have $\la_2=a$. Then $T$ is an automorphism.
Since $\la_2\geq 2$,
by replacing $f$ by a positive iterate, we may suppose that $\la_2>(K\cdot K)$. Let $p$ be a general point in $X$ and $C$ be a point in $Y$ such that $p\in C$, we have $\#(f^{-1}(p)\cap T^{-1}(C))=a=\#f^{-1}(p)$. It follows that $f^{-1}(p)\subseteq T^{-1}(C).$
 If there exists another point $C'\in Y$ containing $p$, then we have $f^{-1}(p)\subseteq T^{-1}(C)\cap T^{-1}(C')$. It follows that $(K\cdot K)=(T^{-1}(C)\cdot T^{-1}(C'))\geq \#(T^{-1}(C)\cap T^{-1}(C'))\geq\la_2$ which contradicts our assumption. Then there are only one $C\in Y$ containing $p$. In other words, $p_2$ is birational. Then $S$ is a surface and $Y$ is a curve which conclude our Proposition.
 \endproof

Finally, we prove a technical result which shows that how to use Theorem \ref{thmasequencecurves} to construct a sequence of curves satisfying the conditions in
Proposition \ref{procurveboundeddegfibration}.

Let $f:\mathbb{A}^2_k\rightarrow \mathbb{A}^2_k$ be a dominate polynomial endomorphism on $\mathbb{A}^2_k.$
Let $X\in \mathcal{C}$ be a compactification of $\mathbb{A}^2_k$ and we extend $f$ to a rational endomorphism of $X.$
There exists a smooth projective surface $Y$ a birational morphism $\pi_1:Y\to X$ and a morphism $\pi_2:Y\to X$ satisfying $f=\pi_2\circ\pi_1^{-1}$. Set $f^*=\pi_{1*}\circ \pi_2^*:\Div(X)\to \Div(X)$.

\begin{defi}Let $E'$ be any irreducible curve in $X\setminus \A^2_k$.  If $\pi_2$ contracts all irreducible curves in $\pi_2^*(E')$ except $\pi_1^{\#}E'$, then we say that $E'$ is {\em totally invariant} by $f$.
\end{defi}

\begin{rem}In fact, $E'$ is totally invariant if and only if $v_{E'}$ is totally invariant by $f_{\d}$. Moreover we have $f^*E'=d(f,v_{E'})E'.$
\end{rem}

\begin{defi}Let $C$ and $C'$ be two distinct irreducible curves in a projective surface $X$, and $B$ a set of points in $X$. Denote by $(C \cdot C' \setminus B)$ the sum of local intersection numbers of $C$ and $C'$ outside $B$.
\end{defi}

Let $\{C^j\}_{j\leq 0}$ be a sequence of curves with $s=1$ or $2$ branches at infinity satisfying
$f(C^j)=C^{j+1}$ for $j\leq -1$.
Let $C^j_i$'s be branches of $C^j$ and suppose that $f(C^j_i)=C_i^{j+1}$ for $j\leq -1$ and $1\leq i\leq s.$

\begin{thm}\label{prototinvanonexcepdml}
Let $E$ be the union of all totally invariant curves lying in the divisor at infinity $X\setminus \A^2_k$ and suppose that $f^*E=dE$ for some $d\geq 1$.

Denote by $\mathcal{E} $ the subset of points $q \in E$ whose orbit under $f|_E$ is periodic, and contains either a point of indeterminacy of $f$, or a singular point of $X_\infty$, or a
critical point for $f|_E$.

Let $G$ be the set of index $1\leq i\leq s$ such that for all $j \leq 0$ the center $q^j_i$ of $C^j_i$ are contained in $E\setminus I(f)$.
Let $D$ be an effective ample divisor supposed by $X\setminus \A^2_k$ and $D_E$ part of $D$ supported by $E$.

If $C^0\cap (E\setminus \mathcal{E})\neq\emptyset$ and  $\sum_{i\in G}(D_E\cdot C^j_i)\geq \varepsilon(D\cdot C^j)$ for some $\varepsilon >0$ and all $j\leq 0$, then $\sup_{i\in \Z_-}(C^i\cdot D)$ is bounded.
\end{thm}

%
%
%
%
%
%
%
%
\proof[Proof of Proposition \ref{prototinvanonexcepdml}]Let $E_1,\cdots,E_m$ be all irreducible components of $E.$
For any $i=1,\cdots,m$ and $y\in E_i$, set $U(y)$ open set in $V_{\infty}$ consisting by the valuations presented the vector corresponding to $y$.
Since $v_{E_i}$ is totally invariant under $f_{\d}$, for any valuation $v\in U(t)$ satisfying $d(f,v)>0$, we have $f_{\d}(v)\in U(f|_{E_i}(y))$.
Set $q^j_i:=C^j_i\cap X_{\infty}.$

Let $E'$ be an irreducible component of $E$. Observer that if $q^{j}_i\in E,$ then $q^{j+1}_i=f|_E(q^{j}_i).$
Since $E$ is totally invariant, we have $q^{j}_i\in E$ if and only if $q^{0}_i\in E$. If $q^{j}_i\in E,$ then $q^{j+1}_i=f|_E(q^{j}_i).$

We may suppose that $q^{0}_1\in E_1\setminus \mathcal{E}.$
By replacing $C$ by $C^{-l}$ for $l$ large enough, we may suppose that for all $j\leq 0$, we have
$f|_{E_1}$ is not ramified at $q^{0}_1$.

 Pick a neighborhood $U_{j,1}$ of $q^{j}_1$ for $j\leq 0$ and $i\in G$, we may suppose that in some local coordinate $f: U_{j,1}\rightarrow U_{j+1,1}$ has form $(x,y)\mapsto (x,y^d)$. In these coordinates, $E_1=\{y=0\}$. It follows that $\deg f|_{C^j_1}$ is at most $d$. Since $C$ is irreducible, we have $\deg f|_{C^j}=\deg f|_{C^j_1}\leq d.$

 Pick $E'$ an irreducible component of $E$. If $q^0_i\in E'$,
then we have $$(C_i^j\cdot E')=1/d(C_i^j\cdot f^* E')=(\deg(f|_{C_i^j})/d)(C_i^{j+1}\cdot E')\leq (C_i^{j+1}\cdot E')$$ for all $i\leq -1.$

 If $q^0_i\not\in E'$, then $(C_i^j\cdot E')=0$ for all $j\leq 0.$

We have $$(D\cdot C^{j})\leq1/\varepsilon\sum_{i\in G}(D_E,C^{j}_i)\leq 1/\varepsilon\sum_{i\in G}(D_E,C^{j+1}_i)\leq 1/\varepsilon\sum_{i\in G}(D_E\cdot C^0_i)$$ for all $j\leq -1.$ By Proposition \ref{procurveboundeddegfibration}, we conclude our Proposition.
\endproof

%
%

%
%
\section{The proof of Theorem \ref{thmdmlmsv}}
In this section, we denote by $k:=\overline{\mathbb{Q}}$ the field of algebraic numbers.

\medskip

We first recall the setting:

Let $f:=(F_1(x_1),\cdots,F_m(x_m))$ be an endomorphism on $\mathbb{A}^m$ defined over $k$.
Let $C$ be any irreducible curve in $\mathbb{A}^m$ defined over $k$ and $p$ be any point in $\mathbb{A}^m(k)$. We need to show that the set $\{n\geq 0|f^n(p)\in C\}$ is a finite union of arithmetic progressions.

\medskip

When $m=1$, the statement is trivial.

\medskip
\subsection{The case $m=2$.}
When $m=2$, Theorem \ref{thmdmlmsv} immediately comes from our Main theorem. Here we give a direct proof of it to see how can we use the results in Part 3 to the Dynamical Mordell-Lang Conjecture.

\smallskip

Since $F_1,F_2$ can extend to endomorphisms of $\P^1_k$, $f$ extends to an endomorphism on $X:=\P^1_k\times \P^1_k$. Then $f$ preserves the two projection $\pi_i$, $i=1,2$ the the $i-$th coordinate. Denote by $d_i$ the degree of $F_i$ for $i=1,2$.  Suppose that $C$ is irreducible and the set $\{n|\,\,f^n(p)\in C\}$ is infinite.

\smallskip

We first treat the case $d_1\neq d_2.$
We may suppose that $d_1>d_2.$

If $C$ is a fiber of $\pi_1$ or $\pi_2$, the conclusion is trivial.
So we may assume that $\pi_i|_C$ is dominate with degree $c_i>0$ for $i=1,2.$ Set $x_i^n:=F_i^n(x_i^0)$ and $p^n:=f^n(p)=(x_1^n,x_2^n)$ for  $i=1,2$ and $n\geq 0.$

If there exists one $i=1,2$ such that $x_i^0$ is preperiodic, by replacing $f$ by some positive iterate and $p$ by some $p^k$ for $k\geq 0$, we may suppose that $x_i^0$ is fixed by $F_i$. Then we conclude our theorem by induction hypothesise.

We suppose that $x_i^0$ is not $F_i$ preperiodic for $i=1,2$. The set $\{n| f^n(p)\in C\}$ can be written as an increase sequence $\{n_k\}_{k\geq 1}$.

Denote by $h$ the naive height function on $\mathbb{P}^1$ which is a Weil height with respect to the ample line bundle $L:=O_{\mathbb{P}^1}(1)$. Then $h\circ\pi_i$ is a Weil height with respect to the line bundle $L_i=\pi_i^*L|_C$ which has degree $c_i$ for $i=1,2.$ Then we have $$h\circ\pi_i(p^{n})=h(x_i^{n})=h(F^{n}(x_i^0))$$ for all $i=1,2$ and $n\geq 0.$

For any $i=1,2,$ $x_i^0$ is not $F_i$-preperiodic, hence here exists $C_i>0$ ,$D_i>0$ such that $$C_i(d_i-1/3)^n-D_i\leq h(F_i^{n}(x_1^0))\leq C_i(d_i+1/3)^n+D_i.$$
Since $d_1>d_2$, we have $d_1-1/3>d_2+1/3$, so we have $$\lim_{k\rightarrow \infty}h\circ\pi_1(p_{n_k})/h\circ\pi_2(p_{n_k})=+\infty.$$
This contradicts the following

\begin{lem}[\cite{Hindry2000}]\label{lemhoverh}Let $C$ be a projective curve over a number field $K$ and $L_1,L_2$ be two ample line bundles on $X$ over $K$ with degrees $d_1$ and $d_2.$ If $h_1$, $h_2$ are Weil heights with respect to $L_1$ and $L_2$ and $\{x_n\}_{n\geq 0}$ is an infinite set of points in $C(K)$, then we have $\lim_{n\rightarrow\infty}h_1(x_n)/h_2(x_n)=d_1/d_2$.
\end{lem}

\smallskip

Then we treat the case $d:=d_1=d_2$. If $d=1$, we have that $f$ is an automorphism. Then we may conclude our Theorem by \cite{Bell2010} in this case.  So we may suppose that $d\geq 2.$  Let $E_i$ be the section of $\pi_i$ at infinity for $i=1,2$. Then $f|_{E_i}=F_i$ for $i=1,2$ and $X\setminus\A^2_k=E_1\cup E_2.$

If $C$ is a fiber of $\pi_1$ or $\pi_2$, the
conclusion trivially holds.  So we may suppose that $C\cap E_1\neq\emptyset$. If $C$ passes the point $O:=E_1\cap E_2$, we conclude by the following
\begin{lem}\label{lemddtotinv}Let $f:\A^2_k\to \A^2_k$ be a polynomial endomorphism on $\A^2_k$ and $C$ be a curve in $\A^2_k$. let $X$ be a compactification of $\A^2_k$ in $\mathcal{C}$ such that $f$ extends to an endomorphism on $X$. Suppose that $f^*(X\setminus \A_k^2)=d(X\setminus \A_k^2)$ for some $d\geq 2.$
Let $q$ be a point in $X\setminus \A^2_k$ which is totally invariant and locally $f$ takes form $(x^d,y^d)$.  If $C$ passes through $q$, then $(X,f)$ satisfies the DML property for the curve $C$.
\end{lem}

So we may suppose that $C\cap (E_1\setminus\{O\})\neq\emptyset$.

By Theorem \ref{thmasequencecurves}, we construct a sequence of rational curves $\{C_i\}_{i\in \mathbb{Z}}$ with at most two places at infinity such that
\begin{points}
\item $C^0=C$;
\item $f(C^i)=C^{i+1}$;
\item for all $i\in \mathbb{Z}$, the set $\{n\geq 0|f^n(p)\in C^i\}$ is infinite.
\end{points}

By replacing $C$ by $C^j$ for some $j\leq 0$, we may suppose that for all $j\leq 0$, number of places of $C^j$ at infinity are the same number $s\in\{1,2\}$. Let $C^j_i$'s be branches of $C^j$, we may suppose that $f(C^j_i)=C_i^{j+1}$ for $j\leq -1$ and $1\leq i\leq s.$

If $C$ passes through a $f|_{E_1}$-critical periodic point $q\in E\setminus \{O\},$ by replacing $f$ by a positive iterate, we may suppose that $q$ is fixed by $f$. In a suitable local coordinate at $q$, $f$ takes form $(x,y)\mapsto (x^s,y^d)$ where $2\leq s\leq d$. When $s=d$, we conclude by Lemma \ref{lemddtotinv}. When $2\leq s\leq d-1$, we conclude by Corollary \ref{cormonomial}.

Then we may suppose that there exists a point $q_1\in C\cap (E_1\setminus \{O\})$ which is not a critical periodic point for $f$. Set $D=E_1+E_2$. Observe that $D$ is ample.
By  Proposition \ref{prototinvanonexcepdml}, $\sup_{i\in \Z_-}(C^i\cdot D)$ is bounded. Then we conclude the proof in this case by Theorem \ref{thmseqcurvebounddegdml}.

\proof[Proof of Lemma \ref{lemddtotinv}]
By Theorem \ref{thmasequencecurves}, we construct a sequence of rational curves $\{C_i\}_{i\in \mathbb{Z}}$ with at most two places at infinity such that
\begin{points}
\item $C^0=C$;
\item $f(C^i)=C^{i+1}$;
\item for all $i\in \mathbb{Z}$, the set $\{n\geq 0|f^n(p)\in C^i\}$ is infinite.
\end{points}

Since $q$ is totally invariant, $C^j$ passes through $q$ for all $j\leq 0$

By replacing $C$ by $C^j$ for some $j\leq 0$, we may suppose that for all $j\leq 0$, number of places of $C^j$ at infinity are the same number $s\in\{1,2\}$. Let $C^j_i$'s be branches of $C^j$, we may suppose that $f(C^j_i)=C_i^{j+1}$ for $j\leq -1$ and $1\leq i\leq s.$

\smallskip

Let $C_1$ be a branch of $C$ at $q$. Let $E_1$, $E_2$ be the formal curve locally defined by $\{x=0\}$ and $\{y=0\}$. Since $E_1,E_2$ are fixed by $f$, we may suppose that $C_1$ is different from both $E_1$ and $E_2$.

We define a sequence of surfaces $\pi_k:X_k\rightarrow X$ by induction:
\begin{points}\item Set $X_0:=X$ and $\pi_0:=\id$.
\item Suppose that we have $X_0,\cdots, X_k$. If $C_1$ does not pass through any singular point of $\pi_k^{-1}(E_1\cup E_2)$, we stop our progression.
\item If $C_1$ is passing through one singular point of $\pi_k^{-1}(E_1\cup E_2)$, let $X_{k+1}$ be the surface defined by blowup at this point in $X_k$.
\item Denote by $C_1$ the
strict transform of $C_1$ in $X_k$.
Then return to (i).
\end{points}
This progression terminates in finitely many steps and we get surfaces $X_0,\cdots, X_l$ for $l\geq 0.$

It is easy to see that $f$ is an endomorphism on $X_l.$ At any singular point of $\pi_k^{-1}(E_1\cup E_2)$, $f$ locally conjugates to $(x,y)\mapsto (x^d,y^d).$ Let $E$ be the unique exceptional curve of $\pi_l$ which intersects $C_1$ at one point. We see that $E$ is totally invariant and $f|_E$ can be written as $z\rightarrow z^d.$ All the ramified points of $f|_{E}$ are singular in $\pi_k^{-1}(E_1\cup E_2)$.

Let $D$ be any ample divisor supported by $X_l\setminus \A^2_k$,
by Proposition \ref{prototinvanonexcepdml}, $\sup_{i\in \Z_-}(C^i\cdot D)$ is bounded. Then we conclude our Lemma by Theorem \ref{thmseqcurvebounddegdml}.
\endproof

\subsection{The higher dimensional case}
In the case $m\geq 3$,  we prove this theorem by induction.
Suppose that $C$ is irreducible and the set $\{n|\,\,f^n(p)\in C\}$ is infinite.
Write $p=(p_1,\cdots, p_m)$ and denote by $\pi_i$ the projection from $\mathbb{A}^m$ to the $i$-th coordinate.

If there exists $1\leq i\leq m$ such that $p_i$ is $F_i$ preperiodic, by replacing $f$ by some positive iterate $f^l$ and $p$ by $f^l(p),$ we may suppose that $p_i$ is fixed. Then $f^n(p)\in \pi_i^{-1}(p_i)$ for all $n\geq 0.$ If $C$ is not contained in $\pi_i^{-1}(p_i),$ we have $C\cap \pi^{-1}_i(p_i)$ is finite and then $p$ is preperiodic.  If $C\subseteq \pi^{-1}(p_i)$. By replacing $\mathbb{A}^m$ by $\pi^{-1}(p_i)\simeq \mathbb{A}^{m-1}$, we conclude our theorem by the induction hypotheses.

So we may suppose that for all $1\leq i\leq m$, $p_i$ is non preperiodic by $F_i.$ It follows that $\pi_i(C)$ can not be a point for all $i=1,\cdots,m$.

The fibration $\pi_{1,2}:=\pi_1\times \pi_2:\mathbb{A}^m\rightarrow \mathbb{A}^2$ is persevered by $f.$ By our hypotheses, $\pi_{1,2}(C)$ is periodic. By replacing $f$ by some suitable positive iterate, we suppose that $\pi_{1,2}(C)$ is fixed by $f$. Observe that $\pi_{1,2}^{-1}(\pi_{1,2}(C))$ is a divisor on $\mathbb{A}^m$.

The fibration $\pi_{2,\cdots,m}:=\pi_2\times \cdots\times\pi_m:\mathbb{A}^m\rightarrow \mathbb{A}^{m-1}$ is persevered by $f.$ By the induction hypotheses, $\pi_{2,\cdots,m}(C)$ is periodic. By replacing $f$ by some suitable positive iterate, we suppose that $\pi_{2,\cdots m}(C)$ is fixed by $f$. Observer that $\pi_{2,\cdots,m}^{-1}(\pi_{2,\cdots,m}(C))$ is a surface.

If $\pi_{1,2}^{-1}(\pi_{1,2}(C))$ contains $\pi_{2,\cdots,m}^{-1}(\pi_{2,\cdots,m}(C))$, then we have $$\pi_{1,2}(\pi_{2,\cdots,m}^{-1}(\pi_{2,\cdots,m}(C)))\subseteq\pi_{1,2}(C).$$ Observe that $$\pi_{1,2}(\pi_{2,\cdots,m}^{-1}(\pi_{2,\cdots,m}(C)))=\pi_{1,2}(\pi_2^{-1}(\pi_2(C)))=\mathbb{A}^2.$$ Since $\pi_{1,2}(C)$ is a curve, this is a contradiction.

So we have that $\pi_{1,2}^{-1}(\pi_{1,2}(C))$ does not contain $\pi_{2,\cdots,m}^{-1}(\pi_{2,\cdots,m}(C))$, and then $D:=\pi_{1,2}^{-1}(\pi_{1,2}(C))\cap \pi_{2,\cdots,m}^{-1}(\pi_{2,\cdots,m}(C))$ is dimensional $1$ and it is fixed by $f.$ Since $C$ is an irreducible component of $D,$ we have that $C$ is periodic.
\endproof

\newpage

\part{The resonant case $\la_1^2=\la_2$}
In this part, we prove the main theorem in the case that $\la_1^2=\la_2.$ By \cite[Theorem C]{Favre2011}, we have either  $\deg (f^n)\asymp n\la_1^n$ or
$\deg (f^n)\asymp \la_1^n$. We will treat these two cases separately.



%
%

\section{The case $\la_1^2=\la_2$ and $\deg (f^n)\asymp n\la_1^n$}
In this section, denote by $k:=\overline{\mathbb{Q}}$, the field of algebraic numbers.

By \cite[Theorem C]{Favre2011}, we may suppose that $f$ takes form
$$f=(F(x),G(x,y))=(F(x),\sum_{i=0}^dA_i(x)y^i)$$ where $d=\deg F$ and $\deg A_d\geq 1.$ In this case $\la_1=d$ and $\la_2=d^2.$

The aim in this section is to show
\begin{thm}\label{thmlaoseqlatfibrdml}
If $\la_1^2=\la_2$ and $\deg (f^n)\asymp n\la_1^n$, then
the pair $(\mathbb{A}^2_k,f)$ satisfies the DML property.
\end{thm}

If $d=1$, then $f$ is birational. By \cite[Theorem A]{Xie2014}, Theorem \ref{thmlaoseqlatfibrdml} holds. So we may suppose that $d\geq 2$ in the rest of this section.


\subsection{Find an algebraically stable model}
Our aim is to make $f$ to be algebraically stable in a suitable Hirzebruch surface $\mathbb{F}_n$\index{$\mathbb{F}_n$} for some $n\geq 0.$

It is convenient to work with the presentation of these surfaces as a quotient by $(\mathbb{G}_m)^2,$ as in \cite{AntonioLaface}. By definition, the set of closed point $\mathbb{F}_n(k)$ is the quotient of $\mathbb{A}^4(k)\setminus (\{x_1=0 \text{ and } x_2=0\}\cup \{x_3=0 \text{ and } x_4=0\})$ by the equivalence relation generated by $$(x_1,x_2,x_3,x_4)\sim (\lambda x_1,\lambda x_2,\mu x_3,\mu/\lambda^n x_4)$$ for $\lambda,\mu \in k^*$. Denote by $[x_1,x_2,x_3,x_4]$ the equivalence class of $(x_1,x_2,x_3,x_4)$. We have a natural morphism
$\pi_n:\mathbb{F}_n\rightarrow \mathbb{P}^1$ given by $\pi_n([x_1,x_2,x_3,x_4])=[x_1:x_2]$ which makes $\mathbb{F}_n$ into be a locally trivial $\mathbb{P}^1$ fibration.

We shall look at the embedding $$i_n:\mathbb{A}^2\hookrightarrow \mathbb{F}_n:(x,y)\mapsto [x,1,y,1].$$ Then $\mathbb{F}_n\setminus \mathbb{A}^2$ is union of two lines: one is the fiber at infinity $F_{\infty}$ of $\pi_n$, and the other one is a section of $\pi_n$ which is denoted by $L_{\infty}$.
\bigskip

For each $n\geq \max\{\deg A_i|\,\, i=1,\cdots,d\}+1$, the map $f$ extends to a rational transformation $$f_n:[x_1,x_2,x_3,x_4]\mapsto [x_2^dF(x_1/x_2),x_2^d,x_2^{nd+\deg A_d}x_4^d(\sum_{i=1}^dA_i(x_1/x_2)(\frac{x_3}{x_2^nx_4})^i),x_2^{\deg A_d}x_4^d]$$ on $\mathbb{F}_n$.
We have $$I(f_n)=\{[1,0,0,1]\}\cup\{[r,1,1,0]|\,\, A_d(r)=0\}.$$ The unique curve contracted by $f_n$ is $F_{\infty}=\{x_2=0\}$ and its image is $f_n(F_{\infty})=[1,0,1,0]$. It implies the following:

\begin{pro}\label{propolycontra}\label{corpolyem}For $n\geq \max\{\deg A_i|\,\, i={1,\cdots,d}\}+1$, $f_n$ is algebraically stable on $\mathbb{F}_n$ and contracts the curve $F_{\infty}$ to the point $[1,0,1,0]$.
\end{pro}

\subsection{Dynamics on $V_{\infty}$}
Denote by $v_*$ the unique valuation in $V_1$ such that $f_{\d}(v_*)=v_*$ as in \cite[Proposition 5.1]{Favre2011}. Set $W(f):=\{v\in V_{\infty}|\,\, v\geq v_*\}$\index{$W(f)$}.

\begin{pro}\label{proproperoutside}For all $v\in V_{\infty}\setminus W(f)$, we have $d(f,v)\geq \la_1\alpha(v\wedge v_*)>0.$
\end{pro}
\proof Write $F(x)=a\prod_{i=1}^d(x-r_i)$ where $a>0$. Denote by $v_i$ the curve valuation defined by the unique branch of $\{x-r_i=0\}$ at infinity. Observe that $v_i\in W(f)$.

By definition, we have $d(f,v)=-\min\{v(F),v(G)\}$. It follows that $$d(f,v)\geq -v(F)=\sum_{i=1}^d\alpha(v_i\wedge v)=\la_1\alpha(v\wedge v_*)>0.$$
\endproof

By Proposition \ref{prologqbranch}, the function $\log|G|:v\mapsto -v(G)$ on $V_{\infty}$ can be written as $$\log|G|(v)=\sum_{i=1}^{l}m_i\alpha(v_i\wedge v)$$ where $v_i$'s are all curve valuations associated to the branches at infinity of $\{G(x,y)=0\}$ and $m_i\geq 1$. Suppose that $v_i\geq v_*$ for $i=1,\cdots,l_1$ and $v_i\not\geq v_*$ for $i=l_1+1,\cdots,l.$

There exists $v'\in V_{\infty}$ such that
\begin{points}
\item $v'<v^*$;
\item set $U:=\{v\in V_{\infty}|\,\,v'<v<v_*\}$, we have $f_{\d}$ maps $U$ strictly into itself and is order-preserving there;
\item $v_i\not\in U$ for all $i=1,\cdots,l$;
\item for all $v\in U$, we have $f^n_{\d}v\rightarrow v^*.$
\end{points}

Set $G_n:=G\circ f^{n-1}$ for $n\geq 1$ and write $\log |G_n|$ as $v\mapsto \sum_{i=1}^{l_n}m^n_i\alpha(v^n_i\wedge v).$ We may suppose that $v^n_i\geq v_*$ for $1\leq i\leq l'$ and $v^n_i\not\geq v_*$ for $l'+1\leq i\leq l_n$. Since $f_{\d}(U)\subseteq U$, we have $v_i^n\not\in U$ for all $i=1,\cdots,l_n.$  So $v_i^n\wedge v_*\leq v'$ for $i=l'+1,\cdots,l_n.$

Let $G_n^+$ be the function defined by $v\mapsto \sum_{i=1}^{l'}m^n_i\alpha(v^n_i\wedge v)$ and $G^-_n$ be the function defined by $v\mapsto \sum_{i=l'+1}^{l_n}m^n_i\alpha(v^n_i\wedge v)$. Then we have
$$\log |G_n|=G_n^++G_n^-.$$
Since $v_*(F^n)=\la_1\alpha(v_*)=0$, we have $$\la_1^n=d(f^n,v_*)=-v_*(G_n)=G_n^-(v_*).$$ Since $v_i^n\wedge v_*\leq v'$ for $i=l'+1,\cdots,l_n,$ we have
$G_n^-(v_*)=G_n^-(v')\geq \alpha(v') G_n^-(-\deg).$ It follows that $G_n^+(-\deg)\geq \deg G_n-\alpha(v')^{-1}\la_1^n.$
By \cite[Proposition 5.1]{Favre2011}, there exists $c'\geq 0$, such that $\deg G_n\geq c'n\la_1^n.$ Then we have the following
\begin{pro}\label{probounddeggnplus}There exists $c\geq 0$ such that $G_n^+(-\deg)\geq cn\la_1^n$ for all $n\geq 1.$
\end{pro}

For any $M\leq 0$, set $W_M:=\{v\in V_{\infty}|\,\,\alpha(v)\geq M\}$.
\begin{pro}\label{prowmsmwftou}There exists $N\geq 0$ such that $$f^N_{\d}(W_M\setminus W(f))\subseteq U.$$
\end{pro}

\proof Let $c$ be the number defined in Proposition \ref{probounddeggnplus}. Let $N$ be an integer at least $(c\alpha(v')^2)^{-1}(1-M)+1$. For any valuation $v\in W_M\setminus (W(f)\cup U)$, we have $$d(f^N,v)=-\min\{v(F^N),v(G_N)\}\geq -v(G_N)=G_N^+(v)+G_N^-(v).$$ Since $v^N_i\geq v_*$ for $1\leq i\leq l'$, we have $G_N^+(v)=\alpha(v_*\wedge v)G_N^+(-\deg)\geq c\alpha(v')N\la_1^N.$ On the other hand $G_N^-(v)\geq \alpha(v)G_N^-(-\deg)\geq \alpha(v')^{-1}M\la_1^N$. Then we have
$$d(f^N,v)\geq c\alpha(v')N\la_1^N+\alpha(v')^{-1}M\la_1^N>\la_1^N/\alpha(v').$$ Since $\la_1^N\alpha(v_*\wedge v)=(f^{*N}Z_{v^*}\cdot Z_v)=(Z_{v^*}\cdot f_*^NZ_v)=d(f,v)\alpha(v_*\wedge f_{\d}^Nv),$ then we have $\alpha(v_*\wedge f_{\d}^Nv)\leq \la_1^N/d(f,v)<\alpha(v').$ It follows that $f_{\d}^Nv\in U.$
\endproof

\subsection{Apply the Local dynamical Mordell-Lang Theorem}
\begin{pro}\label{prola1sela2localdml}Let $C$ be a curve in $\mathbb{A}^2_k$ admitting a branch at infinity which associates to a curve valuation in $U$ and let $p\in X$ be a closed point.
Then either $p$ is preperiodic or
the set $\{n\in \mathbb{Z}^+|\,\, f^n(p)\in C\}$ is finite.
\end{pro}
\proof
Fix an algebraically stable model $X:=\mathbb{F}_n$ for $n$ large enough, we see that $v_{L_{\infty}}=v_*$ and $v_{F_{\infty}}<v_*$. We may suppose that $v_{F_{\infty}}>v'$ and $v_{F_{\infty}}>v_*\wedge v_i$ for $i=1,\cdots, m$. Denote by $O:=[1,0,1,0]$ the intersection of $L_{\infty}$ and $F_{\infty}$.
We may check that $df|_O^2=0$, so $f$ is supperattracting at $O$. By replacing $C$ by $f^n(C)$ for $n$ large enough, we may assume $O\in C.$ Observe that the eigenvaluation in the local tree $V_{O}$ is a curve valuation.
Then by Theorem \ref{thmlodml}, we conclude our Proposition.
\endproof

\subsection{Curves with one place at infinity}
\begin{pro}\label{prolaoselatfoplace}If $C$ is a curve with one place at infinity and $p$ is a closed point in $\mathbb{A}^2$. If the set $\{n\geq 0|f^n(p)\in C\}$ is infinite, then either $p$ is preperodic or $C$ periodic.
\end{pro}

\proof Let $v_C$ be the curve valuation associated to the  unique branch at infinity of $C$. Pick an algebraically stable model $X:=\mathbb{F}_m$ for $m$ large enough. Either $C=\{x=a\}$ for some $a\in k$ or $C\cap F_{\infty}\neq \emptyset.$ In the forme case, our proposition trivially holds. Then we may suppose that $v_C\not\in W(f).$

By Proposition \ref{prola1sela2localdml}, we may suppose that $f^n_{\d}(v_C)\not\in U$ for all $n\geq 0.$ By Proposition \ref{prowmsmwftou}, there exists $N>0$ such that $W_{-1}\setminus W(f)\subseteq f_{\d}^{-N}(U)$. The boundary $\partial f_{\d}^{-N}(U)$ of $f^{-N}(U)$ is finite and for every point $v\in \partial f^{-N}(U)\setminus \{v_*\}$, we have $\alpha(v)\leq -1.$ Since for all $n\geq 0$, $f^n_{\d}(v_C)\not\in U$, there exists $v^n\in \partial f^{-N}(U)\setminus \{v_*\}$ such that $f^n_{\d}(v_C)\geq v^n.$ It follows that there exists $n_1>n_2\geq 0$ such that $v^{n_1}=v^{n_2}$.

If $f^{n_1}_{\d}(v_C)\neq f^{n_2}_{\d}(v_C),$ we have
\begin{align*}
\deg(f^{n_1}(C))\deg(f^{n_1}(C))=&(f^{n_1}(C)\cdot f^{n_2}(C))\nonumber\\
         \geq &\deg(f^{n_1}(C))\deg(f^{n_1}(C))(1-\alpha(f^{n_1}_{\d}(v_C)\wedge f^{n_2}_{\d}(v_C)))\\
         \geq & 2\deg(f^{n_1}(C))\deg(f^{n_1}(C)).
\end{align*}
It is impossible, so $f^{n_1}_{\d}(v_C)= f^{n_2}_{\d}(v_C),$ and then $C$ is preperiodic.

If $C$ is not periodic, there exist $m>0$ such that $f^m(C)$ is periodic. Then $\cup_{i=m}^{\infty}f^i(C)$ is a union of finitely many irreducible curves and  $f^n(p)\in \cup_{i=m}^{\infty}f^i(C)$ for $n\geq m.$ Since $C$ is not periodic, $C\cap \cup_{i=m}^{\infty}f^i(C)$ is finite. It follows that $p$ is preperiodic.
\endproof

\subsection{Curves with two places at infinity}
The aim of this section is to prove the following
\begin{pro}\label{prolaoselatftplace}If $C$ is a curve with two places at infinity and $p$ is a closed point in $\mathbb{A}^2$. If the set $\{n\geq 0|f^n(p)\in C\}$ is infinite, then either $p$ is preperodic or $C$ periodic.
\end{pro}
\proof[Proof of Proposition \ref{prolaoselatftplace}]
Let $C_1$ and $C_2$ be the two branches at infinity of $C$ and $v_{C_i}$ the curve valuation associated to $C_i$ for $i=1,2.$

Pick an algebraically stable model $X:=\mathbb{F}_m$ for $m$ large enough. Either $C=\{x=a\}$ for some $a\in k$ or $C\cap F_{\infty}\neq \emptyset.$ It follows that there exists $i=1,2$ such that $v_{C_i}\not\in W(f).$ So we may suppose that $v_{C_2}\not\in W(f).$

As in Theorem \ref{thmasequencecurves}, we have a sequence of curves $\{C^i\}_{i\in \mathbb{Z}}$ with at most two branches at infinity. By Proposition \ref{prolaoselatfoplace}, we may suppose that $C^i$ has exactly two branches at infinity for all $i\in \mathbb{Z}$. For $j=1,2,$ denote by $C^i_j$ the unique branch of $C^i$ such that $f^{-i}(C^i_j)=C_j$ for $i\leq 0$ and $f^i(C_j)$ for $i>0$.

\begin{lem}\label{lemlaoselattplessv}If $v_{C_i}\not\in W(f)$ for $i=1,2$, then Proposition \ref{prolaoselatftplace} holds.
\end{lem}

By Lemma \ref{lemlaoselattplessv}, we suppose that $v_{C_1}\in W(f)$, $v_{C_2}\not\in W(f)$.

\begin{lem}\label{lemlaoselattpolog}If $v_{C_1}\in W(f)$, $v_{C_2}\not\in W(f)$ and there are infinitely many $n\in \mathbb{Z}$ such that $(C^n_2\cdot l_{\infty})\geq (C^n_1\cdot l_{\infty})$, then Proposition \ref{prolaoselatftplace} holds.
\end{lem}

\begin{lem}\label{lemnonpassperiodic}Suppose that $v_{C_1}\in W(f)$, $v_{C_2}\not\in W(f)$, and $q=C_1\cap L_{\infty}$ satisfying one of the following
\begin{points}\item either $q$ is not periodic of $f|_{L_{\infty}}$;
\item or $q$ is $r$-periodic for some $r\geq 1$, $q\not\in I(f^r)$ and $f^r|_{L_{\infty}}$ is not ramified at $q$.
\end{points}
then Proposition \ref{prolaoselatftplace} holds.
\end{lem}

By replacing $f$ by a suitable positive iterate and $C$ by $C^j$ for some $j\leq 0$, we may suppose that there exists a point $q\in L_{\infty}$ satisfying
\begin{points}
\item $f|_{L_{\infty}}(q)=q$;
\item $q=C_1^j\cap L_{\infty}$ for all $j\leq 0$;
\item either $q\in I(f)$ or $f|_{L_{\infty}}$ is ramified at $q.$
\end{points}

\begin{lem}\label{lemlaoselatpfp}If there exists a point $q\in L_{\infty}\setminus I(f)$ such that $f(q)=q$ and $q\in C$, then Proposition \ref{prolaoselatftplace} holds.
\end{lem}

We may suppose that $f|_{L_{\infty}}(q)=q$ and $q\in I(f)$.
Then we conclude our Proposition by the following

\begin{lem}\label{lemlaoselatpoindfix}If there exists a point $q\in L_{\infty}\cap I(f)$ such that $f|_{L_{\infty}}(q)=q$ and $q\in C$, then either the set $\{n\geq 0|\,\,f^n(p)\in C\}$ is finite or $p$ is preperiodic.
\end{lem}
\endproof

\proof[Proof of Lemma \ref{lemlaoselattplessv}]By Proposition \ref{prola1sela2localdml}, we may suppose that $f^n_{\d}(v_{C_i})\not\in U$ for $i=1,2$ and all $n\geq 0.$ By Proposition \ref{prowmsmwftou}, there exists $N>0$ such that $W_{-7}\setminus W(f)\subseteq f_{\d}^{-N}(U)$. The boundary $\partial f_{\d}^{-N}(U)$ of $f^{-N}(U)$ is finite and for every point $v\in \partial f^{-N}(U)\setminus \{v_*\}$, we have $\alpha(v)\leq -7.$ Since for all $n\geq 0$, $f^n_{\d}(v_{C_i})\not\in U$ for $i=1,2$, there exists $v^n_i\in \partial f^{-N}(U)\setminus \{v_*\}$ such that $f^n_{\d}(v_{C_i})\geq v^n.$ Set $v^n=v^n_1$ if $(f^n(C_1)\cdot l_{\infty})\geq (f^n(C_2)\cdot l_{\infty})$ and $v^n=v^n_2$ otherwise.
There exists $n_1>n_2\geq 0$ such that $v^{n_1}=v^{n_2}$.

If $f^{n_1}(C)\neq f^{n_2}(C),$ we have

$$\deg(f^{n_1}(C))\deg(f^{n_1}(C))=(f^{n_1}(C)\cdot f^{n_2}(C))$$
$$\geq 4^{-1}\deg(f^{n_1}(C))\deg(f^{n_1}(C))(1-\alpha(v^{n_1}\wedge v^{n_2}))
$$$$\geq  2\deg(f^{n_1}(C))\deg(f^{n_1}(C)).$$

It is impossible, so $f^{n_1}_{\d}(C)= f^{n_2}_{\d}(C),$ and then $C$ is preperiodic.

If $C$ is not periodic, there exist $m>0$ such that $f^m(C)$ is periodic. Then $\cup_{i=m}^{\infty}f^i(C)$ is a union of finitely many irreducible curves and  $f^n(p)\in \cup_{i=m}^{\infty}f^i(C)$ for $n\geq m.$ Since $C$ is not periodic, $C\cap \cup_{i=m}^{\infty}f^i(C)$ is finite. It follows that $p$ is preperiodic.
\endproof

\proof[Proof of Lemma \ref{lemlaoselattpolog}]By Proposition \ref{prola1sela2localdml}, we may suppose that $f^n_{\d}(v_{C_2})\not\in U$ for all $n\geq 0.$ By Proposition \ref{prowmsmwftou}, there exists $N>0$ such that $W_{-7}\setminus W(f)\subseteq f_{\d}^{-N}(U)$. The boundary $\partial f_{\d}^{-N}(U)$ of $f^{-N}(U)$ is finite and for every point $v\in \partial f^{-N}(U)\setminus \{v_*\}$, we have $\alpha(v)\leq -7.$ Since for all $n\in \mathbb{Z}$, $v_{C^n_2}\not\in U$, there exists $v^n\in \partial f^{-N}(U)\setminus \{v_*\}$ such that $v_{C^n_2}\geq v^n.$ Set $A:=\{n\geq 0|\,\,(v_{C^n_2}\cdot l_{\infty})\geq (v_{C^n_1}\cdot l_{\infty})\}$. Since $A$ is infinite,
there exists different elements $n_1,n_2\in A$ such that $v^{n_1}=v^{n_2}$.

If $v_{C^{n_1}_2}\neq v_{C^{n_2}_2},$ we have

$$\deg(C^{n_1})\deg(C^{n_2})=(C^{n_1}\cdot C^{n_2})
       \geq (C^{n_1}_2\cdot C^{n_2}_2)$$$$
         \geq 4^{-1}\deg(C^{n_1})\deg(C^{n_2})(1-\alpha(v^{n_1}\wedge v^{n_2}))$$$$
         \geq 2\deg(C^{n_1})\deg(C^{n_2}).$$

It is impossible, so $v_{C^{n_1}_2}= v_{C^{n_2}_2},$ and then $C$ is preperiodic.
\endproof

\proof[Proof of Lemma \ref{lemnonpassperiodic}]By Lemma \ref{lemlaoselattpolog}, we may suppose that $(C^i_2\cdot l_{\infty})\leq (C^i_1\cdot l_{\infty})$ for $i\leq 0$. Set $D:=L_{\infty}+(m+1)F_{\infty}$ which is ample on $X=\mathbb{F}_m$.
Observe that for $i\leq 0$, $(C_1^i\cdot l_{\infty})=b_{L_{\infty}}(C_1^i\cdot L_{\infty})=(C_1^i\cdot L_{\infty})$ and $(C^i_2\cdot F_{\infty})\leq (C^i_2\cdot l_{\infty}).$ It follows that $$(C_2^i\cdot D)=(m+1)(C_2^i\cdot F_{\infty})=(m+1)(C^i_2\cdot l_{\infty})$$$$\leq (m+1)(C^i_1\cdot l_{\infty})\leq (m+1)(C_1^i\cdot L_{\infty})=(m+1)(C_1^i\cdot D).$$
Observe that $v_*=v_{L_{\infty}}$ is totally invariant.
Then Proposition \ref{prototinvanonexcepdml} shows that $(C^i\cdot D)$ is bounded for $i\leq 0.$ Then we conclude our Lemma by Proposition \ref{procurveboundeddegfibration}.
%
%
%
\endproof

\proof[Proof of Lemma \ref{lemlaoselatpfp}]By \cite{Favre2000}, $f$ locally conjugates to $(x,y)\mapsto (x^s,y^d)$  where $2\leq s\leq d$ with respect to any nontrivial norm $|\cdot|$ of $k$.

If $2\leq s\leq d-1,$ we conclude our lemma by Corollary \ref{cormonomial}.

Then we treat the case $m=d$.
We define a sequence of surfaces $\pi_k:X_k\rightarrow X$ by induction:
\begin{points}\item Set $X_0:=X$ and $\pi_0:=\id$.
\item Suppose that we have $X_0,\cdots, X_k$. If $C_1$ does not pass through any singular point of $\pi_k^*L_{\infty}$, we stop our progression.
\item If $C_1$ is passing through one singular point of $\pi_k^{-1}L_{\infty}$, let $X_{k+1}$ be the surface defined by blowup at this point in $X_k$.
\item Denote by $C_1$ the
strict transformation of $C_1$ in $X_k$.
Then return to (i).
\end{points}
This progression terminates in finitely many steps and we get surfaces $X_0,\cdots, X_l$ for $l\geq 0.$

It is easy to see that $f$ is regular on $\pi_l^{-1}(q).$ At any singular point of $\pi_l^{-1}(L_{\infty}),$ $f$ locally conjugates to $(x,y)\mapsto (x^d,y^d).$ Let $E$ be the unique exceptional curve of $\pi_l$ which intersects $C_1$ at one point. We have $E$ is totally invariant and $f|_E$ can be written as $z\rightarrow z^d.$ All the ramified points of $f|_{E}$ is singular in $\pi_l^{-1}(L_{\infty})$. Use the same method in the proof of Lemma \ref{lemnonpassperiodic}, we conclude our Lemma.
\endproof

\proof[Proof of Lemma \ref{lemlaoselatpoindfix}]
By changing coordinates, we suppose that $q=[0,1,1,0]\in \mathbb{F}_m$. Then $F(x)$ has form $x^sE(x)$ where $1\leq s\leq d$ and $E(0)\neq 0$ and $A_d(x)$ has form $x^rB(x)$ where $r\geq 1$. Let $\{n_i\}_{i\geq 1}$ be a increase sequence such that $f^{n_i}(p)\in C$. Set $f^n(p)=(x_n,y_n)$ and suppose that $p$ is not preperiodic.

\smallskip

Let $K$ be a number field such that $X$, $f$, $p$ and $C$ are all defined over $K$.
\begin{lem}\label{lemheigpair}There exists a place $v\in \mathcal{M}_K$ such that by replacing $n_i$ by a subsequence, we have $\log\max\{|y_{n_i}|_v,1,|x_{n_i}|_v^m\}-\log\max\{1,|x_{n_i}|_v^m\}\geq cd^{n_i}$ for some $c>0.$
\end{lem}
Then we suppose that $\log\max\{|y_{n_i}|_v,1,|x_{n_i}|_v^m\}-\log\max\{1,|x_{n_i}|_v^m\}\geq cd^{n_i}$ for some $c>0.$ Since $C\cap L_{\infty}$ is just one point $q$, we have that $(x_{n_i},y_{n_i})\rightarrow q$ as $i\rightarrow \infty$ with respect to $|\cdot|_v.$ It follows that $|x_{n_i}^m|_v\rightarrow 0$ and $|y_{n_i}|_v\rightarrow \infty$ as $i\rightarrow \infty.$ It follows that $\log(|y_{n_i}|_v)\geq cd^{n_i}$ for some $c>0.$ Since $(x_{n_i},y_{n_i})\in C$, and $C$ is not vertical, there exists $0<r'<1$ such that $|y_{n_i}|_v^{-1}\geq |x_{n_i}|^{r'}_v$ for $i$ large enough.

\smallskip

At first we treat the case $s=d.$ In a suitable coordinate, we have $F(x)=x^d$. By replacing $p$ by $f^n(p)$ for a suitable $n\geq 0,$ we suppose that $|x_0|_v<1.$ We have $|x_{n}|_v=|x_0|_v^{d^n}$ and $|y_{n+1}|_v\leq a|x_n|_v^r|y_n|_v^d+b|y_n|_v^{d-1}$ for some $a,b>0.$
\begin{lem}\label{lemmbigsed}There exists $N\geq 0$, such that for all $n\geq N$, $a|x_n|^r|y_n|\geq b.$
\end{lem}
By replacing $p$ by $f^N(p)$, we suppose that $N=0.$ Then we have $|y_{n+1}|_v\leq a|x_n|_v^r|y_n|_v^d+b|y_n|_v^{d-1}\leq 2a|x_n|_v^r|y_n|_v^d$ for all $n\geq 0.$ Set $Y_n:=\log(|y_n|_v)$, $A:=\log (2a)$ and $U:=\log (x_0)$, we have $$Y_{n+1}\leq A+rd^nU+dY_{n}$$ for all $n\geq 0.$ Then we have $Y_{n+1}/d^{n+1}-Y_{n}/d^n\leq A/d^{n+1}+rU/d$ for $n\geq 0.$ It follows that $$Y_n/d^n\leq \sum_{i=1}^{\infty}|A|/d^i+nrU/d+Y_0=|A|/(d-1)+Y_0+nrU/d$$ for $n\leq 0.$ Since $U<0$, we have $\log(|y_{n}|_v)/d^n\rightarrow -\infty$. It contradicts to the fact that $\log(|y_{n_i}|_v)\geq cd^{n_i}$ for some $c>0.$

\smallskip

Then we treat the case $2\leq s\leq d-1.$ Since $0$ is an attracting fixed point of $F$, we have $|x_n|_v\rightarrow 0$ as $n\rightarrow \infty.$ We may suppose that for all $n\geq 0$, we have $|x_n|_v<1$ and  $a'|x_n|_v^r|y_n|_v^d-b|y_n|_v^{d-1}\leq |y_{n+1}|_v\leq a|x_n|_v^r|y_n|_v^d+b|y_n|_v^{d-1}$ for some $a>a'>0$ and $b>0$. There exists $e>0$ such that $|x_{n+1}|_v\geq e|x_n|_v^s$ for $n\geq 0$. There exists $c_1>c_2>0$ and $u>0$ such that $c_1u^{s^n}>|x_n|_v>c_2u^{s^n}$ for all $n\geq 0.$

\begin{lem}\label{lemnlargerfirstterm}There exists $N\geq 0$, such that for all $n\geq N$, $a'|x_n|^r|y_n|\geq 2b.$
\end{lem}
By replacing $p$ by $f^N(p)$, we may suppose that $N=0.$ Then we have $$|y_{n+1}|_v\geq a'|x_n|_v^r|y_n|_v^d-b|y_n|_v^{d-1}\geq a'/2|x_n|_v^r|y_n|_v^d\geq a'c_2/2u^{rs^n}|y_n|_v^d$$ for $n\geq 0.$ Set $Y_n:=\log (|y_{n+1}|_v)$, $A:=\log (a'c_2/2)$ and $U:=\log u$. We have $$Y_{n+1}\geq dY_n+s^nU+A$$ for $n\geq 0.$ It follows that $$Y_n/d^n\geq Y_0+\sum_{i=0}^{\infty}(s/d)^iU/d-\sum_{i=0}^{\infty}|A|/d^{i+1}=Y_0+U/(d-s)-|A|/(d-1).$$ Since $Y_{n_i}\geq d^{n_i}+\log(c)$ and $d>s$, by replacing $p$ by $f^{n_i}(p)$ and $u$ by $u^{s^{n_i}}$ for some $i\geq 1$, we may suppose that $Y_0+U/(d-s)-|A|/(d-1)>0$. Then there exists $B>1$ such that $|y_n|_v\geq B^{d^n}.$ Since $|x_n|_v>c_2u^{s^n}$, for any $r'>0$, $$|y_n|_v^{-1}\leq B^{-d^n}<c_2^{r'}u^{r's^n}<|x_n|^{r'}_v$$ for $n$ large enough. It contradicts that fact that there exists $0<r'<1$ such that $|y_{n_i}|_v^{-1}\geq |x_{n_i}|^{r'}_v$ for $i$ large enough.

Finally we treat the case $s=1.$ If there exists $i\leq -1$ such that the center $q_i$ of $C^i_1$ is not $q$, then $q_i$ is not a periodic point of $f|_{L_{\infty}}$. Then we conclude our proposition by Lemma \ref{lemnonpassperiodic}. So we may suppose that the center of $C^i_1$ is $q$ for all $i\in \mathbb{Z}$. Since $s=1$, for any point of $C^i$ near $q$ has at most $d$ preimages near $q$. It follows that $\deg (f|_{C^{i-1}})\leq d.$ Then we have $$(C_1^i\cdot L_{\infty})=1/d(C_1^i\cdot f^* L_{\infty})=\deg(f_{C_1^i})/d(C^{i+1}_1\cdot L_{\infty})\leq (C_1^{i+1}\cdot L_{\infty})$$ for $i\leq -1.$ Then we conclude our Proposition by the same argument in the proof of Lemma \ref{lemnonpassperiodic}.
\endproof

\proof[Proof of Lemma \ref{lemheigpair}] Let $h_1:C(K)\rightarrow \mathbb{R}$ be the function defined by $(x,y)\mapsto \sum_{v\in \mathcal{M}_K}\log\max\{|x|_v,1\}$ and  $h_2:C(K)\rightarrow \mathbb{R}$ be the function defined by $(x,y)\mapsto \sum_{v\in \mathcal{M}_K}(\log\max\{|y|_v,1,|x|_v^m\}-\log\max\{1,|x|_v^m\})$. It follows that $h_1$ is a Weil height function with respect to the divisor $C\cdot {F_{\infty}}$ and $h_2$ is a Weil height function with respect to the divisor $C\cdot{L_{\infty}}.$ If $x_0$ is preperiodic, since $p$ is not preperiodic, we have $C=\{x=x_0\}$. It contradicts the fact that $C$ has two place at infinity.

By Lemma \ref{lemhoverh}, we have $h_2(f^{n_i}(p))\geq c_1h_1(f^{n_i}(p))\geq c_1c_2d^{n_i}$ where $c_1,c_2>0.$ There exists a finite set $S$ of place, such that for all $v\in \mathcal{M}_K\setminus S$, we have $|x_n|_v\leq 1$ and $|y_n|_v\leq 1$ for all $n\geq 0.$ Then we have $\sum_{v\in S}(\log\max\{|y_{n_i}|_v,1,|x_{n_i}|_v^m\}-\log\max\{1,|x_{n_i}|_v^m\})=h_2(f^{n_i}(p))\geq c_1c_2d^{n_i}.$ it follows that there exists $v\in S$ such that there exists infinitely many $i$ such that $\log\max\{|y_{n_i}|_v,1,|x_{n_i}|_v^m\}-\log\max\{1,|x_{n_i}|_v^m\}\geq  (\#S)^{-1}c_1c_2d^{n_i}.$
\endproof

\proof[Proof of Lemma \ref{lemmbigsed}] Set $u:=|x_0|_v<1.$ There exists $N\geq 0$ such that $$u^{rd^n}\leq \frac{a^{d-2}}{2b^{d-1}}$$ for all $n\geq N.$ If there exists $n\geq N$ such that $a|x_n|_v^r|y_n|_v\leq b$, then we have $|y_{n+1}|_v\leq a|x_n|_v^r|y_n|_v^d+b|y_n|_v^{d-1}\leq 2b|y_n|_v^{d-1}$. It follows that

$$a|x_{n+1}|_v^r|y_{n+1}|_v\leq 2abu^{rd^{n+1}}|y_{n}|_v^{d-1}
         = 2abu^{rd^n}(u^{rd^n}|y_n|_v)^{d-1}$$$$
         \leq  2abu^{rd^n}(b/a)^{d-1}
          \leq b.$$

It follows that there exists $N'\geq N$ such that $a|x_n|_v^r|y_n|_v\leq b$ for all $n\geq N'.$ Replacing $p$ by $f^{N'}(p)$, we may suppose that $N'=0.$ Then we have $$|y_{n+1}|_v\leq a|x_n|_v^r|y_n|_v^d+b|y_n|_v^{d-1}\leq 2b|y_n|_v^{d-1}$$ for $n\geq 0.$ It follows that there exists $c_1>0$ such that $|y_n|_v\leq c_1^{(d-1)^n}$ for all $n\geq 0.$ It contradicts the fact that $\log(|y_{n_i}|_v)\geq cd^{n_i}$ for some $c>0.$
\endproof
\proof[Proof of Lemma \ref{lemnlargerfirstterm}]
 Set $M:=\max\{(a^{'-d+2}e^r/2)^{-1/(d-1)},2b\}$. Since $\log(|y_{n_i}|_v)\geq cd^{n_i}$ for some $c>0$, we have $a'|x_{n_i}|^r|y_{n_i}|\rightarrow \infty$ as $i\rightarrow \infty$. So there exists $n\geq 0$ such that $a'|x_n|^r|y_n|\geq M\geq 2b.$ By induction, we only have to show $a'|x_{n+1}|^r|y_{n+1}|\geq M$. We have

$$|y_{n+1}|_v\geq a'|x_n|_v^r|y_n|_v^d-b|y_n|_v^{d-1}
         \geq a'/2|x_n|_v^r|y_n|_v^d.$$

It follows that
$$
a'|x_{n+1}|^r|y_{n+1}|\geq a^{'2}/2|x_{n+1}|^r|x_n|_v^r|y_n|_v^d$$$$
         \geq a^{'2}e^r/2|x_n|_v^{s+1}r|y_n|_v^d
        \geq a^{'2}e^r/2|x_n|_v^dr|y_n|_v^d$$$$
        \geq a^{'2}e^r/2(M/a')^d
        \geq M.
     $$
\endproof

\section{The case $\lambda_1^2=\la_2$ and $\deg (f^n)\asymp \la_1^n$}In this section, denote by $k:=\overline{\mathbb{Q}}$ the field of algebraic numbers.

The aim of this section is to prove the following
\begin{thm}\label{thmlaoqelat}Let $f:\mathbb{A}^2\rightarrow \mathbb{A}^2$ be a polynomial endomorphism define over $k.$ We suppose that $\la_1(f)^2=\la_2(f)$, and $\deg(f^n)/\la_1(f)^n$ is bounded. Let $C$ be a curve in $\mathbb{A}^2$ and $p$ be a closed point in $\mathbb{A}^2(k)$. Then if the set $\{n\in \mathbb{N}|\,\,f^n(p)\in C\}$ is infinite, we have that either $p$ is preperiodic for $f$ or $C$ is periodic for $f$.
\end{thm}

If $\la_1(f)=1$, then $f$ is birational. We conclude Theorem \ref{thmlaoqelat} by \cite{Xie2014}.

In the rest of this section, suppose that $\la_1(f)>1$.

\begin{defi}We define \index{$\mathcal{T}_f$}$$\mathcal{T}_f:=\{v\in V_1|\,\, f_{\d}(v)=v \}.$$
\end{defi}
Recall that $V_1$\index{$V_1$} is the set of valuations $v\in V_{\infty}$ satisfying $\alpha(v)\geq 0$ and $A(v)\leq 0$.
The boundary of $V_{1}$ is the set of  valuations $v\in V_1$ satisfying $\alpha(v)=0$ or $A(v)=0$.

The following proposition is come from \cite[Section 5]{Favre2011}.
\begin{pro}\label{prosecfivefa2011}We have
\begin{points}
\item $f$ is proper;
\item for every valuation $v\in \mathcal{T}_f$, we have that $v$ is totally invariant under $f_{\d}$, $f^*Z_v=\la_1Z_v$ and $d(f,v)=\la_1.$
\item by replacing $f$ by $f^2$, we may assume that $\mathcal{T}_{f^n}=\mathcal{T}_f$ for $n\geq 1$ and
either $T_f$ consists of a single divisorial valuation $v_*\in V_1$ with $\alpha(v_*)>0$
or $T_f$ is a closed segment in $V_1$ whose endpoints are divisorial valuations.
\end{points}
\end{pro}
In the rest of this section, we suppose that $\mathcal{T}_{f^n}=\mathcal{T}_f$ for $n\geq 1$.

\medskip

At first, we need a result of the dynamics on $V_{\infty}.$
For any divisorial valuation $v_E\in \mathcal{T}_f$, denote by $\mathbf{f}:\Tan_{v_E}\to \Tan_{v_E}$ the tangent map. Let $\vec{v_E}$ be a direction at $v_E$ fixed by $\mathbf{f}$. For any valuation $w\in U(\vec{v_E})$ we define $\vec{w}$ to be the direction at $w$ determined by $v_E$ and $U_{v_E,w}$ to be the open set $U(\vec{v_E})\cap U(\vec{w})$.

Then we have the following
\begin{pro}\label{provefixnontotinvathenatr}If $\alpha(v_E)>0$ and $\vec{v_E}$ is not totally invariant under $\mathbf{f}$, then there exists $w\in U(\vec{v_E})$ such that
\begin{points}
\item
$f_{\d}(U_{v_E,w})\subseteq U_{v_E,w};$
\item
for all $v\in U_{v_E,w}$, we have $f^n_{\d}(v)\to v_E$ for $v\to \infty$;
\item for any $M\leq 0$, there exists $N\geq 0$ such that $U(\vec{v_E})\cap \{v\in V_{\infty}|,\, \alpha(v)\geq M\}\subseteq f_{\d}^{-N}(U_{v_E,w}).$
\end{points}
\end{pro}

\proof[Proof of Proposition \ref{provefixnontotinvathenatr}]
By the proof of \cite[Theorem C]{Favre2011}, there exists a projective compactification $X$ of $\mathbb{A}^2_k$  with at most a quotient singularities such that the unique irreducible component of $X\setminus \mathbb{A}^2_k$ is $E$ and $f$ extends to an endomorphism on $X.$ The direction $\vec{v_E}$ determines a point $q\in E$ which is fixed by $f.$ Denote by $m$ the local degree of map $f|_E$ at $q$. Since
$\vec{v_E}$ is not totally invariant under $\mathbf{f}$, $q$ is not totally invariant and then $m< \la_1.$

By embedding $k$ in $\mathbb{C}$, we may view $X$ as a complex variety.
There exists a map $\pi:(\mathbb{C}^2,0)\rightarrow (\mathbb{C}^2,0)/G=(X,q)$ where $\pi$ is the quotient map and $G$ is the cyclic group generated $g:(x,y)\mapsto (e^{\frac{2\pi i}{l}}x,e^{\frac{2\pi is}{l}}y)$, $s,l\in \mathbb{Z}^+$ and $(s,l)=1$. Since $\mathbb{C}^2\setminus \{0\}$ is simply connected,
$f$ lifts to an endomorphism $F:(\mathbb{C}^2,0)\to (\mathbb{C}^2,0).$ Denote by $V_0$ the local valuative tree of $(\mathbb{C}^2,0)$.
\begin{lem}\label{lempipullbackeirre}The pullback $\pi^{-1}E$ is irreducible in $(\mathbb{C}^2,0)$. There exists a valuation $w_0<v^q_{\pi^{-1}(E)}$, such that
$F_{\d}(\{v>w_0\})\subseteq \{v>w_0\}$ and
for all $v\in\{v>w_0\},$ we have $F^n_*(v)\rightarrow v^q_{\pi^{-1}(E)}$ as $n\to \infty$. Further for any $v\in V_0$ satisfying $\alpha^0(v)<\infty$, we have $F^n_*(v)\rightarrow v^q_{\pi^{-1}(E)}$ as $n\to \infty$.
\end{lem}

Set $E':=\pi^{-1}(E).$ Since $\pi(g(E'))=E$, we have $g(E')=E'$. We may suppose that $E'=\{y=0\}$. Denote by $E''$ the curve defined by $x=0$. Let $P_1$ (resp. $P_2$) be germ of analytic function on $(X,q)$ defined by $P_1(\pi((x,y)))=x^l$ (resp. $P_2(\pi((x,y)))=y^l$). Observe that these functions are well defined.

There exists a map $\pi_{\d}:V_0\to \overline{U(\overrightarrow{v_E})}$ defined by $\pi_{\d}(P):=h(\pi,v)v(\pi^*P)$ for any polynomial $P\in \mathbb{C}[x,y]$ where $h(\pi,v)=-v(\pi^*L)$ where $L$ is a general linear form in $\mathbb{C}[x,y].$

The group $G$ acts on $V_0$ and the map $\pi_{\d}$ is the quotient map $V_0\rightarrow V_0/G\simeq \overline{U(\overrightarrow{v_E})}$.
Observe that for any $v\in V_0\setminus ([\ord_0,v^0_{E'}]\cup [\ord_0,v^0_{E''}])$ the orbit $Gv$ has $l$ elements and for any $v\in [\ord_0,v^0_{E'}]\cup [\ord_0,v^0_{E''}]$ the orbit $Gv$ has $1$ elements. Pick $w_0$ as in Lemma \ref{lempipullbackeirre} and $w:=\pi_{\d}(w_0)$, then $\pi_{\d}(\{v\in V_{0}|\,\,v>w_0\}\setminus \{v_{E'}\})=U_{v_E,w}$ which satisfies (i) and (ii) in our proposition.

We claim the following
\begin{lem}\label{lemcontrolalpzebyalp}For any $M\leq 1$, there exists a real number $C_M>0$ such that for all $v\in V_{0}\setminus \{v\in V_{0}|\,\, v>w_0\}$ satisfying $\alpha(\pi_{\d}(v))\geq M$ we have $\alpha^0(v)\leq C_M.$
\end{lem}
For any $M\leq 1$, we have $\pi_{\d}^{-1}(\{v\in U(q)|\,\,\alpha(v)\geq M\}\setminus U_{v_E,w})\subseteq \{v\in V_0|\,\,\alpha^0(v)\leq C_M\}.$
By Lemma \ref{lempipullbackeirre} and the compactness of $\{v\in V_0|\,\,\alpha^0(v)\leq C_M\}$, there exists $N\geq 0$ such that $\{v\in V_0|\,\,\alpha^0(v)\leq C_M\}\subseteq F^{-N}_{\d}(\{v>w_0\}).$ Since $\pi_{\d}$ is surjective, $f_{\d}^N(\{v\in U(q)|\,\,\alpha(v)\geq M\}\setminus\{v\in U_{v_E,w}\})\subseteq U_{v_E,w}.$ Since $f_{\d}(U_{v_E,w})\subseteq U_{v_E,w}$, we have that $f_{\d}^N(\{v\in U(q)|\,\,\alpha(v)\geq M\})\subseteq U_{v_E,w}$ which concludes (iii).
\endproof
\proof[Proof of Lemma \ref{lempipullbackeirre}]By Lemma \ref{lemlocalallessinftytendtocu}, we only have to show that $\pi^{-1}(E)$ is irreducible.
Let $E'$ be an irreducible component of $\pi^{-1}(E)$. Since $f^*E=\la_1E$, we have $F^*(\pi^*E)=\la_1\pi^*E$. It follows that $F^*E'$ is an irreducible component of $\pi^*E$. By replacing $f$ by a suitable positive iterate, we may suppose that $F^*E'=\la_1E'$. Since $\pi|_E'$ is finite, we have $F_*E'=mE'$ locally. Pick $v$ a valuation in $V_0$ satisfying $\alpha^0(v)<\infty$, by Lemma \ref{lemlocalallessinftytendtocu}, we have $F^n_{\d}v\to v^q_{E'}$ as $n\to \infty.$ If $E''$ is another irreducible component of $\pi^{-1}(E)$, the same argument shows that $F^n_{\d}v\to v^q_{E''}$ as $n\to \infty.$ It follows that $E'=E''$ and then $\pi^{-1}(E)$ is irreducible.
\endproof
\proof[Proof of Lemma \ref{lemcontrolalpzebyalp}]
There exists $T\geq 1$, such that for all $v\in V_{0}\setminus \{v\in V_{0}|\,\, v>w_0\}$, $v(y)\leq T.$

Observe that $\pi^*L=y^{-b_E}U(x,y)$, where $U$ is a unit in $\mathbb{C}[[x,y]].$
For any divisorial valuation $v_{D'}^0$, there exist are birational model $Y_0\rightarrow (\mathbb{C}^2,0)$ and $Y\rightarrow(X,q)$ such that $D'$ is an exceptional divisor in $Y'$, the rational map $\pi':Y_0\rightarrow Y$ induced by $\pi$ is a morphism, and $\pi'|_{D'}$ is finite. Denote by $e_{D'}$ the degree of $\pi'|_{D'}$. Set $r_{D'}:=\ord_{D'}(\pi^{'*}D).$
Observe that $r_{D'}\times\#(Gv_{D'})\times e_{D'}=l.$ Set $D:=\pi'(D').$

It follows that
$$-b_E\ord_{D'}(y)=\ord_{D'}(\pi^*L)=r_{D'}\ord_D(L)=-r_{D'}b_D.$$ Then we have $\ord_{D'}(y)=r_{D'}b_D/b_E.$
If $v_{D'}\in V_{0}\setminus \{v\in V_{0}|\,\, v>w_0\}$, we have $T\geq v_{D'}(y)=(b^0_{D'})^{-1}\ord_{D'}(y)=(b^0_{D'})^{-1}r_{D'}b_D/b_E.$
It follows that $b_D/b^0_{D'}\leq Tb_E/r_{D'}.$

\smallskip

Since $g$ is an automorphism on $(\mathbb{C}^2,0)$, we have $c(g,v)=1$ for all $v\in V_0$ and for any valuations $v_1,v_2\in V_0$, we have $\alpha^0(v_1\wedge v_2)=\alpha^0(g_{\d}(v_1)\wedge g_{\d}(v_2)).$ In particular, $\alpha^0(v)=\alpha^0(g_{\d}(v))$ for all $v\in V_0.$

For any $v\in V_0$, we have $\alpha^0(v\wedge g_{\d}(v))=\alpha^0(g_{\d}(v)\wedge g^2_{\d}(v))$. It follows that $v\wedge g_{\d}(v)=g_{\d}(v)\wedge g_{\d}^2(v)=v\wedge g_{\d}(v)\wedge g_{\d}^2(v)$. The same argument for $g^i$, $i=1,\cdots,l-1$, we have $v\wedge g_{\d}^i(v)=\wedge_{i=0}^{l-1}g^i_{\d}(v)$ for all $v\in V_0.$

We suppose first that $\overrightarrow{v_E}$ is not defined by $-\deg.$ Let $v_{D}$ be a divisorial valuation in $U(\overrightarrow{v_E})\setminus U_{v_E,w}$.
There exist are birational model $Y_0\rightarrow (\mathbb{C}^2,0)$ and $Y\rightarrow(X,q)$ such that $D$ is an exceptional divisor in $Y$, the rational map $\pi':Y_0\rightarrow Y$ induced by $\pi$ is a morphism and $g$ is lift to an endomorphism of $Y'$. Denote by $D'$ an irreducible component of $\pi'^{*}D$. Observe that $v^0_{D'}\in V_{0}\setminus \{v\in V_{0}|\,\, v>w_0\}.$ Set $H:=b_D(Z_{v_D}-Z_{v_E}).$ For any exceptional divisor $F$ of $Y\rightarrow (X,q)$, we have $(H\cdot F)=\delta_{F,D}$ and the support of $H$ are contained in the exceptional set of $Y\rightarrow (X,q)$. Then the support of $\pi'^*H$ is contained in the exceptional set of $Y'\rightarrow (\mathbb{C}^2,0)$ and for any irreducible exceptional divisor $F'$ of $Y'\rightarrow (\mathbb{C}^2,0)$, we have
$(\pi'^*H\cdot F')=(H\cdot \pi'_*F')=e_{F'}(H\cdot \pi'(F))=e_{F'}\delta_{\pi'(F'),D}.$ When $\pi'(F')=D$, we have $F'=g^i(D')$ for some $i=1,\cdots,l.$
It follows that $e_{D'}^{-1}\pi'^*H=(b^0_{D'})^2(\sum_{v\in G_{v_{D'}}}Z^0_v).$ It follows that
$$\left((\sum_{v\in G_{v_{D'}}}Z^0_v)\cdot (\sum_{v\in G_{v_{D'}}}Z^0_v)\right)=(b^0_{D'})^{-2}e_{D'}^{-2}(\pi'^*H\cdot \pi'^*H)$$$$=(b^0_{D'})^{-2}e_{D'}^{-2}l(H\cdot H)=(b_D/b^0_{D'})^{2}e_{D'}^{-2}l\left((Z_{v_D}-Z_{v_E})\cdot (Z_{v_D}-Z_{v_E})\right)$$$$=(b_D/b^0_{D'})^{2}e_{D'}^{-2}l(\alpha(v_D)-\alpha(v_E)).$$
Since for any $v,w\in V_0$, we have $(Z^0_{v}\cdot Z^0_w)=-\alpha(v\wedge w)<0$, we have $$\left((\sum_{v\in G_{v_{D'}}}Z^0_v)\cdot (\sum_{v\in G_{v_{D'}}}Z^0_v)\right)\leq \sum_{v\in Gv_{D'}}(Z^0_{v}\cdot Z^0_{v})=-\#(Gv_{D'})\alpha(v_{D'}).$$ Then we have
$$\alpha(v_{D'})\leq (Tb_E/r_{D'})^{2}e_{D'}^{-2}(\#(Gv_{D'}))^{-1}l(\alpha(v_E)-\alpha(v_D))\leq (Tb_E)^{2}(\alpha(v_E)-\alpha(v_D)).$$
Since divisorial valuation is dense in $V_{0}\setminus \{v\in V_{0}|\,\, v>w_0\}$, we have
$$\alpha^0(v)\leq (Tb_E)^{2}(\alpha(v_E)-\alpha(\pi_{\d}(v)))$$ for all $v\in V_{0}\setminus \{v\in V_{0}|\,\, v>w_0\}.$
If $\alpha(\pi_{\d}(v))\geq M$, we have $$\alpha^0(v)\leq (Tb_E)^{2}(\alpha(v_E)-\alpha(\pi_{\d}(v)))\leq (Tb_E)^{2}(\alpha(v_E)-M).$$
Then $C_M:=(Tb_E)^{2}(\alpha(v_E)-M)$ is what we require.

Now we suppose that $\overrightarrow{v_E}$ is defined by $-\deg.$ Let $v_{D}$ be a divisorial valuation in $U(\overrightarrow{v_E})\setminus U_{v_E,w}$.
There exist are birational model $Y_0\rightarrow (\mathbb{C}^2,0)$ and $Y\rightarrow(X,q)$ such that $D$ is an exceptional divisor in $Y$, the rational map $\pi':Y_0\rightarrow Y$ induced by $\pi$ is a morphism and $g$ is lift to an endomorphism of $Y'$. Denote by $D'$ an irreducible component of $\pi'^{*}D$.
Observe that $v^0_{D'}\in V_{0}\setminus \{v\in V_{0}|\,\, v>w_0\}.$ Set $H:=b_D(Z_{v_D}-(\alpha(v_D\wedge v_E)/\alpha(v_E))Z_{v_E}).$ For any exceptional divisor $F$ of $Y\rightarrow (X,q)$, we have $(H\cdot F)=\delta_{F,D}$ and the support of $H$ are contained in the exceptional set of $Y\rightarrow (X,q)$.
The same argument in the previous paragraph shows that $\pi'^*H=e_Db_{D'}^0(\sum_{v\in Gv_{D'}}Z_{v}).$ It follows that
$$\left((\sum_{v\in G_{v_{D'}}}Z^0_v)\cdot (\sum_{v\in G_{v_{D'}}}Z^0_v)\right)=(b^0_{D'})^{-2}e_{D'}^{-2}(\pi'^*H\cdot \pi'^*H)=(b^0_{D'})^{-2}e_{D'}^{-2}l(H\cdot H)$$$$=(b_D/b^0_{D'})^{2}e_{D'}^{-2}l\left((Z_{v_D}-(\alpha(v_D\wedge v_E)/\alpha(v_E))Z_{v_E})\cdot (Z_{v_D}-(\alpha(v_D\wedge v_E)/\alpha(v_E))Z_{v_E})\right)$$
$$=(b_D/b^0_{D'})^{2}e_{D'}^{-2}l\alpha(v_E)^{-1}(\alpha(v_D)\alpha(v_E)-\alpha(v_E\wedge v_D)^2).$$ It follows that
$$\alpha(v_{D'})\leq -(\#(Gv_{D'}))^{-1}\left((\sum_{v\in G_{v_{D'}}}Z^0_v)\cdot (\sum_{v\in G_{v_{D'}}}Z^0_v)\right)$$$$
= (Tb_E/r_{D'})^{2}e_{D'}^{-2}(\#(Gv_{D'}))^{-1}l\alpha(v_E)^{-1}(\alpha(v_E\wedge v_D)^2-\alpha(v_D)\alpha(v_E))$$$$\leq (Tb_E)^{2}\alpha(v_E)^{-1}(\alpha(v_E\wedge v_D)^2-\alpha(v_D)\alpha(v_E)).$$
Since divisorial valuation is dense in $V_{0}\setminus \{v\in V_{0}|\,\, v>w_0\}$, we have
$$\alpha^0(v)\leq (Tb_E)^{2}\alpha(v_E)^{-1}(\alpha(v_E\wedge \pi_{\d}(v))^2-\alpha(\pi_{\d}(v))\alpha(v_E))$$ for all $v\in V_{0}\setminus \{v\in V_{0}|\,\, v>w_0\}.$
If $\alpha(\pi_{\d}(v))\geq M$, we have $$\alpha^0(v)\leq (Tb_E)^{2}\alpha(v_E)^{-1}(\alpha(v_E\wedge \pi_{\d}(v))^2-\alpha(\pi_{\d}(v))\alpha(v_E))
$$$$\leq (Tb_E)^{2}\alpha(v_E)^{-1}(1+(M+1)\alpha(v_E))\leq (Tb_E)^{2}\alpha(v_E)^{-1}(M+2)
.$$
Then $C_M:=(Tb_E)^{2}\alpha(v_E)^{-1}(M+2)$ is what we require.
\endproof

Let $C_i$'s be all branches of $C$ at infinity.
\begin{pro}\label{proallbranchnotend}If for every branch $C_i$ of $C$ at infinity, we have $\alpha(r_{\mathcal{T}_{f}}(v_{C_i}))>0$, then Theorem \ref{thmlaoqelat} holds.

In particular, if for all $v\in \mathcal{T}_f$, we have $\alpha(v)>0$, then Theorem \ref{thmlaoqelat} holds.
\end{pro}
\proof[Proof of Proposition \ref{proallbranchnotend}]Let $s\in \{1,2\}$ be the number of places of $C$ at infinite. Set $v_i:=r_{\mathcal{T}_{f}}(v_{C_i})$ and let $\overrightarrow{v}_i$ be the tangent vector at $v_i$ presented by the segment $[v_i,v_{C_i}].$ Let $\mathbf{f}:\Tan_{v_i}\rightarrow \Tan_{v_i}$  be the tangent map at $v_i$ induced by $f$. By (iii) of Proposition \ref{prosecfivefa2011}, $v_i$ is divisorial.
There exists a projective smooth compactification $X$ of $\mathbb{A}^2$ such that for every $v_i$, there exists an exceptional component $E_i$ in $X\setminus \mathbb{A}^2_k$ satisfying $v_{E_i}=v_i$.

Let $G$ be the set of indexes $i$ such that $\overrightarrow{v_i}$ is not periodic under the tangent map $\mathbf{f}$.  By replacing $f$ by some positive iterate, we may suppose that  $\overrightarrow{v_i}$ is fixed by $\mathbf{f}$ for all $i\not\in G$.
By Theorem \ref{thmasequencecurves}
there exists a sequence of curves $\{C^j\}_{j\leq 0}$ with $s$ places at infinity such that
\begin{points}
\item $C^0=C$;
\item $f(C^j)=C^{j+1}$;
\item for all $j\leq -1$, the set $\{n\geq 0|f^n(p)\in C^j\}$ is infinite.
\end{points}
By replacing $C$ by some $C^j$, we may suppose that for all $j\leq 0$, $C^j$ has exact $s$ branches at infinity.
Let $C^j_i$'s be branches of $C^j$, we may suppose that $f(C^j_i)=C_i^{j+1}$ for $j\leq -1$ and $1\leq i\leq s.$ Since $v_i$ is totally invariant under $f_{\d}$, we have $r_{\mathcal{T}_{f}}(v_{C^j_i})=v_i$. Denote by $q^j_i$ the point the point in $E_i$ determined by the direction defined by $[v_i,v_{C^j_i}].$
By replacing $C$ by some $C^j$, we may suppose that for all $i\not\in G$, $q^j_i=q_i$ and for all $i\in G$ and $j\leq 0$, $E_i$ is the unique irreducible component of $X\setminus \mathbb{A}^2$ containing $q^j_i$, $q^j_i\not\in I(f)$ and $f|_{E_i}$ is not ramified at $q^j_i.$

We first treat the case that there exists $t>0$ such that $\sum_{i\in G}(C^j_i\cdot l_{\infty})\geq t\deg(C^j_i)$ for all $j\leq 0$.
Then we apply Proposition \ref{prototinvanonexcepdml} and \ref{procurveboundeddegfibration} to conclude our proposition in this case.

Then we may suppose that there exists a sequence of nonpositive integers $\{n_1>n_2>\cdots>n_j>n_{j+1}>\cdots\}_{j\geq 0}$ such that $$\sum_{i\in \{1,\cdots,s\}\setminus G}(C^{n_j}_i\cdot l_{\infty})\geq \deg(C^{n_j}_i)/2$$ for all $j\geq 0.$ Since $s\leq 2,$ there exists an index $i'\in \{1,\cdots,s\}\setminus G$ such that there exists infinitely many $j\geq 0$ for which $(C^{n_j}_{i'}\cdot l_{\infty})\geq 1/2\sum_{i\in \{1,\cdots,s\}\setminus G}(C^{n_j}_i\cdot l_{\infty}).$ We may suppose that $i'=1$. By picking subsequence we may suppose that for all $j\geq 0$, $(C^{n_j}_1\cdot l_{\infty})\geq 1/2\sum_{i\in \{1,\cdots,s\}\setminus G}(C^{n_j}_i\cdot l_{\infty})\geq \deg(C^{n_j})/4.$

Observe that $$d(f,v_1)A(f_{\d}(v_1))=A(v_1)+v_1(Jf),$$ then we have $(\la_1-1)A(f_{\d}(v_1))=v_1(Jf).$ If $Jf$ is a constant, then $f$ is nonramified on $\mathbb{A}^2$ and then by \cite[Theorem 1.3]{Bell2010}, our proposition holds. So we suppose that $Jf$ is not a constant.
Since $\alpha(v_1)>0$ and $Jf$ is not a constant, we have $v_1(Jf)<0$. It follows that $A(v_1)<0.$ Since $v_1$ is divisorial, $\alpha(v_1)>0$ and $A(v_1)<0$, $v_1$ is not in the boundary of $V_1.$ It implies that the direction at $v_i$ defined by $q_{1}$ is not totally invariant. By Proposition \ref{provefixnontotinvathenatr},
there exists $w\in U(\vec{v_1})$ such that
\begin{points}
\item $v_{C_1}\not\in U_{v_1,w_1}$;
\item
$f_{\d}(U_{v_1,w_1})\subseteq U_{v_1,w_1};$
\item
for all $v\in U_{v_1,w_1}$, we have $f^n_{\d}(v)\to v_1$ for $v\to \infty$;
\item there exists $N\geq 0$ such that $U(\vec{v_E})\cap \{v\in V_{\infty}|,\, \alpha(v)\geq -16\}\subseteq f_{\d}^{-N}(U_{v_1,w_1}).$
\end{points}
We may assume that $n_0\leq -N.$

The boundary $\partial f_{\d}^{-N}(U_{v_1,w_1})$ of $f^{-N}(U_{v_1,w_1})$ is finite and for every point $v\in \partial f^{-N}(U_{v_1,w_1})\setminus \{v_1\}$, we have $\alpha(v)<-16.$

Since $v_{C^{n_j}_1}\not\in f^{-N}(U_{v_1,w_1})$, there exists $w^{n_j}\in \partial f^{-N}(U_{v_1,w_1})\setminus \{v_1\}$ such that $v_{C^{n_j}_1}\geq w^{n_j}.$ Since the set $\partial f^{-N}(U_{v_1,w_1})\setminus \{v_1\}$ is finite, there exists two distinct number $l>k\geq 0$, such that $w^{n_l}=w^{n_k}.$
If $v_{C^{n_l}_1}\neq v_{C^{n_k}_1},$ we have
$$
\deg(C^{n_l})\deg(C^{n_k})=(C^{n_l}\cdot C^{n_k})$$$$
         \geq (C^{n_l}_1\cdot C^{n_k}_1)
         = (C^{n_l}\cdot l_{\infty})(C^{n_k}\cdot l_{\infty})(1-\alpha(v_{C^{n_l}_1}\wedge v_{C^{n_k}_1}))$$$$
         \geq 16^{-1}\deg(C^{n_1})\deg(C^{n_2})\times 17
         >  \deg(C^{n_l})\deg(C^{n_k}).
$$
It is impossible, so $v_{C^{n_l}_1}= v_{C^{n_k}_1},$ and then $C$ is periodic.
\endproof


\subsection{Proof of Theorem \ref{thmlaoqelat}}By Proposition \ref{proallbranchnotend}, to prove Theorem \ref{thmlaoqelat}
we may suppose that there exists $v_*\in T_f$ such that $\alpha(v_*)=0.$ By (iii) of Proposition \ref{prosecfivefa2011}, we have that $v_*$ is divisorial. It follows that $v_*$ is a rational pencil valuation. By Line Embedding Theorem, $f$ takes form $f=(F(x),G(x,y)).$ Set $d=\la_1$, we have $\deg F=d.$ Since $\la_1^2=\la_2$, $\la_1\geq 2$ and $\deg(f^n)/\la_1^n$ is bounded, by changing coordinates we may suppose that $G$ takes form $$G(x,y)=y^d+\sum_{i=0}^{d-1}a_i(x)y^i.$$

Set $m$ be an integer at least $\deg_x G+1$. Then $f$ extends to a rational morphism on $\mathbb{F}_m$ which takes form
$$f=[x_2^dF(x_1/x_2),x_2^d,x_3^d+x_2^{md}x_4^d\sum_{i=0}^{d-1}a_i(x_1/x_2)(x_3/x_2^nx_4)^i,x_4^d].$$ By calculation, we see that $f$ is an endomorphism on $\mathbb{F}_m.$
Let $L_{\infty}$ be the irreducible component of $\mathbb{F}_m\setminus \A^2$ such that $v_{L_{\infty}}=v_*$ and $F_{\infty}$ the fiber of $\pi_m$ at infinity. Set $O:=L_{\infty}\cap F_{\infty}.$

By Proposition \ref{proallbranchnotend}, we may suppose that there exists a branch $C_1$ of $C$ satisfying $v_{C_1}>v_*.$ If $C$ is a fiber of $\pi_m$, then $\pi_m(C)$ is periodic. It follows that $C$ is periodic. Otherwise, there exists a branch $C_2$ of $C$ such that the center of $C_2$ is contained in $F_{\infty}.$ It follows that $v_{C_2}\in V_{\infty}\setminus \{v\in V_{\infty}|\,\,v\geq v_*\}$. By taking $m$ large enough, we may suppose that $O\not\in C.$

Set $q_1:=C_1\cap L_{\infty}.$ If $q_1$ is not a periodic point of $f|_{L_{\infty}}$ or $q_1$ is $r$-periodic for some $r\geq 1$ and $f^r|_{L_{\infty}}$ is not ramified at $q_1$, by Proposition \ref{prototinvanonexcepdml} and then by Proposition \ref{procurveboundeddegfibration} we conclude Theorem \ref{thmlaoqelat} in this case.

Now we may suppose that $q_1$ is fixed by $f|_{L_{\infty}}$ and in some local coordinate at $q_1$, $f$ takes form $(x,y)\mapsto (x^s,y^d)$ where $2\leq s\leq d$. If $s<d$, by Corollary \ref{cormonomial}, we conclude Theorem \ref{thmlaoqelat}. If $s=d$, we conclude our Theorem by Lemma \ref{lemddtotinv}.
\endproof

\newpage

\part{Valuative dynamics in the case $\la_1^2>\la_2$}
Let $f:\mathbb{A}^2\to \A^2$ be a dominant polynomial endomorphism defined over an algebraically closed field satisfying  $\la_1^2>\la_2$. In this part, we study the dynamics of $f_{\d}$ on the valuative tree $\V_{\infty}$ at infinity.

 At first we introduce the Green function $\theta^*$ of $f$ in Section \ref{sectionbasicgreenfunction}.
 This function is a nonnegative subharmonic function on $V_{\infty}$. This function gives us many information of the dynamics of $f_{\d}$. For example, for any valuation $v\in V_{\infty}$ satisfying $\alpha(v)>-\infty$ and $\theta^*(v)>0$, we have $f^n_{\d}(v)\to v_*$ as $n\to \infty$. Next in Section \ref{sectiondyva},  we prove Theorem \ref{thmmostimportantstep} which is a strong version of Theorem \ref{thmmostimportantstepsimple}. Theorem \ref{thmmostimportantstep} is a key technique tool in the proof of our main theorem in the case $\la_1^2>\la_2$. This theorem is more useful in the case that $\#J(f)\geq 3$. In the case $\#J(f)\leq 2$ or more generally $\#J(f)<\infty$, the Green function is continuous and  piece linear. So in Section \ref{sectionjffinite} we analyzes the valuative dynamics in this case more carefully. In particular, we prove that all nondivisorial valuations in $J(f)$ are repelling periodic points. At last in Section \ref{sectionjffinitediv}, we treat the case that all valuations in $J(f)$ are divisorial. We prove that in this case either $f$ is \'etale or $f$ preserves a fibration.

\section{Basic properties of the Green function of $f$}\label{sectionbasicgreenfunction}
Let \index{$\theta^*$}$\theta^*\in \mathbb{L}^2(\mathfrak{X})$ be the  Weil divisor defined as in Appendix A of \cite{Favre2011}. In fact $\theta^*$ is contained in $\Nef_{\infty}(V_{\infty}).$
Recall that there exists an isomorphism $i:\SH(V_{\infty})\rightarrow \Nef_{\infty}(V_{\infty})$ defined in Section \ref{subsectionclassesandva}.
In the rest of this paper, we identify $\SH(V_{\infty})$ with $\Nef_{\infty}(V_{\infty})$ by $i$.
Then $\theta^*$ can be view as a function in $\mathbb{L}^2(V_{\infty})\cap \SH(V_{\infty})$. Observe that on the set $\{v\in V_{\infty}|\,\,\alpha(v)>-\infty\}$ by $\theta^*(v)=(\theta^*\cdot Z_v)$ when $\alpha(v)>-\infty$ and $\theta^*(v)=\lim_{v'<v,v'\rightarrow v}\theta^*(v')$ when $\alpha(v)=-\infty.$ Moreover we have the following
\begin{pro}\label{probasictheta}We have
\begin{points}\item $\theta^*$ is contained in $\SH^+(V_{\infty})$;
\item $\theta^*$ is decreasing;
\item $\langle\theta^*, \theta^*\rangle=0$.
\end{points}
\end{pro}
We normalize $\theta^*$ such that $\theta^*(-\deg)=1,$ and call it the {\em Green function of $f$}\index{Green function of $f$}.

\medskip

Set \index{$W(\theta^*)$}$W(\theta^*):=\{v\in V_{\infty}|\,\,\theta^*(v)=0\}$. In general, $\theta^*$ is not continuous and $W(\theta^*)$ is not closed.


But we have the following

\begin{pro}\label{corthetacontinuousgeqm}For any $M\leq 1$, $\theta^*$ is continuous in the set  $\{v\in V_{\infty}|\,\,\alpha(v)\geq M\}$.
In particular  the set $W(\theta^*)\cap \{v\in V_{\infty}|\,\,\alpha(v)\geq M\}$ is compact.
\end{pro}

%
%

To proof Proposition \ref{corthetacontinuousgeqm}, we first prove the following
\begin{pro}\label{proltwouniforminball}Let $M$ be a real number at most $1$, and $\phi$ be a function in $\mathbb{L}^2(V_{\infty})$. For any $\epsilon>0$, there exists a continuous function $\psi$ in $\mathbb{L}^2(V_{\infty})$ such that $|\phi(v)-\psi(v)|\leq \epsilon$ for all $v\in \{v\in V_{\infty}|\,\,\alpha(v)\geq M\}.$
\end{pro}
\proof[Proof of Proposition \ref{proltwouniforminball}]There exists $X\in \mathcal{C}_0$ such that $$\langle\phi-R_{\Gamma_X}\phi,\phi-R_{\Gamma_X}\phi\rangle\leq (1-M)^{-2}\epsilon^2 .$$ Set $\psi=R_{\Gamma_X}\phi$, then for all $v\in \{v\in V_{\infty}|\,\,\alpha(v)\geq M\}$  we have
$$(\phi(v)-\psi(v))^2=\langle(\phi-\psi),Z_v\rangle^2 \leq \langle(\phi-\psi),(\phi-\psi)\rangle(1-\alpha(v))\leq \epsilon^2$$ by \cite[Proposition 3.18]{Xieb}. It follows that $|\phi(v)-\psi(v)|\leq \epsilon$.
\endproof
\proof[Proof of Proposition \ref{corthetacontinuousgeqm}]By Proposition \ref{proltwouniforminball}, $\theta^*|_{\{v\in V_{\infty}|\,\,\alpha(v)\geq M\}}$ can be uniformly approximated by continuous functions on $\{v\in V_{\infty}|\,\,\alpha(v)\geq M\}$. Then itself is continuous on $\{v\in V_{\infty}|\,\,\alpha(v)\geq M\}$.
\endproof

\medskip

Define $J(f):=\Supp\Delta\theta^*$.\index{$J(f)$} Observe that $J(f)$ is a closed subset in $V_{\infty}$.
Then we have the following

\begin{pro}\label{provaluinjfnotcompare} Let $T$ be a finite closed subtree of $V_{\infty}$ containing $-\deg$ and $m_T$ the number for maximal points in $T$.
Let $\phi$ be a subharmonic function in $\SH^+(V_{\infty})$ satisfying $\langle\phi,\phi\rangle=0$. If $m_T<\# \Supp\Delta\phi$, then $\Supp\Delta\phi$ is not contained in $T$.
\end{pro}
\proof[Proof of Proposition \ref{provaluinjfnotcompare}]Otherwise we have $\Supp\Delta\phi\subseteq T$, it follows that $$R_T\phi=\phi.$$
Then for any $r> 0$, set $W_{r}:=\{v\in T|\,\,\phi(v)<r\}.$ We have $$0=\int_{T}\phi\Delta\phi\geq \int_{T\setminus W_r}\phi\Delta \phi\geq r\int_{T\setminus W_r}\Delta\phi.$$ It follows that $\Supp \Delta \phi\subseteq T\setminus(\bigcup_{r>0}W_r)=\{v\in T|\,\,\phi(v)=0\}.$ Since $\phi(v)$ is decreasing, we have  $\#\{v\in T|\,\,\phi(v)=0\}\leq m_T$ which is a contradiction.
\endproof

%
%
%
%
%
%
%
%
%
%
%
%
%
%

\medskip

For any $v_1,v_2\in V_{\infty}$, the distance is defined by
$$d(v_1,v_2):=2\alpha(v_1\wedge v_2)-\alpha(v_1)-\alpha(v_2).$$

As in Section \ref{subsectionvaldy}, denote by $v_*$ is the eigenvaluation of $f$.
\begin{pro}\label{prothetageztendstovs}Let $M$ be a real number at most one and $r$ be a positive real number. If $v\in V_{\infty}$ is a valuation satisfying $\alpha(v)\geq M$ and $\theta^*(v)\geq r$, then there are $\delta,C>0$ such that  for all $n\geq 0$, we have $d(f^n,v)>\delta$  and $$d(f_{\d}^n(v),v_*)\leq C(\frac{\la_2}{\la_1^2})^{n}.$$

In particular $f^n_{\d}(v)\rightarrow v_*$ as $n\rightarrow \infty.$

\end{pro}

\proof[Proof of Proposition \ref{prothetageztendstovs}]Let $v$ be any valuation in $V_{\infty}$ satisfying $\alpha(v)>-\infty$ and $\theta^*(v)>0$. By  \cite[Lemma A.6]{Favre2011}, $f^n_*Z_v=d(f^n,v)Z_{f_{\d}(v)}$. Then we have $d(f^n,v)\theta^*(f_{\d}(v))=(f^n_*Z_v\cdot \theta^*)=\la_1^n\theta^*(v)>0.$ It follows that $d(f^n,v)>0.$

Set $K_{r,M}:=\{v\in V_{\infty}|\,\,\alpha(v)\geq M, \theta^*(v)\geq r\}.$
By Proposition \ref{corthetacontinuousgeqm}, $K_{r,M}$ is compact. For any $n\geq 0$, set $\delta_n:=\inf_{v\in K_{r,M}}d(f^n,v).$ Since $d(f^n,v)$ is continuous and $K_{r,M}$ is compact, we have $\delta_n>0$ for all $n\geq 0.$

\medskip

Set $L:=Z_{-\deg}\in \mathbb{L}^2(\mathfrak{X})$ and $\theta_*:=Z_{v_*}.$ By Theorem \cite[Theorem A.8]{Favre2011}, we have $(\theta^*\cdot\theta_*)>0.$

As in \cite{boucksomfavrejonsson}, there exists a norm $\|\cdot\|$ on $\mathbb{L}^2(\mathfrak{X})$ defined by $\|\psi\|^2:=2(\psi\cdot L)^2-(\psi\cdot\psi)$ which makes $\mathbb{L}^2(\mathfrak{X})$ to be a Hilbert space. Observe that $\|Z_v\|=\left(2-\alpha(v)\right)^{\frac{1}{2}}\geq 1$ for all $v\in \{w\in V_{\infty}|\,\,\alpha(w)>-\infty\}$. It is easy to check that for $v_1,v_2\in \{w\in V_{\infty}|\,\,\alpha(w)>-\infty\}$, we have $$d(v_1,v_2)=\|Z_{v_1}-Z_{v_2}\|^2.$$

By \cite{boucksomfavrejonsson}, we have that for any $\psi\in \mathbb{L}^2(\mathfrak{X})$ satisfying $(\psi\cdot \theta^*)\neq 0$, we have $$\|\la_1^{-n}(\theta^*\cdot \theta_*)(\psi\cdot \theta^*)^{-1}f_*^n\psi-\theta_*\|\leq B(\psi\cdot \theta^*)^{-1}\|\psi\|(\frac{\la_2}{\la_1^2})^{\frac{n}{2}}$$ for some $B>0.$ It follows that
$$\|\la_1^{-n}(\theta^*\cdot \theta_*)(Z_v\cdot \theta^*)^{-1}f_*^nZ_v-\theta_*\|\leq Br^{-1}(2-M)^{\frac{1}{2}}(\frac{\la_2}{\la_1^2})^{\frac{n}{2}}$$
and then $$|\la_1^{-n}(\theta^*\cdot \theta_*)(Z_{v}\cdot \theta^*)^{-1}d(f^n,v)-1|=|\left((\la_1^{-n}(\theta^*\cdot \theta_*)(Z_v\cdot \theta^*)^{-1}f_*^nZ_v-\theta_*)\cdot L\right)|$$$$\leq \|\la_1^{-n}(\theta^*\cdot \theta_*)(Z_v\cdot \theta^*)^{-1}f_*^nZ_v-\theta_*\|\|L\|\leq Br^{-1}(2-M)^{\frac{1}{2}}(\frac{\la_2}{\la_1^2})^{\frac{n}{2}}.$$

Because $\frac{\la_2}{\la_1^2}<1$, there exists $N\geq 0$ such that $Br^{-1}(2-M)^{\frac{1}{2}}(\frac{\la_2}{\la_1^2})^{\frac{N}{2}}<1/2$.
For all $n\geq N$, we have  $|\la_1^{-n}(\theta^*\cdot \theta_*)(Z_{v}\cdot \theta^*)^{-1}d(f^n,v)-1|<1/2$.
It follows that $d(f^n,v)>\frac{1}{2}\la_1^n(\theta^*\cdot\theta_*)r>\frac{1}{2}(\theta^*\cdot\theta_*)r$ when $n\geq N.$

Set $\delta:=\frac{1}{2}\min\{\frac{1}{2}(\theta^*\cdot\theta_*)r, \delta_0,\cdots,\delta_{N}\}$. We have $d(f^n,v)>\delta>0$ for all $n\geq 0.$

\medskip

When $n\geq N$, we have $\la_1^{-n}(\theta^*\cdot \theta_*)(Z_{v}\cdot \theta^*)^{-1}d(f^n,v)>1/2$. It follows that
$$\frac{1}{2}\|Z_{f_{\d}^nv}\|\leq \|\theta_*\|+\|\la_1^{-n}(\theta^*\cdot \theta_*)(Z_v\cdot \theta^*)^{-1}d(f^n,v)Z_{f_{\d}^nv}-\theta_*\|\leq \|\theta_*\|+1/2.$$ It follows that $\|Z_{f_{\d}^nv}\|\leq 1+2\|\theta_*\|$ when $n\geq N$.

When $n\leq N,$ we have $$\|\la_1^{-n}(\theta^*\cdot \theta_*)(Z_v\cdot \theta^*)^{-1}d(f^n,v)Z_{f_{\d}^nv}\|
\leq \|\theta_*\|+\|\la_1^{-n}(\theta^*\cdot \theta_*)(Z_v\cdot \theta^*)^{-1}d(f^n,v)Z_{f_{\d}^nv}-\theta_*\|$$$$\leq \|\theta_*\|+Br^{-1}(2-M)^{\frac{1}{2}}.$$ It follows that
$$\|Z_{f_{\d}^nv}\|\leq \la_1^{n}(\theta^*\cdot \theta_*)^{-1}\theta^*(v)d(f^n,v)^{-1}(\|\theta_*\|+Br^{-1}(2-M)^{\frac{1}{2}})$$$$\leq \la_1^{N}(\theta^*\cdot \theta_*)^{-1}\theta^*(v)\delta^{-1}(\|\theta_*\|+Br^{-1}(2-M)^{\frac{1}{2}}).$$

Set $C_1:=\max\{2\|\theta_*\|+1, \la_1^{N}(\theta^*\cdot \theta_*)^{-1}\theta^*(v)\delta^{-1}(\|\theta_*\|+Br^{-1}(2-M)^{\frac{1}{2}})\}$, then we have $\|Z_{f_{\d}^nv}\|\leq C_1$ for all $n\geq 0.$
 Then we have
$$\|Z_{f^n_{\d}v}-\theta_*\|\leq \|\la_1^{-n}(\theta^*\cdot \theta_*)(Z_{v}\cdot \theta^*)^{-1}d(f^n,v)-1\|\|Z_{f^n_{\d}v}\|+\|\la_1^{-n}(\theta^*\cdot \theta_*)(Z_{v}\cdot \theta^*)^{-1}f_*^nZ_v-\theta_*\|$$
$$\leq (C_1+1)Br^{-1}(2-M)^{\frac{1}{2}}(\frac{\la_2}{\la_1^2})^{\frac{n}{2}}.$$ Set $C:=(C_1+1)^2B^2r^{-2}(2-M)$,
then we have $$d(f_{\d}(v),v_*)=|2\alpha(f_{\d}^n(v)\wedge v_*)-\alpha(f_{\d}^n(v))-\alpha(v_*)|=\|Z_{f^n_{\d}v}-\theta_*\|^2\leq C(\frac{\la_2}{\la_1^2})^{n}.$$
\endproof

\begin{cor}\label{corfdwelldefvinsmw}For any $v\in V_{\infty}$ satisfying $\theta^*(v)>0$, we have $d(f^n,v)>0$ and $\theta^*(f_{\d}^nv)>0$ for all $n\geq 0.$
\end{cor}
\proof[Proof of Corollary \ref{corfdwelldefvinsmw}]If $\alpha(v)>-\infty$, we conclude our corollary by Proposition \ref{prothetageztendstovs}. If $\alpha(v)=-\infty$, by \cite[Proposition 7.2]{Favre2007}, for any $n\geq 0$, there exists $w<v$ such that $d(f^n,v)=d(f^n,w).$ Since $\theta^*(w)\geq \theta^*(v)>0$ and $\alpha(w)>-\infty$, we have $d(f^n,v)=d(f^n,w)>0.$ Then we have $\theta^*(f_{\d}^nv)=d(f^n,v)^{-1}\theta^*(v)>0.$
\endproof

\bigskip

\section{Valuative dynamics of polynomial endomorphisms with $\la_1^2>\la_2$}\label{sectiondyva}

The aim of this section is to prove the following theorem.
\begin{thm}\label{thmmostimportantstep}Let $f$ be a dominant polynomial endomorphism on $\mathbb{A}^2$ defined over an algebraically closed field satisfying $\la_1^2>\la_2.$
Let $l$ be a positive integer strictly less than $\#J(f)$, $W$ be an open neighborhood of $v_*$ in $V_{\infty}$ and $k$ be a non negative integer. There exists a real number $r>0$, a finite set of polynomials $\{P_i\}_{1\leq i\leq s}$ and a positive integer $N$ such that for any finite set of valuations $\{v_1,\cdots,v_t\}$ with $t\leq l$ satisfying $\{v_1,\cdots,v_t\}\subseteq V_{\infty}\setminus (\cap_{j=0}^kf_{\d}^{-N-j}(W))$, there exists an index $i\in \{1,\cdots,s\}$ such that $v_j(P_i)> r$ for all $j\in\{1,\cdots,t\}$.
\end{thm}

Observe that Theorem \ref{thmmostimportantstepsimple} in Introduction is a direct corollary of Theorem \ref{thmmostimportantstep}.
%
%

For this purpose, we first need the following
\begin{lem}\label{lemcompactlrichtolponerich}Let $l$ be a nonnegative integer and let $W$ be a compact subset of $V_{\infty}$ such that any subset $S$ of $W$ containing at most $l$ elements is rich. Then there exists an open set $U$ containing $W$ such that for any positive integer $s$ there exists $M_s\leq 1$ such that for any subset $S_1$ of $U$ with at most $l$ elements and any subset $S_2$ of $\{v\in V_{\infty}|\,\, \alpha(v)\leq M_s\}$ with at most $s$ elements, the set $S_1\cup S_2$ is rich.
\end{lem}
\proof[Proof of Lemma \ref{lemcompactlrichtolponerich}]
Let $w=(v_1,\cdots, v_l)$ be a point in $W^l\subseteq V_{\infty}^l$. The set $\{v_1,\cdots,v_l\}$
is rich then by Proposition \ref{prodimtwothengeqz}, there exists $v_i'<v_i$ such that the set $\{v_1',\cdots,v_l'\}$ is rich. Then there exists $\phi_{w}\in \SH^+(V_{\infty})$ satisfying $\phi_w(v)=0$ for $v\in B(\{v_1',\cdots,v_l'\})$ and $\langle\phi_w,\phi_w\rangle>0.$ Set $U_w:=B(\{v_1',\cdots,v_l'\})^{\circ }$. Observe that $w\in U_w^l.$

Since $W^l$ is compact, there are finitely many points $w_1,\cdots, w_L\in W^l$ such that $W^l\subseteq \cup_{i=1}^LU^l_{w_i}$. We rename $U_{w_i}$ be $U_i$ and $\phi_{w_i}$ by $\phi_{i}.$ By Lemma \ref{lemgeqm}, there exists $M^i_s$ such that for any subset $S$ of $V_{\infty}$ satisfying
\begin{points}
\item[$\d$]$S\subseteq U_i\cup \{v|\,\,\alpha(v)\leq M^i_s\};$
\item[$\d$]$\#(S\setminus U_i)\leq s$;
\end{points}
we have that $S$ is rich.

For any point $x\in W$, set $I_x:=\{i\, |\,\, x\in U_i\}.$
Set $M_s:=\min \{M_s^i\}_{i=1,\cdots,L}$ and $U:=\cup_{x\in W}(\cap_{i\in I_x}U_i)$.

For any point $(y_1,\cdots,y_l)\in U^l$, there exists $(x_1,\cdots,x_l)\in W^l$ such that $y_i\in \cap_{j\in I_{x_i}}U_j$ for all $i=1,\cdots,l$. Since  $W^l\subseteq \cup_{i=1}^LU^l_{w_i}$, there exists $t=1,\cdots,L$ such that $(x_1,\cdots,x_l)\in U_t^l$. It follows that $t\in I_{x_i}$ for all $i=1,\cdots,l$.
Then $y_i\in \cap_{j\in I_{x_i}}U_j\subseteq U_t$ for all $i=1,\cdots,l$. Then we have $(y_1,\cdots,y_l)\in U_t^l$. It follows that $U^l\in \cup_{i=1}^LU^l_i.$ It follows that $U$ and $M_s$ are what we need.
\endproof

\begin{lem}\label{lemwflrich}Let $l$ be a positive integer strict less than $\#J(f)$. Let $S$ be a subset of $W(\theta^*)$ containing at most $l$ elements, then $S$ is rich.
\end{lem}
\proof[Proof of Lemma \ref{lemwflrich}]
Let $T$ be the subtree of $V_{\infty}$ generated by $S$ and $-\deg$. Since $\#J(f)\geq l+1$, by Proposition \ref{provaluinjfnotcompare}, $\Delta \theta^*$ is not supported by $T.$ By \cite[Proposition 3.21]{Xieb}, we have $R_T(\theta^*)\in \SH^+(V_{\infty})$ and $\langle R_T(\theta^*),R_T(\theta^*)\rangle>0$ and $R_T(\theta^*)(v)=0$ for all $v\in B(S).$ By Proposition \ref{prodimtwothengeqz}, the set $S$ is rich.
\endproof

Then we have the following
\begin{pro}\label{proanopensetrivhconwfin}For any integer $l\geq 1$ strict less than $\#J(f)$, there exists an open set $U$ and a number $M\leq 1$ such that $U$ contains $W(\theta^*)\cup \{v\in V_{\infty}|\,\,\alpha(v)\leq M\}$ and for any subset $S\subseteq U$ with $\#S\leq l$, we have that $S$ is rich.
\end{pro}
\proof[Proof of Proposition \ref{proanopensetrivhconwfin}]
By Lemma \ref{lemcompactlrichtolponerich}, there exists $M_1\leq 0$, such that for any subset $S$ of $\{v\in V_{\infty}|\,\,\alpha(v)<M_1\}$ containing at most $l$ elements, we have that $S$ rich.

By Lemma \ref{lemwflrich}, there exists $U_1$ containing $W(\theta^*)\cap\{v\in V_{\infty}|\,\,\alpha(v)\geq M_1\}$ and $M_2\leq M_1$ such that
for any subset $S$ of $U_1\cup \{v\in V_{\infty}|\,\,\alpha(v)<M_2\}$ containing at most $l$ elements, we have that $S$ rich.

By induction, we get a sequence of numbers $M_1>M_2>\cdots,>M_l$ and open set $U_1,\cdots,U_l$ satisfying $U_i$ contains $W(\theta^*)\cap\{v\in V_{\infty}|\,\,\alpha(v)\geq M_i\}$ for $i=1,\cdots,l$; for any subset $S$ of $U_i\cup \{v\in V_{\infty}|\,\,\alpha(v)<M_{i+1}\}$ containing at most $l$ elements, we have that $S$ is rich for $i=1,\cdots,l-1$ and for any subset $S$ of $U_l$ containing at most $l$ elements, we have that $S$ rich.

Set $V_0:=\{v\in V_{\infty}|\,\,\alpha(v)<M_{1}\}$; $V_i:=U_i\cup \{v\in V_{\infty}|\,\,\alpha(v)<M_{i+1}\}$ for $i=1,\cdots,l-1$ and $V_l:=U_l$.
We claim that $W(\theta^*)^l\subseteq \cup_{i=0}^lV_i^l$. Otherwise, suppose that there exists a point $w:=(v_1,\cdots,v_l)\in W(\theta^*)^l\setminus (\cup_{i=0}^lV_i^l)$. We may suppose that $\alpha(v_i)\geq \alpha(v_{i+1})$ for $i=1,\cdots,l-1$. Since $w\not\in V_l^l$, we have $\alpha(v_l)< M_l$. There exists $t$ minimal in $\{1,\cdots,l\}$ satisfying $\alpha(v_t)< M_l.$ It follows that the set $\{v_1,\cdots,v_{t-1}\}\subseteq \{v\in W(\theta^*)|\,\,\alpha(v)\geq M_l\}=\coprod_{i=1}^l\{v\in W(\theta^*)|\,\, M_{i-1}>\alpha_v\geq M_i\}$ where $M_0:=2.$ Since $t-1<l$, there exists $i\in \{1,\cdots,l\}$ such that $\{v\in W(\theta^*)|\,\, M_{i-1}>\alpha_v\geq M_i\}\cap \{v_1,\cdots,v_l\}=\emptyset.$ It follows that $\{v_1,\cdots,v_l\}\subseteq V_{i-1}$ and then  $w\subseteq V_{i-1}^l$ which contradicts our assumption.

For any $i=1,\cdots,l+1$ set $W_i:=\cap_{j\neq i-1} V_j$, then we have $\{v\in W(\theta^*)|\,\, M_{i-1}>\alpha_v\geq M_i\}\subseteq W_i$ for $i=1,\cdots,l$ and $\{v\in V_{\infty}|\,\,\alpha(v)<M_l\}\subseteq W_{l+1}$. Set $U:=\cup_{i=1}^{l+1} W_i$ and $M:=M_l-1$, we have that $(W(\theta^*)\cup \{v\in V_{\infty}|\,\,\alpha(v)\leq M\})\subseteq U$ and $U^l\subseteq \cup_{i=0}^lV_i^l.$ Then for any subset $S\subseteq U$ with $\#S\leq l$, we have that $S$ is rich.
\endproof

\proof[Proof of Theorem \ref{thmmostimportantstep}]
Pick $U$ and $M\leq -1$ as in Proposition \ref{proanopensetrivhconwfin}. For any point $v\in V_{\infty}\setminus U$, we have $\alpha(v)\geq M$ and $\theta^*(v)>0.$
Since $V_{\infty}\setminus U$ is compact and contained in $\{v\in V_{\infty}|\,\,\alpha(v)\geq M\}$, there exists $r>0$ such that $\theta^*(v)>r$ for all $v\in V_{\infty}\setminus U$. There exists $t>0$, such that the set $\{v\in V_{\infty}|\,\,2\alpha(v\wedge v_*)-\alpha(v)-\alpha(v_*)\leq t\}$ is contained in $W\cap \{v\in V_{\infty}|\,\,\alpha(v)\geq M\}.$
By Proposition \ref{prothetageztendstovs}
there exists $N$, such that we have $$f^{n}_{\d}(V_{\infty}\setminus U)\subseteq \{v\in V_{\infty}|\,\,2\alpha(v\wedge v_*)-\alpha(v)-\alpha(v_*)\leq t\}$$ for all $n\geq N.$ If follows that $V_{\infty}\setminus \cap_{j=0}^kf_{\d}^{-N-j}(W)\subseteq U.$ If follows that for any subset $S$ of $V_{\infty}\setminus \cap_{j=0}^kf_{\d}^{-N-j}(W)$ containing at most $l$ elements, the set $S$ is rich. Observe that $V_{\infty}\setminus \cap_{j=0}^kf_{\d}^{-N-j}(W)$ is compact, then we conclude our theorem by the following
\begin{lem}\label{lemrichcompact}Let $l$ be an positive integer and $Z$ be a compact subset of $V_{\infty}$ such that for any subset $S$ of $Z$ with at most $l$ elements, $S$ is rich. Then there exists a real number $r$, a finite set of polynomials $\{P_i\}_{1\leq i\leq s}$ such that for any subset $S$ of $Z$ with at most $l$ elements, there exists $i\in \{1,\cdots,s\}$ such that $v(P_i)>r$ for all $v\in S$.
\end{lem}
\endproof

\proof[Proof of Lemma \ref{lemrichcompact}]
For any point $w=(v_1,\cdots,v_l)\in Z^l$, there exists a real number $r_w>0$  and a non constant polynomial $P_w$ satisfying $v_i(P)>r_w>0$ for $i=1,\cdots,l.$ Set $U_w:=\{v\in V_{\infty}|\,\, v(P_w)>r_w\}$.
Since $Z^l$ is compact, there exist $w_1,\cdots,w_s\in Z^l$ such that $Z^l\subseteq \cup_{i=1}^sU^l_{w_i}.$ Set $U_i:=U_{w_i}$, $r=\min\{r_{w_i}\}$ and $P_i:=P_{w_i}$ for $i=1,\cdots,s$.

Let  $\{v_1,\cdots,v_t\}$ be a finite subset of $Z$ with $t\leq l$.
Set $w:=(v_1,\cdots,v_t,\cdots,v_t)\in V^l_{\infty}$, we have $w\in Z^l$. Then there exists $U_j$ for some $j=1,\cdots,s$ such that $w\in U_j^l$ and then $P_j(v_i)> r$ for $i=1,\cdots,t$.
\endproof

\section{Dynamics on $V_{\infty}$ when $J(f)$ is finite}\label{sectionjffinite}
In this section, we denote by $k$ an algebraically closed field. Let $f:\mathbb{A}^2_k\rightarrow \mathbb{A}^2_k$ be a dominant endomorphism defined over $k$ with $\la_1^2>\la_2.$ Moreover, we suppose that $\#J(f)$ is finite.

Set $J(f)=\Supp\Delta\theta^*=\{v_1,\cdots,v_s\}$ where $s$ is a positive integer.
By the definition of subharmonic functions, we may write $\theta^*=\sum_{i=1}^sr_iZ_{v_i}$ where $r_i>0$ for $i=1,\cdots,s$, $\sum_{j=1}^{s}r_i\alpha(v_i\wedge v_j)=0$ and $\sum_{i=1}^sr_i=1.$

In this situation, we have that $\theta^*$ is continuous in $V_{\infty}$ and then $$W(\theta^*)=\{v\in V_{\infty}|\,\,\theta^*(v)=0\}=B(\{v_1,\cdots,v_s\})$$ is compact. By the continuity of $f_{\d}|_{\{v\in V_{\infty}|\,\, d(f,v)>0\}}$, we have that $f_{\d}(V_{\infty}\setminus W(\theta^*))\subseteq V_{\infty}\setminus W(\theta^*)$ and for all $v\in W(\theta^*)$ satisfying $d(f,v)\neq 0$, we have $f_{\d}(v)\in W(\theta^*).$

\begin{pro}\label{proreducetofix}There exists $n\geq 1$ such that $f^n_{\d}(v_i)=v_i$ and $d(f^n,v_i)=(\la_2/\la_1)^n$ for all $i=1,\cdots,s.$
\end{pro}
\proof[Proof of Proposition \ref{proreducetofix}] For $i=1,\cdots,s,$ if $d(f,v_i)\neq 0$, we have $f_{\d}(v_i)\in W(\theta^*).$ If $f_{\d}(v_i)\in W(\theta^*)^{\circ}$, then by the continuity of $f_{\d}$ there exists $v<v_i$ such that $f_{\d}(v)\in W(\theta^*).$ Then we get a contradiction. It follows that $f_{\d}(v_i)=v_{j_i}$ for some $j_i\in \{1,\cdots, s\}.$ If $d(f,v_i)= 0$, set $j_i:=1$. Then we have $f_*Z_{v_i}=d(f,v_i)Z_{v_{i_j}}.$
Since $f_*\theta^*=\la_2/\la_1\theta^*$, we have $\la_2/\la_1(\sum_{i=1}^sr_iZ_{v_i})=\sum_{i=1}^sr_id(f,v_i)Z_{v_{j_i}}.$ Since $\Delta Z_{v_i}=\delta_{v_i}$, then we have that $Z_{v_i}$'s are linear independence. It follows that $d(f,v_i)\neq 0$ for all $i=1,\cdots,s$ and the map $i\mapsto j_i$ is a permutation of $\{1,\cdots,s\}$. By replacing $f$ by some positive iterate, we may suppose that $j_i=i$ for all $i=1,\cdots,s.$ Then $f_*Z_{v_i}=d(f,v_i)Z_{v_i}$, it follows that $d(f,v_i)=\la_2/\la_1.$
\endproof

Up to a positive iterate, we may suppose that $f_{\d}(v_i)=v_i$ and $d(f,v_i)=\la_2 /\la_1$ for all $i=1,\cdots,s.$

The following proposition shows that $f_{\d}$ is repelling at $v_i$ in the direction determinate by $[v_i,-\deg]$. Moreover it is repelling at $v_i$, if $v_i$ is irrational.

\begin{figure}
\centering
\includegraphics[width=10.5cm]{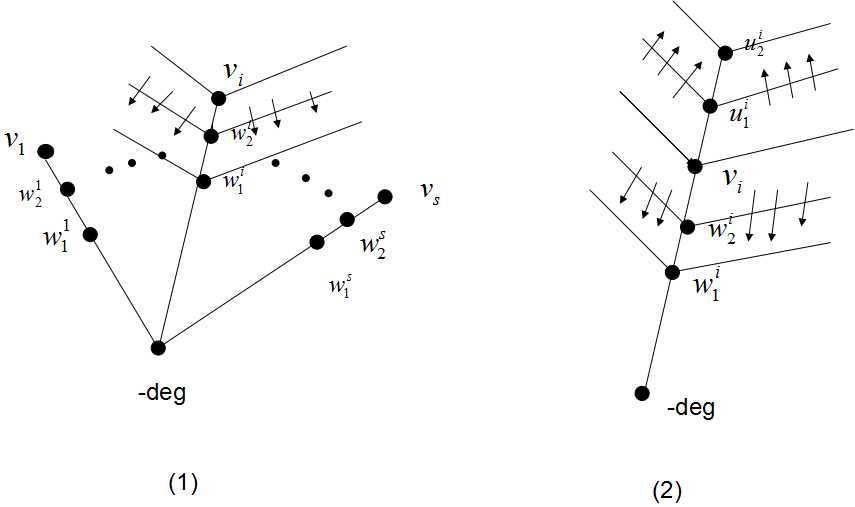}
\caption{}
\label{figc1}
\end{figure}

\begin{pro}\label{prosegmentrep}For all $i=1,\cdots,s$, there are two valuations $w^i_1<w^i_2<v_i$ as in \emph{[(1),Figure \ref{figc1}]} such that
\begin{points}
\item $f^{-1}_{\d}(\{v\in V_{\infty}|\,\,w^i_1< v\wedge v_i<v_i\})=\{v\in V_{\infty}|\,\,w^i_2< v\wedge v_i<v_i\};$
\item $f_{\d}|_{\{v\in V_{\infty}|\,\,w^i_2< v\wedge v_i<v_i\}}$ is order-preserving;
\item for all valuation $w\in [w^i_1,v_i]$, $f^{-1}_{\d}(w)$ is one point in $[w^i_1,v_i];$
\item for all valuation $w\in\{v\in V_{\infty}|\,\,w^i_1< v\wedge v_i<v_i\}$, there exists $N\geq 1$ such that $f^n_{\d}(w)\in V_{\infty}\setminus \{v\in V_{\infty}|\,\,v\wedge v_i\geq w^i_1\}$ for all $n\geq N.$
\end{points}
Moreover if $v_i$ is irrational, then there are two valuations $v_i<u^{i}_1<u^i_2$ as in \emph{[(2),Figure \ref{figc1}]}
such that
\begin{points}
\item[\text{(1).}] $f^{-1}_{\d}(\{v\in V_{\infty}|\,\,v_i< v\wedge u^i_2<u^i_2\})=\{v\in V_{\infty}|\,\,v_i< v\wedge u^i_1<u^i_1\};$
\item[\text{(2).}] $f_{\d}|_{\{v\in V_{\infty}|\,\,v_i< v\wedge u^i_1<u^i_1\}}$ is order-preserving;
\item[\text{(3).}] for all valuation $w\in [v_i,u^i_2]$, $f^{-1}_{\d}(w)$ is one point in $[v_i,u^i_1]$;
\item[\text{(4).}] for all valuation $w\in\{v\in V_{\infty}|\,\,v_i< v\wedge u^i_2<u_2^i\}$, there exists $N\geq 1$ such that for all $n\geq N$, either $d(f^n,w)=0$ or $f^n_{\d}(w)\in \{v\in V_{\infty}|\,\,u_2^i<v\}$.
\end{points}
\end{pro}
\rem The valuations $u^i_1$, $u^i_2$ and $w^i_1$, $w^i_2$ can be chosen to be arbitrarily closed to $v_i.$
\endrem
\proof[Proof of Proposition \ref{prosegmentrep}] Set $V:=V_{\infty}\setminus B(\{v_1,\cdots,v_s\})^{\circ},$ observe that $V$ is compact. By Corollary \ref{corfdwelldefvinsmw}, $f_{\d}$ is well defined on $V$ and $f_{\d}(V)\subseteq V.$ Denote by $T$ the convex hull of $\{v_1,\cdots,v_s\}\cup\{-\deg\}$.
For any $i\in \{1,\cdots,s\}$, there exists $v_i'<v_i$ such that $\{v\in V_{\infty}|\,\,v\geq v_i'\}\cap T=[v_i',v_i]$. Since $f_{\d}(v_i)=v_i$, we may further suppose that $\{v\in V_{\infty}|\,\,v\geq f_{\d}(v_i')\}\cap T=[f_{\d}(v_i'),v_i]$. For any $v\in V$ satisfying $v\geq v_i'$, we have
$$\alpha(f_{\d}(v)\wedge v_i)-\alpha(v_i)=((Z_{f_{\d}(v)}-Z_{v_i})\cdot Z_{v_i})$$$$=r_i^{-1}((Z_{f_{\d}(v)}-Z_{v_i})\cdot \theta^*)=r^{-1}_i(Z_{f_{\d}(v)}\cdot \theta^*)$$$$=r^{-1}_id(f,v)^{-1}(f_{*}Z_v\cdot \theta^*)=r^{-1}_i(\la_1/d(f,v))(Z_v\cdot \theta^*)$$$$=(\la_1/d(f,v))(\alpha(v\wedge v_i)-\alpha(v_i)).$$

Since $d(f,v_i)=\la_2/\la_1$, we have $\la_1/d(f,v_i)=\la_1^2/\la_2>1$. By assuming $v_i'$ close enough to $v_i$, we have $\la_1/d(f,v_i')>C$ for some constant $C>1$ and then $\la_1/d(f,v)>\la_1/d(f,v_i')>C.$ It follows that $$\alpha(f_{\d}(v)\wedge v_i)-\alpha(v_i)>C(\alpha(v\wedge v_i)-\alpha(v_i))\eqno (*).$$

 By \cite[Proposition 7.2]{Favre2007}, there exists a finite subtree $\mathcal{T}_f$ of $V_{\infty}$ such that $d(f,\cdot)$ is locally constant on $V_{\infty}\setminus \mathcal{T}_f$. Then $f_{\d}$ preserves the ordering on $V_{\infty}\setminus (\{v\in V_{\infty}|\,\,d(f,v)=0\}\cup \mathcal{T}_f).$
 By assuming $v_i'$ closed enough to $v_i$, we may suppose that the set $\{v\in V_{\infty}|\,\,v_i'\leq v\wedge v_i'<v_i\}\setminus [v_i',v_i]\subseteq V_{\infty}\setminus \mathcal{T}_f$. Set $t(v):=\alpha(v\wedge v_i)-\alpha(v_i)$. Since $d(f,v)$ is a decreasing piece-linear function and $d(f,v_i)=\la_2/\la_1$, there exists a constant $a\in \mathbb{Q}_{\geq 0}$, such that $$t(f_{\d}(v))=\frac{\la_1t(v)}{at(v)+\la_2/\la_1} \eqno (**)$$ for $v\in \{v\in V|,\,v_i'\leq v\wedge v_i'\}$ by assuming $v_i'$ closed enough to $v_i.$  It follows that $t(f_{\d}(v))$ strictly increases in the segment $[v_i',v_i]$ and then we have that $f_{\d}$ maps $\{v\in V|,\,v_i'\leq v\wedge v_i'\}$ onto $\{v\in V|,\,f_{\d}(v_i')\leq v\wedge f_{\d}(v_i')\}$ and it preserves the ordering.

Since $f^{-1}_{\d}(v_i)=\{v_i\}$ and $f_{\d}(V\setminus \{v\in V|,\,v_i'< v\wedge v_i'\})$ is compact, there exists $w^i_1<v_i$ such that $w^i_1>v_i'$ and $\{v\in V| w^i_1\leq v\}\cap f_{\d}(V\setminus \{v\in V|,\,v_i'< v\wedge v_i'\})=\emptyset .$ There exists $w^i_2\in (w^i_1,v_i)$ such that $\{w^i_2\}=f^{-1}_{\d}(\{w^i_1\}).$ Then the pair $(w^i_1,w^i_2)$ satisfies the conditions (i),(ii) and (iii) immediately. The inequality $(*)$ implies the condition (iv).

Now we suppose that $v_i$ is irrational. We claim the following
\begin{lem}\label{lemlocallinear}There are two valuations $w_1,w_2$ satisfying $w_1<v_i<w_2$ such that  for any $v\in \{v\in V_{\infty}|\,\,w_1<v\wedge w_2< w_2\}$ we have $d(f,v)>0$ and $$\alpha(f_{\d}(v)\wedge w_2)-\alpha(v_i)=\frac{A(\alpha(v\wedge w_2)-\alpha(v_i))}{C(\alpha(v\wedge w_2)-\alpha(v_i))+D}$$ where $A,C,D\in \mathbb{R}$.
\end{lem}
Pick valuations $w_1,w_2$ as in Lemma \ref{lemlocallinear}. Since $A,C,D$ are constants, then equation $(**)$ implies $$\alpha(f_{\d}(v)\wedge w_2)-\alpha(v_i)=\frac{\la_1(\alpha(v\wedge w_2)-\alpha(v_i))}{a(\alpha(v\wedge w_2)-\alpha(v_i))+\la_2/\la_1} \eqno (***)$$ for $v\in \{v\in V_{\infty}|\,\,w_1<v\wedge w_2< w_2\}.$
Set $V_i:=\{v\in V_{\infty}|\,\,v\geq v_i\}\setminus \{v\in V_{\infty}|\,\,d(f,v)=0\}^{\circ}$. For every valuation $v\in V_1$ satisfying $d(f,v)>0$, we have $f_{\d}(v)\in \{v\in V_{\infty}|\,\,v\geq v_i\}$. By \cite[Theorem 7.1]{Favre2007}, $f_{\d}$ extends to a continuous map $f_{\d}:V_i\rightarrow \{v\in V_{\infty}|\,\,v\geq v_i\}$. Since $\{v\in V_i|\,\,v\geq w_2\}$ is compact and $f^{-1}_{\d}(\{v_i\})=v_i$, there exists $u^i_2\in (v_i,w_2)$, such that $\{v\in V_{\infty}|\,\,v_i\leq v\leq u^i_2\}\cap f_{\d}(\{v\in V_i|\,\,v\geq w_2\})=\emptyset$. There exists $u^i_1\in (v_i,u^i_2)$ such that $\{u^i_1\}=f^{-1}_{\d}(\{u^i_2\}).$ Then equation $(***)$ implies that the pair $(u^i_1,u^i_2)$ satisfies the conditions (i),(ii), (iii) and (iv).
\endproof
\proof[Proof of Lemma \ref{lemlocallinear}]There exists a nonconstant polynomial $P$ such that there exists a branch $C_1$ of $\{P=0\}$ satisfying $v_{C_1}>v_i$ and there exists a branch $D_1$ of $\{P=0\}$ satisfying $v_{D_1}\wedge v_i<v_i$. Set $C_1,\cdots,C_s$ be all branch of $\{P=0\}$ satisfying $v_{C_j}>v_i$ and set $D_1,\cdots,D_t$ be all branch of $\{P=0\}$ satisfying $v_{D_j}\wedge v_i<v_i$. Since $v_i$ is irrational, $C_1,\cdots,C_s,D_1,\cdots,D_t$ are all branches of $\{P=0\}.$ It follows that there exists $m_1,\cdots,m_s,n_1,\cdots,n_t\in \mathbb{Z}^+$ such that for all $v\in V_{\infty}$, $v(P)=\sum_{j=1}^sm_j\alpha(v_{C_j}\wedge v)+\sum_{j=1}^tn_j\alpha(v_{D_j}\wedge v).$ Similarly, write $v(f^*P)=\sum_{j=1}^rm_j'\alpha(v_{C_j'}\wedge v)+\sum_{j=1}^ln_j'\alpha(v_{D_j'}\wedge v)$ where $n_j',m_j'\in \mathbb{Z}^+$, $v_{C_j'}>v_i$ and $v_{D_j'}\wedge v_i<v_i$.

Set $M:=\sum_{j=1}^sm_j$ and $M':=\sum_{j=1}^rm_j'$. We have $M,N,M',N'\in \mathbb{Z}_{\geq 0}$ and $M,N>0.$
Set $w_1':=\max (\{v_{C_j}\wedge v_i\}\cup \{v_{C_j'}\wedge v_i\})$ and $w_2':=(\wedge_{j=1}^t v_{D_j})\wedge (\wedge_{j=1}^{l} v_{D_j'}).$ Since $v_1$ is irrational, we have $w_1'<v_i<w_2'.$ For any $v\in \{v\in V_{\infty}|\,\,w_1'<v\wedge w_2'< w_2'\}$, we have $$v(P)=M\alpha(v\wedge w_2')+T$$ where $T=\sum_{j=1}^tn_j\alpha(v_{D_j}\wedge w_1')$ and  $$v(f^*P)=M'\alpha(v\wedge w_2')+L$$ where $T=\sum_{j=1}^ln_j'\alpha(v_{D_j'}\wedge w_1')$. On the other hand, we have $d(f,v)f_{\d}(v)(P)=v(f^*P)$. Since $v_i$ is irrational, $d(f,v_i)=\la_2/\la_1$ and $d(f,\cdot)$ is decreasing, there exist $w_1\in [w_1',v_i)$ and $w_2\in (v_i,w_2']$ such that for all $v\in \{v\in V_{\infty}|\,\,w_1<v\wedge w_2< w_2\}$,
\begin{points}
\item[$\d$] we have $d(f,v)=\la_2/\la_1+K(\alpha(v\wedge w_2)-\alpha(v_i))$ for some constant $K\in \mathbb{Q}_{\geq 0}$;
\item[$\d$]$d(f,v)>0$;
\item[$\d$]$f_{\d}(v)\in \{u\in V_{\infty}|\,\,w_1'<u\wedge w_2'< w_2'\}$.
\end{points}
It follows that for all $v\in \{u\in V_{\infty}|\,\,w_1<u\wedge w_2< w_2\}$, we have $$\left(\la_2/\la_1+K(\alpha(v\wedge w_2)-\alpha(v_i))\right)(M\alpha(f_{\d}(v)\wedge w_2)+T)=M'\alpha(v\wedge w_2)+L  \eqno (1).$$
Set $v=v_i$ in equation $(1)$, then we have $$\la_2/\la_1(M\alpha(v_i)+T)=M'\alpha(v_i)+L  \eqno (2).$$
Set $t(v):=\alpha(v\wedge w_2)-\alpha(v_i)$. By taking difference $(1)-(2)$, we have $$(\la_2/\la_1)Mt(f_{\d}(v))+K(M\alpha(f_{\d}(v)\wedge w_2)+T)t(v)=M't(v).$$ It follows that $$\left((\la_2/\la_1)M+KMt(v)\right)t(f_{\d}(v))+K\left(M\alpha(v_i)+T\right)t(v)=M't(v)$$ and then we have
$$t(f_{\d}(v))=\frac{\left(M'-K\left(M\alpha(v_i)+T\right)\right)t(v)}{(\la_2/\la_1)M+KMt(v)}$$ for $v\in \{u\in V_{\infty}|\,\,w_1<u\wedge w_2< w_2\}$ by taking $w_1,w_2$ closed enough to $v_i.$
\endproof

\section{When $J(f)$ is a finite set of divisorial valuations}\label{sectionjffinitediv}
In this section $k$ is an algebraically closed field.
Let $f:\mathbb{A}^2_k\rightarrow \mathbb{A}^2_k$ be a dominant endomorphism defined over $k$ with $\la_1^2>\la_2$ such that $J(f)=\Supp\Delta \theta^*$  is a finite set of divisorial valuations.

We first fix the setting.
Write $\theta^*=\sum_{i=1}^sr_iZ_{v_i}$ where $r_i>0$ and $v_i$ is divisorial for $i=1,\cdots,s$. The coefficients $r_i$'s satisfy the following conditions:
\begin{points}
\item $\sum_{j=1}^{s}r_i\alpha(v_i\wedge v_j)=0$ for $i=1,\cdots,s$;
\item $\sum_{i=1}^sr_i=1.$
\end{points}

By Proposition \ref{proreducetofix}, we may suppose that $f_{\d}(v_i)=v_i$ and $d(f,v_i)=\la_2 /\la_1$ for all $i=1,\cdots,s.$
The aim of this section is the following

\begin{thm}\label{thmsuppfinitdivcarzfi}If $Jf$ is not a constant, then $f$ preserves a nontrivial fibration  $G\in k[x,y]\setminus k$ i.e. there there exists a polynomial morphism $G:\mathbb{A}^1_k\rightarrow \mathbb{A}^1_k$ such that $P\circ f=G\circ P.$ Moreover we have $R_{\{v_1,\cdots,v_s\}}=k[P]$.
\end{thm}

\begin{rem}In the preparing work \cite{Jonssona}, Jonsson, Wulcan and I show that $Jf$ can not be constant in this case.
\end{rem}

\proof[Proof of Theorem \ref{thmsuppfinitdivcarzfi}]For every polynomial $Q\in k[x,y]$, set $\theta^*(Q):=\sum_{i=1}^sr_iv_i(Q).$
Recall that the function $\log|Q|:v\mapsto -v(Q)$ on $V_{\infty}$ can be written as $$\log|Q|(v)=\sum_{i=1}^{l}m_i\alpha(Q_i\wedge v)$$ where $Q_i$'s are all curve valuations associated to the branches at infinity of $\{Q(x,y)=0\}$ and $m_i\geq 1$. Then we have
$$\theta^*(Q)=\sum_{i=1}^sr_iv_i(Q)=-\sum_{i=1}^sr_i\log|Q|(v_i)$$$$=-\sum_{i=1}^sr_i(\sum_{j=1}^{l}m_j\alpha(Q_j\wedge v_i))
=-\sum_{i=1}^sr_i\sum_{j=1}^{l}m_j(Z_{Q_j}\cdot Z_{v_i})$$$$=-\sum_{j=1}^{l}(Z_{Q_j}\cdot \sum_{i=1}^sr_iZ_{v_i})=-\sum_{j=1}^{l}(Z_{Q_j}\cdot \theta^*)\leq 0.$$
For all element $Q\in R_{\{v_1,\cdots,v_s\}}$, we have $0\geq \theta^*(Q)=\sum_{i=1}^sr_iv_i(Q)\geq 0$. It follows that $v_i(Q)=0$ for all $i=1,\cdots,s.$
By Proposition \ref{prodimtwothengeqz}, the transcendence degree of $\Frac(R_{\{v_1,\cdots,v_s\}})$ is at most one.

On the other hand, we claim the following
\begin{pro}\label{prosuppdfinitedivfibration}If $v_i$ is divisorial for all $i=1,\cdots,s$, then either $Jf$ is a constant or there exists a polynomial $P\in k[x,y]\setminus k$ such that $v_i(P)\geq 0$ for all $i=1,\cdots,s$.
\end{pro}
By Proposition \ref{prosuppdfinitedivfibration}, the transcendence degree of $\Frac(R_{\{v_1,\cdots,v_s\}})$ is at least one and then equals to one. By \cite[Proposition 5.8]{Xieb}, there exists a nonconstant polynomial $P\in k[x,y]$ such that
$R_{\{v_1,\cdots,v_s\}}= k[P]$.

Observe that $v_i(P\circ f)=(f_*v_i)(P)=\la_2/\la_1v_i(P)=0$ for all $i=1,\cdots,s.$ Then we have $P\circ f\in R_{\{v_1,\cdots,v_s\}}= k[P]$. Then there exists $G\in k[t]$ such that $P\circ f=G\circ P$.
\endproof

\proof[Proof of Proposition \ref{prosuppdfinitedivfibration}]
Since $v_1$ is divisorial, we have $\la_2/\la_1=d(f,v_1)\in \mathbb{Z}^+$, it follows that $\la_2\geq \la_1.$

We define $A(\theta^*):=\sum_{i=1}^sr_iA(v_i).$ As in the beginning of the proof of Theorem \ref{thmsuppfinitdivcarzfi}, we set $\theta^*(Q):=\sum_{i=1}^sr_iv_i(Q)$ for all polynomial $Q\in k[x,y]$ and we have $\theta^*(Q)\leq 0$ for all $Q\in k[x,y].$ Observe that $$\la_2/\la_1A(\theta^*)=\sum_{i=1}^sr_i(\la_2/\la_1A(v_i))$$$$=\sum_{i=1}^sr_i(A(v_i)+v_i(Jf))=A(\theta^*)+\theta^*(Jf).$$
It follows that $(\la_2/\la_1-1)A(\theta^*)=\theta^*(Jf)\leq 0.$

We claim the following
\begin{lem}\label{lemlaoelatsuppfinite}If $\la_1=\la_2$, then either $Jf$ is a constant or there exists a polynomial $P\in k[x,y]\setminus k$ such that $v_i(P)\geq 0$ for all $i=1,\cdots,s$.
\end{lem}

By Lemma \ref{lemlaoelatsuppfinite}, we may suppose that $\la_2>\la_1$, it follows that $A(\theta^*)\leq 0.$
Then we conclude our Proposition by \cite[Proposition 5.6]{Xieb}.
\endproof
\proof[Proof of Lemma \ref{lemlaoelatsuppfinite}]We may suppose that $Jf$ is not a constant. For all $i=1,\cdots,s$, we have formula $$d(f,v_i)A(v_i)=A(v_i)+v_i(Jf).$$ Since $d(f,v_i)=\la_2/\la_1=1$, we have $v_i(Jf)=0$ for all $i=1,\cdots,s.$ Set $P=Jf$, then we conclude of lemma.
\endproof

\newpage

\part{The non-resonant case $\la_1^2>\la_2$}\label{sectiongreenf}
In this part,  $f$ is a dominant polynomial endomorphism defined over an algebraically closed field satisfying $\la_1^2>\la_2.$
The aim of this part is to prove the main theorem in the case $\la_1^2>\la_2$ which completes the proof of Theorem \ref{thmdmlpoly}.

\begin{thm}\label{thmlaonesqsllatwo}Let $f:\mathbb{A}^2_{\overline{\mathbb{Q}}}\rightarrow \mathbb{A}^2_{\overline{\mathbb{Q}}}$ be a dominant polynomial endomorphism on $\mathbb{A}^2_{\overline{\mathbb{Q}}}$ satisfying $\la_1^2>\la_2.$ Then the pair $(\mathbb{A}^2_{\overline{\Q}},f)$ satisfies the DML property.
\end{thm}

\section{The case $\la_1^2>\la_2$ and $\#J(f)\geq 3$}
Our aim of this part is to prove the following
\begin{thm}\label{thmdmlsuppgeqthree}Set $k=\overline{\mathbb{Q}}$. Let $f$ be a dominant polynomial endomorphism on $\mathbb{A}^2_k$ defined over $k$ with $\la_1(f)^2>\la_2(f)$ and $\#J(f)\geq 3$. Then the pair $(\mathbb{A}^2_k,f)$ satisfies the DML property.
\end{thm}
\medskip
We first fix the notations.

Let $C$ be an irreducible curve in $\mathbb{A}^2_k$ and $p$ be a closed point in $\mathbb{A}^2_k$.
We suppose that  $\{n\in \mathbb{N}|\,\,f^n(p)\in C\}$ is infinite and
$p$ is not preperiodic.
By Theorem \ref{thmasequencecurves}, we may suppose that there exists a sequence of curves $\{C_i\}_{i\in \mathbb{Z}}$ with $s\in\{1,2\}$ places at infinity such that
\begin{points}
\item[$\d$] $C^0=C$;
\item[$\d$] $f(C^i)=C^{i+1}$;
\item[$\d$] for all $i\in \mathbb{Z}$, the set $\{n\geq 0|f^n(p)\in C^i\}$ is infinite.
\end{points}
Let $C^j_i$'s be branches of $C^j$, we may suppose that $f(C^j_i)=C_i^{j+1}$ for $j\leq -1$ and $1\leq i\leq s.$

The following lemma is a key ingredient of our proof which is a direct application of Section \ref{sectiondyva}
\begin{lem}\label{lemexwtowcninw}If there exists an open set $W$ of $V_{\infty}$ containing $v_*$ and a nonnegative integer $L\geq 0$, such that for infinitely many $j\leq 0$ we have
$v_{C^j_i}\not\in \cap_{k=0}^Lf_{\d}^{-k}(W)$ for all $i=1,\cdots,s,$
then the pair $(\mathbb{A}^2,f)$ satisfies the DML property for $C$.
\end{lem}

\proof[Proof of Lemma \ref{lemexwtowcninw}]
Since $\#J(f)\geq 3>s$, by Theorem \ref{thmmostimportantstep},
there exists a finite set of polynomials $\{P_i\}_{1\leq i\leq l}$ and a positive integer $N$ such that for any set of valuations $\{v_1,\cdots,v_s\}$ of $s$ elements satisfying $v_i\not\in \cap_{k=0}^Lf^{-N-k}_{\d}(W)$ for all $i=1,\cdots,s$, there exists an index $i\in \{1,\cdots,l\}$ such that $v_j(P_i)> 0$ for all $j\in\{1,\cdots,s\}$. Let $S$ be the infinite set of index $j\leq 0$ such that $v_{C^j_i}\not\in \cap_{k=0}^Lf^{-k}_{\d}(W)$ for all $i=1,\cdots,s.$ Denote by $S^{-N}$ the set of index $j$ such that $j+N\in S$.

Since $v_{C^j_i}\not\in W$ for all $j\in S$, we have $v_{C^{j}_i}\not\in f^{-N}_{\d}(W)$ for all $j\in S^{-N}.$  Denote by $R$ the finite set of irreducible polynomials which is a factor of one polynomial $P_i$, $i\in\{1,\cdots,l\}$.

 For any $j\in S^{-N}$, there exists an index $k\in \{1,\cdots,l\}$ such that $v_{C^j_i}(P_k)> 0$ for all $i\in\{1,\cdots,s\}$. Then $P_k$ has no poles but zeros in the Zariski closure of $C^j$ in $\P^2$.
 It follows that we have $P_k|_{C^j}=0$ and then $C^j$ is defined by $Q_j=0$ where $Q_j$ is an irreducible polynomial in $R.$ Since $R$ is finite, there exists $j_1<j_2\in S^{-N}$ such that $C^{j_1}=C^{j_2}.$ It follows that $C$ is periodic.
\endproof

%
%
%

In the rest of this section we present our proof in the situation $s=2$ and we will give a remark for the situation $s=1$ in every case.
\subsubsection*{\textbf{1) The case $v^*$ is not divisorial}}

\smallskip



By \cite[Theorem 3.1]{Favre2011}, there exists an open set $W$ of $V_{\infty}$ containing $v_*$ such that
\begin{points}
\item[$\d$]
$v_{C^0_i}\not\in W$ for $i=1,2;$
\item[$\d$]
$f_{\d}(W)\subseteq W$.
\end{points}
Then we have $W\subseteq f_{\d}^j(W)$ for all $j\leq 0$. It follows that $v_{C^j_i}\not\in W$ for all $j\leq 0$ and $i=1,2$. By applying Lemma \ref{lemexwtowcninw}, we conclude our proposition in this situation.

\rem
When $s=1$, the proof is the same.
\endrem

\medskip

%
\subsubsection*{\textbf{2) The case $v_*$ is divisorial}}

 There exists a smooth projective compactificaition $X$ of $\mathbb{A}^2$ containing a divisor $E$ satisfying $v_E=v_*$.
 By \cite[Lemma 4.6]{Favre2011} we may suppose that for any point $t$ in $I(f)\cap E$, $t$ is not a periodic point of $f|_E$.

\subsubsection*{\textbf{2.1) The case $\deg(f|_E)=\id$.}}
The proof of this case is similar to Case 1).

There exists a compactification $X\in \mathcal{C}$ such that $E$ is an irreducible component of $X\setminus \A^2$ and $I(f)\cap E=\emptyset$. If follows that there exists an open set $W$ of $V_{\infty}$ containing $v_*$ such that
\begin{points}
\item[$\d$]
$v_{C^0_i}\not\in W$ for $i=1,2;$
\item[$\d$]
$f_{\d}(W)\subseteq W$.
\end{points}
Then we have $W\subseteq f_{\d}^j(W)$ for all $j\leq 0$. It follows that $v_{C^j_i}\not\in W$ for all $j\leq 0.$ Apply Lemma \ref{lemexwtowcninw} and we conclude our proposition in this situation.

\rem
When $s=1$, the proof is the same.
\endrem

\subsubsection*{\textbf{2.2) The case $\deg(f|_E)=1$ and $f|_E^n\neq \id$ for all $n\geq 0$}}

Since $\deg(f|_E)=1$, $f|_E$ has at most two periodic points. By replacing $f$ by a positive iterate, we may suppose that all periodic points of $f|_E$ are fixed.

\medskip

In the case 1) and the case 2.1), there exists a system of invariant neighborhood of $v_*$. Unfortunately, such a system does not exist in this case.
But there exists a system of neighborhood $W$ of $v_*$ which is not invariant but play a similar role as invariant neighborhood of $v_*$ play in the case 1) and 2).

\begin{defi}\label{definice}A neighborhood $W$ of $v_*$ is said to be a {\em nice neighborhood of $v_*$} \index{nice neighborhood of $v_*$} if it satisfies the following properties:
\begin{points}
\item for all valuation $v\in W$, $d(f,v)>0$ and the center of $v$ is contained in $E$;
\item for any point $t\in E$, we have $f_{\d}(U(t)\cap W)\subseteq U(f|_{E}(t))$;
\item for all $j\leq 0$ such that there exists a branch $C^j_i$ of $C^j$ at infinity satisfying $v_{C^j_i}\in W$,
we have $\deg f|_{C^j}\leq \la_1$ for all $j\leq -1$;
\item its boundary $\partial W$ is finite;
\item for any fixed point $x\in E$, $f_{\d}(U(x)\cap W)\subseteq U(x)\cap W$.
\end{points}
\end{defi}

\begin{lem}\label{lemdegeeonecjcjpolo} If $v_*=v_E$ is divisorial, $\deg f|_E=1$ such that all periodic points of $f|_E$ are fixed.
Let $U$ be any neighborhood of $v_*$, there exists nice neighborhood $W$ of $v_*$ contained in $U$.
\end{lem}

\proof[Proof of Lemma \ref{lemdegeeonecjcjpolo}]
There exists a compactification $Y\in \mathcal{C}$ dominates $X$ with morphisms $\pi_1:Y\rightarrow X$, $\pi_2:Y\to X$ where $\pi_1$ is birational and $\pi_2\circ\pi_1^{-1}=f.$ Set $E'$ the strict transform of $E.$ We may suppose that for every irreducible component $F\neq E'$ of $Y\setminus \A^2$  satisfying $F\cap E'\neq\emptyset$, we have that $\pi_1(F)$ is a point, $\pi_2(F)=f|_E(\pi_1(F))$ and $\pi_2$ at every point in $E'$ is locally monomial ( see \cite[Theorem 3.2]{Cutkosky2002}). Denote by $W_Y$ the set of all valuations whose centers on $Y$ are contained in $E'.$

Then we pick a neighborhood $W'$ satisfying conditions (i) and (ii) in Definition \ref{definice} and contained in $W_Y\cap U$.

\smallskip
We will first show that $W'$ satisfies condition (iii) in Definition \ref{definice}.

Fix $j\leq -1,$ we may suppose that $v_{C^j_1}\in W'$. By condition (ii), the center of $f(v_{C^{j+1}_1})$ is contained in $E.$ Write $c^j$ for the center of  $v_{C^j_1}$ on $Y$ and $c^{j+1}$ for the center of $f(v_{C^{j+1}_1})$ on $X$. By condition (iii), $c^j$ is contained in $E_Y$. There exists a local coordinate $U^{j}$ at $c^j$ satisfying $E_Y=\{y=0\}$ in $U^j$ and  a local coordinate $U^{j+1}$ of $c^{j+1}$ satisfying $E=\{y=0\}$ in $U^{j+1}.$ Since $\pi_2|_{E_Y}$ is linear and $d(f,v_E)=\la_1$. We may ask that the map $\pi_2:U^j\rightarrow U^{j+1}$ to take form $(x,y)\mapsto (x,x^my^{\la_1})$ for some $m\geq 0.$ It follows that for a general point in $U^{j+1}\setminus \{(0,0)\}$, it has at most $\la_1$ preimages by $\pi_2$ in $U^j.$ Pick a general point in $C^j$ near $c^{j+1}$, it has at most $\la_1$ preimages by $f|_{C^{j}}$ near the center of $v_{C^j_1}$. It follows that $\deg f|_{C^j}\leq \la_1$ which shows that $W'$ satisfies condition (iii) in Definition \ref{definice}.

\smallskip

Observe that all neighborhoods of $v_*$ contained in $W'$ satisfies conditions (i), (ii) and (iii).

By  replacing $W'$ by a neighborhoods of $v_*$ contained in $W'$, we may suppose that it also satisfies condition (iv).

If $f|_E\neq \id$,
denote by $F$ the set of fixed points of $f|_E$. Then we have  $\#F\leq 2$. By Lemma \ref{lemlocalallessinftytendtocu}, for any $x\in F$, there exists a valuation $w_x\in U(x)$ such that $\{v\in V_{\infty}|\,\,v_E<v\wedge v_E<w_x\}\subseteq W'$ and $f_{\d}(\{v\in V_{\infty}|\,\,v_E<v\wedge v_E<w_x\})\subseteq \{v\in V_{\infty}|\,\,v_E<v\wedge v_E<w_x\}$.
Set $W:=W'\setminus (\cup_{x\in F}(U(x)\setminus\{v\in V_{\infty}|\,\,v_E<v\wedge v_E<w_x\}))$, then $W$ is a nice neighborhood $W$ of $v_*$ contained in $U$.

If $f_E=\id$, the argument in 2.1) shows that $v_*$ is attracting i.e. there exists a neighborhoods $W$ of $v_*$ satisfying $f_{\d}(V)\subseteq V$. Moreover we may suppose that the boundary of $W$ is finite and $W\subseteq W'$.  Then $W$ is a nice neighborhood $W$ of $v_*$ contained in $U$.
\endproof
\smallskip

\textbf{
In the rest of the proof of the case 2.2), we take $W$ to be a nice neighborhood of $v_*$.}

\smallskip

By Lemma \ref{lemexwtowcninw} and by replacing $C$ by some $C^{j_0}$, $j_0\leq 0$, we may suppose that for all $j\leq 0$, there exists a branch $C^j_i$ of $C^j$ at infinity  such that $v_{C^j_i}\in W$.
Now we may suppose that $\deg f|_{C^j}\leq \la_1$ for all $j\leq -1$.

\smallskip

For any branch $C^{j}_i$ of $C^j$ at infinity $j\leq 0$, denote by $m^j_i$ the intersection number $(C^j_i\cdot l_{\infty})$.
Then we want to study the growth of the intersection number $m^j_i$ when $v_{C^j_i}$ is contained in $W$.

Since $d(f,v)$ is locally constant outside a finite tree, there are finitely many directions $\vec{w}_1,\cdots,\vec{w_d}$ at $v_E$ such that $d(f,v)=d(f,v_E)=\la_1$ on $V_{\infty}\setminus \cup_{i=1}^dU(\vec{w}_i).$ Denote by $t_i$ the point in $E$ determined by $\vec{w}_i.$

Since $d(f,v)$ is continuous on $V_{\infty}$, by replacing $W$ by some small open set, we may suppose that for all $v\in W$, $d(f,v)\in (2^{-1/d}\la_1,2^{1/d}\la_1).$

By Lemma \ref{lembrafact}, we have $$m^j_id(f,v_{C^j_i})=\deg(f|_{C^j_i})m^{j+1}_i=\deg(f|_{C^j})m^{j+1}_i$$ for all $j\leq 0$, $i=1,\cdots,s$.
\begin{lem}\label{lemconseqinwbound}
 If there are $i=1,2$, $j\leq 0$ and $k\geq 0$, such that $v_{C^j_i},\cdots,v_{C^{j-k}_i}\in W$ and the centers $q^j_i,\cdots,q^{j-k}_i$ are distinct, then we have $m^{j-k}_i/m^j_i\leq 2$.
\end{lem}
\rem This lemma holds also when $s=1$.
\endrem
\proof[Proof of Lemma \ref{lemconseqinwbound}]
Since $m^j_id(f,v_{C^j_i})=\deg(f|_{C^j})m^{j+1}_i$, we have $$m^j_i/m^{j+1}_i=\deg(f|_{C^j})/d(f,v_{C^j_i})\leq \la_1/d(f,v_{C^j}).$$
When $q^j_i\not\in\{t_1,\cdots,t_d\}$, we have $d(f,v_{C^j})=\la_1$, and then $m^j_i/m^{j+1}_i\leq 1;$
When $q^j_i\in\{t_1,\cdots,t_d\}$, we have $d(f,v_{C^j})\geq2^{-1/d}\la_1$, and then $m^j_i/m^{j+1}_i\leq 2^{1/d}.$
Since $q^j_i,\cdots,q^{j-k}_i$ are distinct, we have $m^{j-k}_i/m^j_i\leq (2^{1/d})^d=2$.
\endproof

\medskip

%
%
%
%
%

%

Observe that
some $v_{C^j_i}$ can be outside $W$ infinitely many times. But however, the following lemma tell us that $m^j_i/(\deg C^j)$ can not be too big.

\begin{lem}\label{lemifvjinotinwthenbound}For any nice neighborhood $W$, there exists $A\geq 0$, such that if there are infinitely many $v_{C^j_1}\not\in W$ satisfying $m^j_1/m^j_2\geq A$, then the pair $(\mathbb{A}^2,f)$ satisfies the DML property for $C$.
\end{lem}

The map $f_{\d}:\overline{\{v\in V_{\infty}|\,\,d(f,v)>0\}}\rightarrow V_{\infty}$ is continuous and the image of any $v\in \partial{\{v\in V_{\infty}|\,\,d(f,v)>0\}}$ is a curve valuation defined by a rational curve with one place at infinity. So there exists $\delta>0$ such that
$v_E\not\in f_{\d}(\{v\in V_{\infty}|\,\,d(f,v)\leq \delta\}).$ So we may take $W$ to be a nice neighborhood of $v_*$ contained in the open set
 $ V_{\infty}\setminus f_{\d}(\{v\in V_{\infty}|\,\,d(f,v)\leq \delta\})$.

By replacing $C$ by some  $C^j$, $j\leq 0$, we may suppose that for all $j\leq 0$, we have $m^j_1/m^j_2<A$  when $v(C^j_1)\not\in W$ and $m^j_2/m^j_1<A$  when $v(C^j_2)\not\in W$.

\medskip

\proof[Proof of Lemma \ref{lemifvjinotinwthenbound}]
By Theorem \ref{thmmostimportantstep}, there exists $r>0$, $N\geq 0$ and a finite set of polynomials $\{P_1,\cdots,P_m\}$ such that for any for any $v\in V_{\infty}\setminus f^{-N}_{\d}(W)$, there exists $i=1,\cdots,m$ such that $v(P_i)>r$.

Set $A:=r^{-1}(\deg(f))^{2N}\max\{\deg(P_1),\cdots,\deg(P_m)\}.$ We claim
\begin{lem}\label{lemcjminkinotfarcji} For any $j\leq 0$ and $k\geq 0$, we have
$$(\deg(f))^{2k}(m^{j}_1/m^j_2+1)\geq m^{j-k}_1/m^{j-k}_2\geq (\deg(f))^{-2k}\frac{m^{j}_1/m^j_2}{1+m^{j}_1/m^j_2}.$$
\end{lem}
If $v_{C^j_1}\not\in W$ and $m^j_1/m^j_2\geq A$,
by Lemma \ref{lemcjminkinotfarcji}, we have $v_{C^{j-N}_1}\not\in f^{-N}_{\d}(W)$ and
$$m^{j-N}_1/m^{j-N}_2\geq (\deg(f))^{-2N}A/(1+A)> (\deg(f))^{-2N}A.$$ There exists $i=1,\cdots,m$
such that $v_{C^{j-N}_1}(P_i)>r$.
Observe that
$$m_1^{j-N}v_{C^{j-N}_1}(P_i)+m_2^{j-N}v_{C^{j-N}_2}(P_i)> m_1^{j-N}r-m_2^{j-N}\deg(P_i)$$$$\geq m_1^{j-N}r-m_2^{j-N}\max\{\deg(P_1),\cdots,\deg(P_m)\}>0.$$

We claim the following
\begin{lem}\label{lemzeroatinftynotmany}Let $C$ be an irreducible curve in $\mathbb{A}^2_{\overline{\mathbb{Q}}}$ and let $C_1,\cdots,C_s$ be all the branches of $C$ at infinity. Let $P$ be any polynomial in $\overline{\Q}[x,y] .$ If $\sum_{i=1}^s(C_i\cdot l_{\infty})v_{C_i}(P)>0$, then $P|_C=0$.
\end{lem}

Lemma \ref{lemzeroatinftynotmany} implies that $P_i|_{C^{j-N}}$ is zero. It follows that $C^{j-N}$ is an irreducible component of $\{\prod_{i=1}^mP_i=0\}.$

If there are infinitely many such $j\leq 0$, there exists $j_1<j_2<0$, such that $C^{j_1-N}=C^{j_2-N}$. It follows that $C$ is periodic, which conclude our Theorem \ref{thmdmlsuppgeqthree}.
\endproof
\proof[Proof of Lemma \ref{lemcjminkinotfarcji}]Observe that $$m^{j-k}_1/m^j_1=\deg(f^k|_{C^{j-k}})/d(f^k,v_{C^{j-k}})\geq 1/(\deg(f))^k.$$
On the other hand, we have $$m_2^{j-k}\leq\deg (C^{j-k})\leq \deg (f^{*k}(C^{j}))=(\deg(f))^k\deg(C^j)=(\deg(f))^k(m^{j}_1+m^{j}_2).$$
It follows that
$$m^{j-k}_1/m^{j-k}_2\geq (\deg(f))^{-2k}m^j_1/(m^{j}_1+m^{j}_2)=(\deg(f))^{-2k}\frac{m^{j}_1/m^j_2}{1+m^{j}_1/m^j_2}.$$

The same we have
$$(\deg(f))^{2k}(m^{j}_1/m^j_2+1)\geq m^{j-k}_1/m^{j-k}_2.$$
 \endproof

\proof[Proof of Lemma \ref{lemzeroatinftynotmany}]
 We extend $C$ to a projective curve in $\mathbb{P}^2.$ By contradiction, we suppose that $P|_C$ is not zero. The pole of the function $P|_{C}$ can only occur in the places at infinity. So the some of all poles and zeros with multiplicities is nonpositive.
 By the definition of curve valuation, this number is $\sum_{i=1}^s(C_i\cdot l_{\infty})v_{C_i}(P)>0$ which is a contradiction.
 It follows that $P|_C=0.$
\endproof

\smallskip
\subsubsection*{\textbf{ 2.2.1)The case  $v_{C^j_1}, v_{C^j_2}\in W$ for all $j\leq 0$.}}

%

If $q^0_1, q^0_2$ are fixed by $f|_E.$  There exists a neighborhood $W'$ of $v_*$ such that $f_{\d}(W')\cap U(q^0_i)\subseteq W'\cap U(q^0_i)$ and $v_{C^0_i}\not\in W'$ for $i=1,2$. Since $\deg f|_E=1$ and $v_{C^j_i}\in W$ for all $j\leq 0$, $i=1,2$, we have $q^j_i=q^0_i$ for all $j\leq 0.$ It follows that $v_{C^j_i}\in U(q^0_i)$ and $v_{C^j_i}\not\in W'$ for all $j\leq 0$, $i=1,2$.
By Lemma \ref{lemexwtowcninw}, we conclude our theorem in this case.

If $q^0_1$ is fixed by $f|_E$ and $q^0_2$ is not fixed by $f|_E$, then $q^0_2$ is not periodic.
There exists a nice neighborhood $W'$ of $v_*$ such that $W'\subseteq W$, $f_{\d}(W')\cap U(q^0_1)\subseteq W'\cap U(q^0_1)$ and $v_{C^0_1}\not\in W'$. Since $\deg f|_E=1$ and $v_{C^j_1}\in W$ for all $j\leq 0$, we have $q^j_1=q^0_1$ for all $j\leq 0.$ It follows that $v_{C^j_1}\in U(q^0_1)$ and $v_{C^j_1}\not\in W'$ for all $j\leq 0$. By applying Lemma \ref{lemifvjinotinwthenbound} for $W'$, we may suppose that there exists $A'>0$ such that $m_1^j/m_2^j\leq A'$ for all $j\leq 0.$
Lemma \ref{lemconseqinwbound} implies that $\{m_2^j\}_{j\leq 0}$ is bounded and then $\deg(C^j)=m^j_1+m^j_2$ is bounded.
Then we conclude our theorem by Proposition \ref{procurveboundeddegfibration}.

Then we may suppose that both $q^0_1, q^0_2$ are not periodic by $f|_E.$ Lemma \ref{lemconseqinwbound} shows that $\{m_i^j\}_{j\leq 0}$ is bounded  for $i=1,2$ and then $\deg(C^j)=m^j_1+m^j_2$ is bounded. Then we conclude our theorem by Proposition \ref{procurveboundeddegfibration}.

\medskip
\subsubsection*{\textbf{2.2.2) The case  $v_{C^j_1}\in W$ for all $j\leq 0$ and there are infinitely many $j\leq 0$ such that $v_{C^j_2}\not\in W$.}}

\smallskip

If $q^0_1$ is fixed by $f|_E$, there exists a neighborhood $W'$ of $v_*$ such that $W'\subseteq W$, $f_{\d}(W')\cap U(q^0_1)\subseteq W'\cap U(q^0_1)$ and $v_{C^0_1}\not\in W'$. Since $\deg f|_E=1$ and $v_{C^j_1}\in W$ for all $j\leq 0$, we have $q^j_1=q^0_1$ for all $j\leq 0.$ It follows that $v_{C^j_1}\in U(q^0_1)$ and $v_{C^j_1}\not\in W'$ for all $j\leq 0$. More over there are infinitely many $j\leq 0$ such that $v_{C^j_2}\not\in W'$, then we conclude our theorem by Lemma \ref{lemexwtowcninw}.

So we may suppose that $q^0_1$ is not fixed by $f|_E$. Then Lemma \ref{lemconseqinwbound} shows that $m^j_1\leq 2m^0_1$ for all $j\leq 0$. By replacing $C$ by some $C^j$, $j\leq 0$, we may suppose that $v_{C^0_2}\not\in W$. For any $j\leq 0$ satisfying $v_{C^j_2}\not\in W$, we have $m^j_2/m^j_1<A$. It follows that $m^j_2<2Am^0_1$ when $v_{C^j_2}\not\in W$. For any $j\leq 0$ satisfying $v_{C^j_2}\in W$, there exists $j\leq j_1\leq -1$ such that $v_{C^{j'}_2}\in W$ for any $j\leq j'\leq j_1$ and $v_{C^{j_1+1}_2}\not\in W$. Since $f(U(x)\cap W)\subseteq U(x)\cap W$ for all $x\in E$ fixed by $f|_E$, $q^{j_1}_2$ is not fixed by $f|_E$.
By Lemma \ref{lemconseqinwbound}, we have $m^j_2\leq 2m^{j_1}_2$. By Lemma \ref{lemcjminkinotfarcji}, we have
$$m^{j_1}_2/m^{j_1}_1\leq (\deg(f))^{2}(m^{j_1+1}_2/m^{j_1+1}_1+1)<(\deg(f))^{2}(A+1).$$ It follows that $m^{j_1}_2\leq (\deg(f))^{2}(A+1)m^{j_1}_1\leq 2(\deg(f))^{2}(A+1)m^{0}_1.$  Then we have $m^j_2\leq 4(\deg(f))^{2}(A+1)m^{0}_1$. It follows that $\{m^j_2\}_{j\leq 0}$ is bounded. Then we have  $\{\deg(C^j)\}_{j\leq 0}$
is bounded and thus we conclude our theorem by Proposition \ref{procurveboundeddegfibration}.

\medskip

\subsubsection*{\textbf{ 2.2.3) The case that for all $i=1,2$, there are infinitely $j_1\leq 0$ such that $v_{C^{j_1}_i}\not\in W$.}}
%
%
%
%
Denote by $S_i$ the set of $j\in \Z_{\leq 0}$ such that $v_{C^j_i}\in W$ for all $i=1,2.$ It follows that $S_1\cup S_2=\Z_{\leq 0}$. In this case we  have
$\Z_{\leq 0}\setminus S_i$ is infinite for all $i=1,2$.

If $S_1$ is finite, then there exists $N'\leq 0$ such that $\{N',N'-1,\cdots\}\subseteq S_2$. It follows that $\Z_{\leq 0}\setminus S_i$ is finite
 which contradicts to our assumption. So $S_1$ is infinite. The same, $S_2$ is infinite also.

 There exists $N_0\leq 0$ such that $\{0,-1,\cdots,N_0\}\cap S_i\neq\emptyset$ for all $i=1,2$.
%
%

For any $n\geq 0$, denote by $O_n$ the set of points $x\in E$ such that $U(x)\cap (\cap_{k=0}^{n} f^{-k}_{\d}(W))=U(x)\cap W$ and $U(x)\cap f^{-n-1}_{\d}(W)\neq U(x)\cap W$. Observe that $O_0$ is finite. Since $O_n=f|_E^{-n}(O_0)$ for all $n\geq 0$, $O_n$ is finite. There are no periodic points in $O_0$, which implies that for any finite subset $B$, $O_n\cap B=\emptyset$ for $n$ large enough.

Set $M:=\min\{-8A-16A^2,-8\deg(f)A/\delta-16\deg(f)^2A^2/\delta^2,-288\}-1$ and let $N_1$ be defined in the following
\begin{lem}\label{lemfarindealpverylar}For any $M\leq 0$, there exists $N_1\geq 0$ such that for all $x\in O_{N_1}$, we have
$\{v\in U(x)|\,\, \alpha(v)\geq M\}\subseteq U(x)\cap f^{-N_1}_{\d}(W)$.
\end{lem}

The following lemma allows us to suppose that for all $i=1,2$ and $j\leq N_0$, if $\{j,j+1,\cdots,j+N_1\}\subseteq S_i$ and $j+N_1+1\not\in S_i$,  we have $m_i^j< (1-M)^{-1/2}\deg (C^j)$.

\begin{lem}\label{lematomnonenotlarge}If there are infinitely many $j\leq 0$ such that $\{j,j+1,\cdots,j+N_1\}\subseteq S_i$, $j+N_1+1\not\in S_i$ and $m_i^j\geq (1-M)^{-1/2}\deg (C^j)$ for some $i=1,2$, then $C$ is periodic.
\end{lem}

\begin{lem}\label{lemvcjonetwoinlomnt} If there are infinitely many $j\leq N_0$ satisfying $\{j,\cdots,j+N_1\}\subseteq S_1\cap S_2$ for all $i=1,2$, then the pair $(\A^2,f)$ satisfies the DML property for $C$.
\end{lem}

By Lemma \ref{lemexwtowcninw}, there exists an infinite sequence $\{j_1>j_2\cdots\}$ such that for all $l\geq 1$, $\{j_l,j_l+1,\cdots,j_l+N_1\}\in S_i$ for some $i=1,2.$
We may suppose that
$\{j_l,j_l+1,\cdots,j_l+N_1\}\in S_1$. For all $l\geq 0$, there exists $n_l\geq 0$ such that  $\{j_l+n_l,j_l+n_l+1,\cdots,j_l+n_l+N_1\}\in S_1$ but $j_l+n_l+N_1+1\not\in S_1$. By Lemma \ref{lemvcjonetwoinlomnt}, we may suppose that $\{j_l+n_l,j_l+n_l+1,\cdots,j_l+n_l+N_1\}\not\in S_2$. It follows that both $v_{C^{j_l+n_l}_2}$ and $v_{C^{j_l+n_l}_1}$ are not contained in $\cap_{k=0}^{N_1+1}f^{-k}_{\d}(W)$ for all $l\geq 0$.

Since $\Z_{\leq 0}\setminus S_1$ is infinite, we may suppose that for all $l\geq 1$, $\{j_{l+1},j_{l+1}+1,\cdots,j_{l}\}\not\subseteq S$. It follows that $j_{l+1}+n_{l+1}<j_{l}<j_{l}+n_l$. Then we have $j_1+n_1>j_2+n_2\cdots$. Then we conclude our Theorem by Lemma \ref{lemexwtowcninw}.

\medskip

\proof[Proof of Lemma \ref{lemfarindealpverylar}]By Proposition \ref{corthetacontinuousgeqm}, the set $W_{M-1}:=W(\theta^*)\cap \{v\in V_{\infty}|\,\,\alpha(v)\geq M-1\}$ is compact. Then there are finitely many open set $U_i$, $i=1,\cdots,l$, taking form $U_i=\{v\in V_{\infty}|\,\, v>v_i\}$ not containing $v_*$ such that $W_{M-1}\subseteq \cup_{i=1}^mU_i$.
Proposition \ref{prothetageztendstovs} shows that there exists $N_2>0$ such that for all $n\geq N_2$, $f^{n}_{\d}(\{v\in V_{\infty}|\,\,\alpha(v)\geq M\}\setminus(\cup_{i=1}^mU_i))\subseteq W$ for all $n\geq N_2.$ It follows that $\{v\in V_{\infty}|\,\,\alpha(v)\geq M\}\setminus(\cup_{i=1}^mU_i)\subseteq f^{-n}_{\d}(W)$ for all $n\geq N_2.$

Denote by $B$ the set of points in $E$ determinate by the direction defined by $[v_*, v_i].$ There exists $N_1\geq N_2$, such that $O_{N_1}\cap B=\emptyset.$
 Then we have $$U(x)\cap (\cup_{i=1}^mU_i)=\emptyset .$$ It follows that $\{v\in U(x)|\,\, \alpha(v)\geq M\}\subseteq U(x)\cap f^{-N_1}_{\d}(W)$.
\endproof

\proof[Proof of Lemma \ref{lematomnonenotlarge}]
We suppose that there are infinitely many $j\leq N_0$ such that $\{j,j+1,\cdots,j+N_1\}\subseteq S_1$, $j+N_1+1\not\in S_1$ and $m_1^j\geq (1-M)^{-1/2}\deg (C^j)$. Then we have $q^j_1\in O_{N_1}$ and $v_{C^{j}_1}\in V_{\infty}\setminus f^{-N_1-1}_{\d}(W)$. By Lemma \ref{lemfarindealpverylar}, $U(q^j_1)\cap\{v\in V_{\infty}|\,\,\alpha(v)\geq M\}\subseteq f^{-N_1-1}_{\d}(W).$

Since $O_{N_1}$ is finite, there exists $t\in O_{N_1}$ such that there exists a sequence of branches $\{C^{j_n}_{1}\}_{n\geq 0}$, $0\geq j_0>j_1>\cdots$ such that $q^{j_n}_{1}=t$ and $m_{1}^{j_n}\geq (1-M)^{-1/2}$.
The boundary $\partial( U(t)\cap f^{-N_1-1}_{\d}(W))$ of $U(t)\cap f^{-N_1-1}_{\d}(W)$ is finite and for all $v\in (\partial (U(t)\cap f^{-N_1-1}_{\d}(W))\setminus \{v_*\}$, we have $\alpha(v)<M$. Since $v_{C^{j_n}_{1}}\in U(t)\setminus f^{-N_1-1}_{\d}(W)$, there exists $v_n\in (\partial (U(t)\cap f^{-N_1-1}_{\d}(W))\setminus \{v_*\}$ satisfying $v_n<v_{C^{j_n}_{1}}$. Then there exists $n_1>n_2>0$ such that $v_{n_1}=v_{n_2}.$ If $v_{C^{j_{n_1}}_{1}}\neq v_{C^{j_{n_2}}_{1}}$, then we have
$$\deg(C^{j_{n_1}}_{1})\deg(C^{j_{n_2}}_{1})=(C^{j_{n_1}}_{1}\cdot C^{j_{n_2}}_{1})
\geq m^{j_{n_1}}_{1}m^{j_{n_2}}_{1}(1-\alpha(v_{C^{j_{n_1}}_{1}}\wedge v_{C^{j_{n_2}}_{1}}))$$$$
>(1-M)^{-1}\deg(C^{j_{n_1}}_{1})\deg(C^{j_{n_2}}_{1})(1-M)=\deg(C^{j_{n_1}}_{1})\deg(C^{j_{n_2}}_{1})$$
which is a contradiction.
It follows that $v_{C^{j_{n_1}}_{1}}= v_{C^{j_{n_2}}_{1}}$ and then $C$ is periodic.
\endproof

\proof[Proof of Lemma \ref{lemvcjonetwoinlomnt}]
Suppose that there are infinitely many $j\leq N_0$ such that $\{j,\cdots,j+N_1\}\subseteq S_1\cap S_2$ for all $i=1,2$. There exists a unique $n_i\geq 0$, such that $\{j+n_i,\cdots,j+n_i+N_1\}\subseteq S_i$ but $j+n_i+N_1+1\not\in S_i$. Then we have
 $q^{j+n_i}_i\in O_{N_1}$ for all $i=1,2$. We may suppose that $n_1\leq n_2$. Since for each $i=1,2$, there are infinitely many $j\leq 0$ such that $v_{C^j_i}\not\in W$, we may suppose that $n_1+j\leq n_2+j\leq N_0$.

We first suppose that $m^{j}_1/m^{j}_2\geq 4((1-M)^{1/2}-1)^{-1}$. By Lemma \ref{lemvonetwoinwbounded}, we have $m^{j+n_1}_1/m^{j+n_1}_2\geq((1-M)^{1/2}-1)^{-1}$ and then we have  $m^{j+n_1}_1\geq (1-M)^{-1/2}\deg(C^{j+n_1}).$ Since $v_{C^{j+n_1}_1}\in W$ and $q^{j+n_1}_1\in O_{N_1}$, it contracts our assumption above Lemma \ref{lematomnonenotlarge}.

Then we may suppose that $m^{j}_1/m^{j}_2< 4((1-M)^{1/2}-1)^{-1}$.

We claim the following
\begin{lem}\label{lemvonetwoinwbounded}
 If there are $j_1<j_0\leq N_0$, such that $v_{C^{j_0}_i},\cdots,v_{C^{j_1}_i}\in W$ for $i=1,2$ and the centers $q^{j_1}_i$, $i=1,2$ are not periodic, then we have $$4^{-1}m^{j_1}_1/m^{j_1}_2\leq m^{j_0}_1/m^{j_0}_2\leq 4m^{j_1}_1/m^{j_1}_2.$$
\end{lem}
If for all $\{j,\cdots, j+n_2\}\in S_1$, then by Lemma \ref{lemvonetwoinwbounded}, we have
$m^{j+n_2}_1/m^{j+n_2}_2< 16((1-M)^{1/2}-1)^{-1}.$ It follows that $$m^{j+n_2}_2> \left(1+16\left((1-M)^{1/2}-1\right)^{-1}\right)^{-1}\deg(C^{j+n_2})$$$$> \left(1+16\left((1+288)^{1/2}-1\right)^{-1}\right)^{-1}\deg(C^{j+n_2})=1/2\deg(C^{j+n_2})$$
$$>17^{-1}\deg(C^{j+n_2})\geq(1-M)^{-1/2}\deg(C^{j+n_2}) .$$
Since $v_{C^{j+n_2}_2}\in W$ and $q^{j+n_2}_2\in O_{N_1}$, it contracts our assumption.

Then we have $n_2\geq n_1+1$ and the set $Y:=\{j+n_1+1,\cdots, j+n_2\}\setminus S_1$ is not empty.

If $v_{C^{j+n_2}_1}\not\in W$, then $m^{j+n_2}_1/m^{j+n_2}_2< A$. It follows that $m^{j+n_2}_1< A/(1+A)\deg(C^{j+n_2}).$ Then we have $$\deg(C^{j+n_2})=m_1^{j+n_2}+m^{j+n_2}_2< \left((1-M)^{-1/2}+A/(1+A)\right)\deg(C^{j+n_2})<\deg(C^{j+n_2})$$ which is a contradiction.

Then we have $j+n_2\not\in Y$. Denote by $n'$ be the maximal number satisfying $j+n'\in Y$. Since $v_{C^{j+n'}_1}\not\in W$, we have
 $m^{j+n'}_1/m^{j+n'}_2< A$. Since $v_{C^{j+n'}_1}\in f^{-1}_{\d}(W)$, we have $d(f,v_{C^{j+n'}_1})>\delta$. Then we have $$m^{j+n'+1}_1/m^{j+n'+1}_2=(d(f,v_{C^{j+n'}_2})/d(f,v_{C^{j+n'}_1}))(m^{j+n'}_1/m^{j+n'}_2)< \deg(f)A/\delta.$$
By Lemma \ref{lemvonetwoinwbounded}, we have $m^{j+n_2}_1/m^{j+n_2}_2< 4\deg(f)A/\delta$ and then
 $m^{j+n_2}_1< (4\deg(f)A/(\delta+4\deg(f)A)\deg(C^{j+n_2}).$ Then we have $$\deg(C^{j+n_2})=m_1^{j+n_2}+m^{j+n_2}_2$$$$< \left((1-M)^{-1/2}+(4\deg(f)A/(\delta+4\deg(f)A))\right)\deg(C^{j+n_2})<\deg(C^{j+n_2})$$ which is a contradiction.
\endproof

\proof[Proof of Lemma \ref{lemvonetwoinwbounded}]
Since $m^j_id(f,v_{C^j_i})=\deg(f|_{C^j})m^{j+1}_i$,
we have $$(m^j_1/m^j_2)/(m^{j+1}_1/m^{j+1}_2)=d(f,v_{C^j_2})/d(f,v_{C^j_1})$$ for all $j_1+1\leq j\leq j_0$.

It follows that $$(m^{j_0}_1/m^{j_0}_2)/(m^{j_1}_1/m^{j_1}_2)=(\prod_{j=j_0}^{j_1+1}d(f,v_{C^{j}_2}))/(\prod_{j=j_0}^{j_1+1}d(f,v_{C^{j}_1})).$$

When $q^j_i\not\in\{t_1,\cdots,t_d\}$, we have $d(f,v_{C^j_i})=\la_1$, for all $j_1+1\leq j\leq j_0$ and $i=1,2$; and
when $q^j_i\in\{t_1,\cdots,t_d\}$, we have $2^{-1/d}\la_1\leq d(f,v_{C^j_i})\leq 2^{1/d}\la_1$ for all $j_1+1\leq j\leq j_0$ and $i=1,2$.
Since $q^{j_0}_i,\cdots,q^{j_1}_i$ are distinct for $i=1,2$, we have $$2^{-1}\la_1^{j_0-j_1-1}\leq \prod_{j=j_0}^{j_1}d(f,v_{C^{j}_i})\leq 2\la_1^{j_0-j_1-1}$$
for $i=1,2$. Then we have $4^{-1}\leq (m^{j_0}_1/m^{j_0}_2)/(m^{j_1}_1/m^{j_1}_2)\leq 4.$
\endproof

\rem
When $s=1$, the proof is much easier than the case $s=2$.

Since $s=1$, we have $v_{C^j_1}\in W$ for all $j\leq 0$.

If $q^0_1$ is a fixed point, then there exists a neighborhood $W'$ of $v_*$ such that $f_{\d}(W')\cap U(q^0_1)\subseteq W'\cap U(q^0_1)$ and $v_{C^0_1}\not\in W'$. Since $\deg f|_E=1$ and $v_{C^j_1}\in W$ for all $j\leq 0$, we have $q^j_1=q^0_1$ for all $j\leq 0.$ It follows that $v_{C^j_1}\in U(q^0_1)$ and $v_{C^j_1}\not\in W'$ for all $j\leq 0$.
By Lemma \ref{lemexwtowcninw}, we conclude our theorem in this case.

If $q^0_1$ is not fixed, then it is not periodic  under $f|_{E}.$ Then any two points in $\{q^j_1\}_{j\leq 0}$ are distinct. By Lemma \ref{lemconseqinwbound}, $\deg (C^j)=m^j_1$ is bounded. Then we conclude by Proposition \ref{procurveboundeddegfibration}.
\endrem

\medskip

\subsubsection*{\textbf{ 2.3) The case $\deg f|_E\geq 2$}}
As in case 2.2), in this case we don't have a system of invariant neighborhood of $v_*$ in general. The key idea in this case is take advantage of the Northcott property. More precisely, since $f|_E$ is defined over a number field with degree at least $2$, then for any point $x\in E$, the set of inverse orbit in a fixed number field of a point $p\in E$ is finite. To do so, we first fix the notations.

There exists a number field $K$ such that $X,E,f,C,p$ are all defined over $K.$ For all $j\leq 0$, since $C^j$
contains infinitely many $K$-points, we have that $C^j$ is defined over $K.$ Since $C^j$ meets infinity at at most two points,  for $i=1,2$ and all $j\leq 0$, we have $[K(c(v_{C_i^j})):K]\leq 2$.

By Lemma \ref{lemexwtowcninw}, we may suppose that there exists $j_0\leq 0$ such that for all $j\leq j_0,$ there exists $i\in \{1,2\}$ satisfying $v_{C^j_i}\in W.$ By replacing $C$ by $C^{j_0}$, we may suppose that $j_0=0.$

\smallskip
\rem
When $s=1$, we are always in the following case 2.3.1) and the argument is the same as in the case $s=2$.
\endrem

\medskip
\subsubsection*{\textbf{ 2.3.1)
The case that there exists $j_0\leq 0$ for which $v_{C^j_i}\in W$ for all $i=1,2$ and $j\leq j_0$}}

\smallskip

By replacing $C$ by $C^{j_0}$, we may suppose that $j_0=0.$

By Northcott property, the set $\{c(v_{C_i^j}),j\leq 0\}$ is finite for all $i=1,2.$ It follows that $c(v_{C^0_i})$ is periodic for $i=1,2.$ By replace $f$ by some positive iterate, we may suppose that there exists $x_i\in E$ which is fixed by $f|_E$ and satisfying $x_i=c(v_{C^j_i})$ for all $j\leq 0.$ Let $W'$ be a neighborhood of $v_*$ in $W$ satisfying
\begin{points}
\item[$\d$]$v_{C^0_i}\not\in W'$ for $i=1,2$;
\item[$\d$]$f_{\d}(U(x_i)\cap W')\subseteq U(x_i)\cap W'$ for $i=1,2$.
\end{points}
It follows that $v_{C^j_i}\not\in W'$ for all $i=1,2$ and $j\leq 0.$ By Lemma \ref{lemexwtowcninw}, we conclude our theorem in this case.

%

\medskip

\subsubsection*{\textbf{2.3.2) The case that
there exists $i_0\in \{1,2\}$ and $j_0\leq 0$ such that $v_{C^j_{i_0}}\in W$ for all $j\leq j_0$}}

\smallskip
We may suppose that $i_0=1$ and by replacing $C$ by $C^{j_0}$, we may suppose that $j_0=0.$

By the argument in the $2.3.1)$, we may suppose that there exists an infinite set $S$ of index $j\leq 0$ such that $v_{C^j_2}\not\in W.$ By the same argument in $2.2.1)$, we may suppose that there exists $x\in E$ which is fixed by $f|_E$ and satisfying $x=c(v_{C^j_1})$ for all $j\leq 0.$ Let $W'$ be a neighborhood of $v_*$ in $W$ satisfying
\begin{points}
\item[$\d$]$v_{C^0_1}\not\in W'$;
\item[$\d$]$f_{\d}(U(x)\cap W')\subseteq U(x)\cap W'$.
\end{points}
It follows that $v_{C^j_1}\not\in W'$ for all $j\leq 0.$ By Lemma \ref{lemexwtowcninw}, we conclude our theorem in this case.

\medskip

 \subsubsection*{\textbf{2.3.3)The case that there exists $j_i\leq 0$ such that $v_{C^{j_i}_i}\not\in W$ for all $i=1,2$}}

 \smallskip

 Since $C^j$ is defined over $K,$ if there exists a point $x\in C^j\cap E$, we have $[K(x):K]\leq 2$ Let $P$
be the set of points $x\in E$ such that $f|_E^n(x)\in I(f)$ for some $n\geq 0$ and satisfying $[K(x):K]\leq 2.$ Observe that for all $x\in P$, $x$ is not periodic. By Northcott property, we have that $P$ is a finite set. Set $L:=\# P$.

 Pick $j_0=\min\{j_1,j_2\}-1$. It follows that for all $i=1,2$ and $j\leq j_0$, if $v_{C^j_i}\in W$, we have $c(v_{C^j_i})\in P.$ By replacing $C$ by $C^{j_0}$, we may suppose that $j_0=0.$

 If there exists $i\in\{1,2\}$ and $j\leq -L$, such that $v_{C^j_i}\in \cap_{k=0}^Lf_{\d}^{-k}(W)$, then we have $\{c(v_{C^{j}_i}),\cdots,c(v_{C^{j+L}_i})\}\subseteq P.$ Since there are not periodic points in $P$, we get a contradiction. It follows that
 $v_{C^j_i}\not\in \cap_{k=0}^Lf_{\d}^{-k}(W)$ for all $j\leq -L$ and $i=1,2$. Then we conclude by Lemma \ref{lemexwtowcninw}.
%
%
%
%
 %
%
%
\endproof

%

\section{The case $\la_1^2>\la_2$ and $\#J(f)\leq 2$}
The aim in this section is to prove Theorem \ref{thmdmlpoly} in the only case left:

Let $f:\mathbb{A}^2_{\overline{\mathbb{Q}}}\rightarrow \mathbb{A}^2_{\overline{\mathbb{Q}}}$ be a dominant endomorphism defined over $\mathbb{\overline{Q}}$ satisfying $\la_1^2>\la_2$ and
$\#J(f)\leq 2$, then we have the following

\begin{thm}\label{thmsupponetwodml}
Let $C$ be an irreducible curve in $\mathbb{A}^2_{\overline{\mathbb{Q}}}$ and $p$ be a closed point in $\mathbb{A}^2_{\overline{\mathbb{Q}}}.$ If the set $\{n\in \mathbb{N}|\,\,f^n(p)\in C\}$ is infinite, then we have either $C$ is periodic or $p$ is preperiodic.
\end{thm}

Write $\theta^*=\sum_{i=1}^sr_iZ_{v_i}$ where $r_i>0$ for $i=1,\cdots,s$, $\sum_{j=1}^{s}r_i\alpha(v_i\wedge v_j)=0$ and $\sum_{i=1}^sr_i=1.$
Further by Proposition \ref{proreducetofix}, we suppose that $f_{\d}(v_i)=v_i$ and $d(f,v_i)=\la_2 /\la_1$ for all $i=1,\cdots,s.$
Let $w^i_1$ and $w_2^i$ for $i=1,\cdots,s$ be valuations defined as in Proposition \ref{prosegmentrep}.

To prove Theorem \ref{thmsupponetwodml}, we need a some new techniques.
In Section \ref{subsectiondgreen}, we introduce the $D$-Green functions for all $\mathbb{R}$-divisor $D$ in $C(\mathfrak{X})$. Then in Section \ref{secanatt} we use these $D$-Green functions to contracts an attracting set in $\A^2$. At last we prove Theorem \ref{thmsupponetwodml} in Section \ref{subsectionproofend}.

\subsection{The $D$-Green functions on $\mathbb{A}^2$}\label{subsectiondgreen}
In this section, $k$ is an algebraically closed field with a nontrivial absolute value $|\cdot|_v$.

\begin{prodefi}\label{prodefidgfunction}Let $D$ be a $\mathbb{R}$-divisor in $C(\mathfrak{X})$. Let $X\in \mathcal{C}$ be a compactification of $\mathbb{A}^2_k$ such that $D$ can be realized as a $\mathbb{R}$-divisor supposed by $X_{\infty}:=X\setminus \A^2.$ A function $\phi:\mathbb{A}^2_k\rightarrow\mathbb{R}$ is said to be a $D-$Green function\index{$D-$Green function} if it is continuous with respect to the the topology induced by $|\cdot|_v$ and there exists a finite set of local coordinate chars $\{U_i\}_{1\leq i\leq l}$ with respect to the topology induced by $|\cdot|_v$ such that
\begin{points}
\item $X_{\infty}\subseteq\cup_{i=1}^lU_i;$
\item for any $i=1,\cdots,l$, $X_{\infty}\cap U_i$ is defined by $x=0$ or $xy=0;$
\item for any $i=1,\cdots,l$, there exists a real number $C_i\geq 0$ such that in $U_i\cap \mathbb{A}^2_k$ we have $$-\ord_{\{x=0\}}D\log|x|_v-C_i\leq\phi\leq -\ord_{\{x=0\}}D\log|x|_v+C_i$$ if $X_{\infty}\cap U_i$ is defined by $x=0$ and $$-\ord_{\{x=0\}}D\log|x|_v-\ord_{\{y=0\}}D\log|y|_v-C_i\leq\phi\leq -\ord_{\{x=0\}}D\log|x|_v-\ord_{\{y=0\}}D\log|y|_v+C_i$$ if $X_{\infty}\cap U_i$ is defined by $xy=0.$
\end{points}
This definition does not depend on the choice of the compactification $X.$
\end{prodefi}
\proof[Proof of Proposition-Definition \ref{prodefidgfunction}]We only have to check that this definition is stable under blowup one point at infinity.

Let $\phi$ be a function on $\A^2$ satisfying the conditions in  Proposition \ref{prodefidgfunction} and $q$ be any point in $X_{\infty}$.

There exists a local coordinate chars $U_i$ of $X$, such that $q\in U_i$. We may suppose that in this coordinate $q=(0,0)$ and $D|_{U_i}$ is defined by $aD_x+bD_y$ where $D_x,D_y$ are divisors of $U_i$ defined by $x=0$ and $y=0$. We may suppose that $D_x$ is contained in $X_{\infty}$. Observe that if $D_y$ is not contained in $X_{\infty}$, then we have $b=0$. Denote by $\pi:Y\to X$ the blowup of $X$ at $q$. We may cover $\pi^{-1}(U_i)$ by two open set $V_1$ and $V_2$ such that $\pi|_{V_1}:(x,y)\to (x,xy)$ and $\pi_{V_2}:(x,y)\to (xy,y)$. Then the Cartier divisor $D$ on $V_1$,$V_2$ takes form $\pi^*D|_{V_1}=(a+b)D_x+bD_y$ and $\pi^*D|_{V_2}=aD_x+(a+b)D_y$. By (iii), we have
$$-a\log|x|_v-bD\log|xy|_v-C_i\leq\phi\circ\pi|_{V_1}\leq -a\log|x|_v-b\log|xy|_v+C_i$$ and $$-a\log|xy|_v-bD\log|y|_v-C_i\leq\phi\circ\pi|_{V_2}\leq -a\log|xy|_v-b\log|y|_v+C_i.$$ Thus we have $$-(a+b)\log|x|_v-bD\log|y|_v-C_i\leq\phi\circ\pi|_{V_1}\leq -(a+b)\log|x|_v-b\log|y|_v+C_i$$ and $$-a\log|x|_v-(a+b)D\log|y|_v-C_i\leq\phi\circ\pi|_{V_2}\leq -a\log|x|_v-(a+b)\log|y|_v+C_i$$ which concludes our proof.
\endproof

Then we have the following basic properties for $D$-Green functions.
\begin{pro}\label{probasicgreenfunction}We have the following properties.
\begin{points}
\item The function $\phi=0$ is a $0$-Green function.
\item Let $D_1,D_2$ be two $\mathbb{R}$-divisors in $C(\mathfrak{X})$. For $i=1,2$, let $\phi_i$ be a $D_i$-Green function on $\mathbb{A}^2_k.$ Then $\phi_1+\phi_2$ is a $(D_1+D_2)$-Green function.
\item Let $D$ be a $\mathbb{R}$-divisor in $C(\mathfrak{X})$ and $\phi$ be a $D$-Green function on $\mathbb{A}^2_k.$ For any $r\in \mathbb{R}$, $r\phi$ is a $rD$-Green function.
\item Let $D$ be a $\mathbb{R}$-divisor in $C(\mathfrak{X})$. Let $\phi_1$ and $\phi_2$ be two $D$-Green functions on $\mathbb{A}^2_k.$ There exists $C\geq 0$ such that $-C\leq\phi_1-\phi_2\leq C.$
\item Let $f:\mathbb{A}^2_k\rightarrow \mathbb{A}^2_k$ be a dominant polynomial endomorphism on $\mathbb{A}^2_k$. Let $D$ be a $\mathbb{R}$-divisor in $C(\mathfrak{X})$ and $\phi$ be a $D$-Green function on $\mathbb{A}^2_k.$ We denote by $f^*D$ the pullback of $D$ as a Cartier class in $C(\mathfrak{X}).$
    Then $\phi\circ f$ is a $f^*D$-Green function.
\end{points}
\end{pro}
\proof[Proof of Definition-Proposition \ref{probasicgreenfunction}]Properties (i),(ii),(iii) and (iv) are directly from the definition of $D$-Green function. So we only need to prove (v). Pick a compactification $X\in \mathcal{C}$ satisfying the conditions in Definition-Proposition \ref{probasicgreenfunction}. Pick a compactification $Y\in \mathcal{C}$ satisfying the conditions in Definition-Proposition \ref{probasicgreenfunction} for $f^*D$ such that the morphism $f:\A^2_k\to \A^2_k$ extends to a morphism $f:Y\to X$. Let $\{U_i\}_{1\leq i\leq l}$ be a system of local coordinate charts satisfying the conditions in Definition-Proposition \ref{probasicgreenfunction}. For every $i=1,\cdots,l$, $D|_{U_i}$ is defined by $aD_x+bD_y$ where $D_x,D_y$ are divisors of $U_i$ defined by $x=0$ and $y=0$. It follows that $f^*D|_{f^{-1}(U_i)}=af^*D_{x}+bf^*D_{y}$. Let $\{V_i\}_{1\leq i\leq m}$ be a system of local coordinate charts of $Y$ for $f^*D$ satisfying the conditions in (i) and (ii) in Definition-Proposition \ref{probasicgreenfunction}. We may further suppose that for any $V_j$, there exists $U_i$ such that $V_j\subseteq f^{-1}(U_i)$. It follows that on $V_j$, we have
$$-\ord_{\{x=0\}}f^*D\log|x|_v-\ord_{\{y=0\}}f^*D\log|y|_v$$$$=-\ord_{\{x=0\}}(af^*D_{x}+bf^*D_{y})\log|x|_v-\ord_{\{y=0\}}(af^*D_{x}+bf^*D_{y})\log|y|_v$$
$$=-a\log|x^{\ord_{x=0}f^*D_x}y^{\ord_{y=0}f^*D_x}|_v-b\log|x^{\ord_{x=0}f^*D_y}y^{\ord_{y=0}f^*D_y}|_v.$$
Since $x\circ f/x^{\ord_{x=0}f^*D_x}y^{\ord_{y=0}f^*D_x}$ and  $y\circ f/x^{\ord_{x=0}f^*D_y}y^{\ord_{y=0}f^*D_y}$ have no zero in $V_j$, we have
$$-a\log|x^{\ord_{x=0}f^*D_x}y^{\ord_{y=0}f^*D_x}|_v-b\log|x^{\ord_{x=0}f^*D_y}y^{\ord_{y=0}f^*D_y}|_v$$$$=-a\log|x\circ f|_v-b\log|y\circ f|_v+O(1)=\phi\circ f+O(1).$$
\endproof

\begin{lem}\label{lemglobalg}Let $|\cdot|_v$ be a nontrivial absolute value of $k$. Let $X$ be a compactification of $\mathbb{A}^2_k$ in $\mathcal{C}$. Let $D$ be an effective divisor supposed by $X_{\infty}.$ We suppose that the line bundle $O_X(D)$ is generated by its global sections. Let $P_1,\cdots, P_s\in k[x,y]$ be a base of $H^0_X(D)$.
Let $\phi_D:X\rightarrow[0,+\infty]$ be a function on $X$ defined by  $\phi_D:=\log\max\{|P_1|_v,\cdots,|P_s|_v,1\}.$

 Then there exists a finite set of local coordinate chars $\{U_i\}_{1\leq i\leq l}$ with respect to the topology induced by $|\cdot|_v$ such that
\begin{points}
\item $X_{\infty}\subseteq \cup_{i=1}^lU_i;$
\item for any $i=1,\cdots,l$, $X_{\infty}\cap U_i$ is defined by $x=0$ or $xy=0;$
\item for any $i=1,\cdots,l$, there exists a real number $C_i\geq 0$ such that $$-\ord_{\{x=0\}}D\log|x|_v-C_i\leq\phi_D\leq -\ord_{\{x=0\}}D\log|x|_v+C_i$$ if $X_{\infty}\cap U_i$ is defined by $x=0$ and $$-\ord_{\{x=0\}}D\log|x|_v-\ord_{\{y=0\}}D\log|y|_v-C_i\leq\phi_D\leq -\ord_{\{x=0\}}D\log|x|_v-\ord_{\{y=0\}}D\log|y|_v+C_i$$ if $X_{\infty}\cap U_i$ is defined by $xy=0.$
\end{points}
In particular $\phi_D|_{\mathbb{A}^2_k}$ is a $D$-Green function.
\end{lem}
\proof[Proof of Lemma \ref{lemglobalg}]Since $|D|$ is base point free, there exists a finite set of local coordinate chars $\{U_i\}_{1\leq i\leq l}$ with respect to the topology induced by $|\cdot|_v$ such that
\begin{points}
\item[$\d$] $X_{\infty}\subseteq \cup_{i=1}^lU_i;$
\item[$\d$] for any $i=1,\cdots,l$, $X_{\infty}\cap U_i$ is defined by $x=0$ or $xy=0;$
\item[$\d$] for any $i=1,\cdots,l$, there exists $j_i\in\{1,\cdots,s\}$ such that $\Supp(\text{div}(P_{j_i})+D)\cap \overline{U_i}=\emptyset.$
\end{points}
In $U_i$, we have $P_j/P_{j_i}\in O(V_i)$ for some open set $V_i\supset \overline{U_i}$. It follows that $\phi_D\leq |\log P_{j_i}|+O(1)$ in $U_i.$
On the other hand, we have $|\phi_D|_v\geq |P_{j_i}|_v$. Then we have $\phi_D= |\log P_{j_i}|+O(1)$. Then we conclude our proposition by the fact that $\text{div}(P_{j_i})|_{U_i}=D|_{U_i}.$
\endproof

\begin{pro}\label{proexistenceofdgreenfunction}Let $D$ be a $\mathbb{R}$-divisor in $C(\mathfrak{X})$, up to a bounded function, there exists a unique $D$-Green function $\phi_D$ on $\mathbb{A}^2_k.$
\end{pro}
\proof[Proof of Proposition \ref{proexistenceofdgreenfunction}]The uniqueness is follows from (iv) of Proposition \ref{probasicgreenfunction}. So we only have to show the existence of the $D$-Green function.

Since $D$ is a $\mathbb{R}$-divisor in $C(\mathfrak{X})$, we may write it as a $\mathbb{R}$ combination of $\mathbb{Z}$ divisors in $C(\mathfrak{X})$. By (ii) and (iii) of Proposition \ref{probasicgreenfunction}, we may suppose that $D$ is a $\mathbb{Z}$ divisors.
Pick a compactification $X$ of $\mathbb{A}^2_k$ such that $D$ can be realized as a divisor supposed by $X_{\infty}.$ Pick two ample $\mathbb{Z}$-divisors $A_1$ and $A_2$ supported by $X_{\infty}$ such that $D=A_1-A_2.$ There exists a positive integer $l\geq 1$ such that for all $i=1,2$, $O_X(lA_i)$ is generated by its global sections. By Lemma \ref{lemglobalg}, for all $i=1,2$ there exists a $lA_i$-Green function $\phi_{A_i}$. Then we have $\phi_D:=l^{-1}(\phi_{lA_1}-\phi_{lA_2})$ is a $D$-Green function.
\endproof

\subsection{An attracting set}\label{secanatt}In this section, $k$ is an algebraically closed field with a nontrivial absolute value $|\cdot|_v$.
\smallskip
\begin{figure}
\centering
\includegraphics[width=10.5cm]{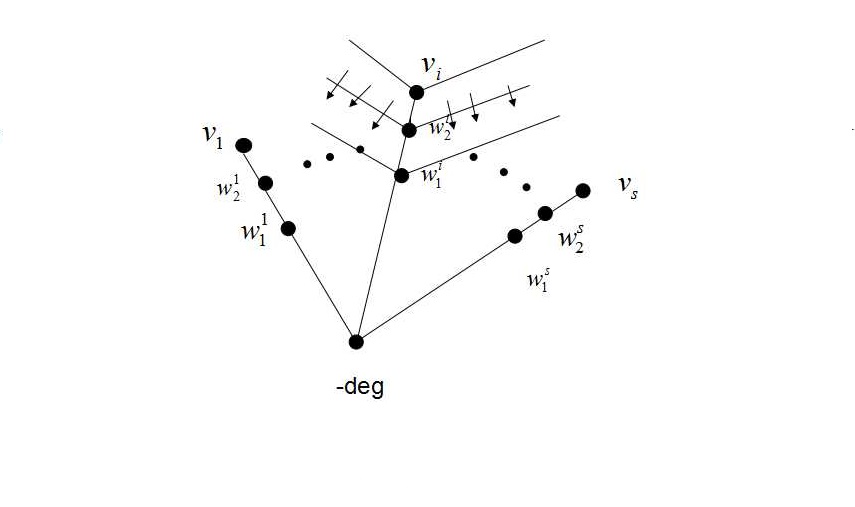}
\caption{}
\label{figc2}
\end{figure}

Let $f:\mathbb{A}^2_k\rightarrow \mathbb{A}^2_k$ be a dominant endomorphism defined over $k$ with $\la_1^2>\la_2$ and $\#J(f)\leq 2$.

we may suppose that $\theta^*=\sum_{i=1}^sr_iZ_{v_i}$ where $r_i>0$ for $i=1,\cdots,s$, $\sum_{j=1}^{s}r_i\alpha(v_i\wedge v_j)=0$ and $\sum_{i=1}^sr_i=1.$
Further by Proposition \ref{proreducetofix}, we suppose that $f_{\d}(v_i)=v_i$ and $d(f,v_i)=\la_2 /\la_1$ for all $i=1,\cdots,s.$

Recall Proposition \ref{prosegmentrep}.
For all $i=1,\cdots,s$, there are two valuations $w^i_1<w^i_2<v_i$ as in (1) of Figure \ref{figc2} such that
\begin{points}
\item $f^{-1}_{\d}(\{v\in V_{\infty}|\,\,w^i_1< v\wedge v_i<v_i\})=\{v\in V_{\infty}|\,\,w^i_2< v\wedge v_i<v_i\};$
\item $f_{\d}|_{\{v\in V_{\infty}|\,\,w^i_2< v\wedge v_i<v_i\}}$ is order-preserving;
\item for all valuation $w\in [w^i_1,v_i]$, $f^{-1}_{\d}(w)$ is one point in $[w^i_1,v_i];$
\item for all valuation $w\in\{v\in V_{\infty}|\,\,w^i_1< v\wedge v_i<v_i\}$, there exists $N\geq 1$ such that $f^n_{\d}(w)\in V_{\infty}\setminus \{v\in V_{\infty}|\,\,v\wedge v_i\geq w^i_1\}$ for all $n\geq N.$
\end{points}

 \begin{figure}
\centering
\includegraphics[width=10.5cm]{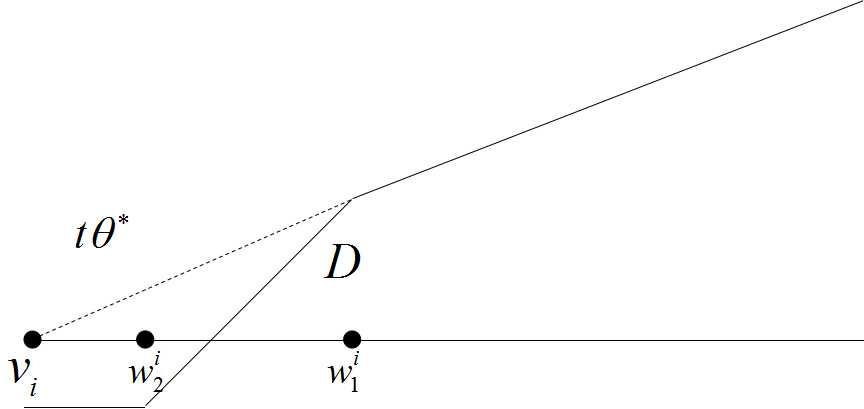}
\caption{}
\label{figc3}
\end{figure}

 Observe that $w^i_2\wedge w^j_2<w^i_2$ if $i\neq j$ and $f^{-1}(w_1^i)=\{w^i_2\}$ for $i=1,\cdots,s.$ Since for all $i=1\cdots,s$, we have $Z_{w^i_2}(w^i_2)<Z_{w^i_1}(w^i_2)$, there exists a positive rational number $t$ satisfying $(1+t)\sum_{i=1}^sr_iZ_{w_2}^i(w^j_2)-\sum_{i=1}^sr_iZ_{w_1}^i(w^j_2)<0$ for all $j=1,\cdots,s.$ Set $$D:=(1+t)\sum_{i=1}^sr_iZ_{w_2}^i-\sum_{i=1}^sr_iZ_{w_1}^i \text{ } \text{ } \eqno (A)$$  see Figue \ref{figc3}. Then $D$ can be viewed as a $\mathbb{R}$-divisor in $C(\mathfrak{X}).$ By Proposition \ref{proexistenceofdgreenfunction}, there exists a $D$-Green function $\phi_D$ on $\mathbb{A}^2_k.$ For any real number $C>0$, set $U_C:= \{p\in \mathbb{A}^2(k)|\,\, \phi_D(p)>C\}$\index{$U_C$} which is an open set of $\mathbb{A}^2_k.$

We have the following
\begin{pro}\label{proconattset}There are real numbers $C,C'>0$ such that for all $p\in U_C$, we have $$\phi_D(f(p))\geq \la_1\phi_D(p)-C'.$$

In particular, for any $B> \max\{C, C'/(\la_1-1)\}$, we have $f(U_B)\subseteq U_B$.
\end{pro}

\proof[Proof of Proposition \ref{proconattset}]Let $X$ be a compactification of $\mathbb{A}^2_k$ in $\mathcal{C}$, such that $D$ and $f^*D$ can be realized as a $\mathbb{R}$-divisor supposed by $X_{\infty}.$ By Proposition \ref{probasicgreenfunction}, $\phi_{f^*D}:=\phi_D\circ f$ a $f^*D$-Green function on $\mathbb{A}^2_k.$

By definition, there exists a finite set of local coordinate chars $\{U_i\}_{1\leq i\leq l}$ with respect to the topology induced by $|\cdot|_v$ such that
\begin{points}
\item $X_{\infty}\subseteq\cup_{i=1}^lU_i;$
\item for any $i=1,\cdots,l$, $X_{\infty}\cap U_i$ is defined by $x=0$ or $xy=0;$
\item for any $i=1,\cdots,l$, there exists a real number $C_i\geq 0$ such that in $U_i\cap \mathbb{A}^2_k$ we have
$$|\phi_D+\ord_{\{x=0\}}D\log|x|_v|\leq C_i,$$
$$|\phi_{f^*D} +\ord_{\{x=0\}}f^*D\log|x|_v|\leq C_i$$
if $X_{\infty}\cap U_i$ is defined by $x=0$;

and $$|\phi_D+\ord_{\{x=0\}}D\log|x|_v+\ord_{\{y=0\}}D\log|y|_v|\leq C_i,$$
$$|\phi_{f^*D}+\ord_{\{x=0\}}f^*D\log|x|_v+\ord_{\{y=0\}}f^*D\log|y|_v|\leq C_i$$
if $X_{\infty}\cap U_i$ is defined by $xy=0.$
\end{points}
For convince we set $\ord_{\{y=0\}}D=\ord_{\{y=0\}}f^*D:=0$ if $\{y=0\}$ is not contained in $X_{\infty}.$
Further, we may suppose that for all $i=1,\cdots,l$, we have $\max\{|x|_v,|y|_v\}<1$ for all points $(x,y)\in U_i.$

The set $X\setminus \cup_{i=1}^lU_i$ is compact in $\mathbb{A}^2_k$, so there exists $B'>0$ such that for all point $p\in X\setminus \cup_{i=1}^lU_i$, we have $\phi_D(p)<B'.$

We may suppose that there exists a $l'\in \{1,\cdots,l\}$ such that an index $i$ is contained in $\{1,\cdots,l'\}$ if and only if there exists an irreducible component $E$ of $X_{\infty}$ such that $\ord_E(D)=b_ED(v_E)>0$ and $E\cap U_i\neq \emptyset.$

For all index $i\geq l'+1$, we have $$\phi_D\leq -\ord_{\{x=0\}}D\log|x|_v-\ord_{\{y=0\}}D\log|y|_v+C_i\leq C_i.$$

Pick $C:=\max\{C_i\}_{1\leq i\leq l}+B'+1$, we have $U_C\subseteq \cup_{i=1}^lU_i$ and $U_C\cap U_i=\emptyset$ for all $i\in \{l'+1,\cdots,l\}.$ It follows that $U_C\subseteq \cup_{i=1}^{l'}U_i$.

Let $E$ be an exceptional divisor of $X$ satisfying $v_E\not\in B(\{w_2^1,\cdots,w_2^s\})^{\circ}$. Then we have $f_{\d}(v_E)\not\in B(\{w_1^1,\cdots,w_1^s\})^{\circ}.$
Then we have $$f^*D(v_E)=(1+t)\sum_{i=1}^sr_i(f^*Z_{w_2^i}\cdot Z_{v_E})-\sum_{i=1}^sr_i(f^*Z_{w_1^i}\cdot Z_{v_E})$$$$=(1+t)\sum_{i=1}^sr_i(Z_{w_2^i}\cdot f_*Z_{v_E})-\sum_{i=1}^sr_i(Z_{w_1^i}\cdot f_*Z_{v_E})$$$$=(1+t)\sum_{i=1}^sr_i(Z_{v_i}\cdot f_*Z_{v_E})-\sum_{i=1}^sr_i(Z_{v_i}^i\cdot f_*Z_{v_E})$$$$=t(f^*\theta^*\cdot Z_{v_E})=\la_1t(\theta^*\cdot Z_{v_E})$$$$=\la_1\left((1+t)\sum_{i=1}^sr_i(Z_{w_2^i}\cdot Z_{v_E})-\sum_{i=1}^sr_i(Z_{w_2^i}\cdot Z_{v_E})\right)$$$$=\la_1((1+t)\sum_{i=1}^sr_i(Z_{w_2^i}\cdot Z_{v_E})-\sum_{i=1}^sr_i(Z_{w_1^i}\cdot Z_{v_E})+\sum_{i=1}^sr_i(Z_{w_1^i}\cdot Z_{v_E})-\sum_{i=1}^sr_i(Z_{w_2^i}\cdot Z_{v_E}))$$$$=\la_1D(v_E)+\la_1(\sum_{i=1}^sr_i(Z_{w_1^i}\cdot Z_{v_E})-\sum_{i=1}^sr_i(Z_{w_2^i}\cdot Z_{v_E}))$$$$=\la_1D(v_E)+\la_1(\sum_{i=1}^sr_i(\alpha(w^i_1\wedge v_E)-\alpha(w^i_2\wedge v_E))).$$ Set $\psi(v_E):=\la_1(\sum_{i=1}^sr_i(\alpha(w^i_1\wedge v_E)-\alpha(w^i_2\wedge v_E))).$ We have $\psi(v_E)\geq 0$.

Fix an index $i\in \{1,\cdots,l'\}.$
If $U_i\cap X_{\infty}=\{x=0\}$, then we have $\ord_{\{x=0\}}D>0$. Set $E:=\{x=0\}$ in $U_i.$ It follows that $v_E\not\in B(\{w_2^1,\cdots,w_2^s\}).$
It follows that for all points $p$ in $U_i\cap \mathbb{A}^2_k$ with local coordinate $(x,y)$, we have $$\phi_D(f(p))=\phi_{f^*D}(p)\geq -\ord_{E}f^*D\log|x|_v-C_i=b_Ef^*D(v_E)\log|x|_v-C_i$$$$=-b_E(\la_1D(v_E)+\psi(v_E))\log|x|_v-C_i$$
$$=-\la_1\ord_{E}D\log|x|_v-b_E\psi(v_E)\log|x|_v-C_i$$$$\geq -\la_1\ord_{E}D\log|x|_v+\la_1C_i-(1+\la_1)C_i$$$$\geq \la_1\phi_D(p)-(1+\la_1)C_i.$$

If $U_i\cap X_{\infty}=\{xy=0\}$, we set $E_1:=\{x=0\}$ and $E_2:=\{y=0\}$ in $U_i.$ We may suppose that $D(v_{E_{1}})>0$ and then $v_{E_1}\not\in B(\{w_2^1,\cdots,w_2^s\})$. Since $w_2^i$'s are  valuations defined by an exceptional divisor in $X,$ we have $v_{E_2}\not\in B(\{w_2^1,\cdots,w_2^s\})^{\circ}.$ It follows that for all points $p$ in $U_i\cap \mathbb{A}^2_k$ with local coordinate $(x,y)$, we have $$\phi_D(f(p))=\phi_{f^*D}(p)\geq -\ord_{E_1}f^*D\log|x|_v-\ord_{E_2}f^*D\log|y|_v-C_i$$$$=-b_Ef^*D(v_{E_1})\log|x|_v-b_Ef^*D(v_{E_2})\log|y|_v-C_i$$
$$=-b_E(\la_1D(v_{E_1})+\psi(v_{E_1}))\log|x|_v-b_E(\la_1D(v_{E_2})+\psi(v_{E_2}))\log|y|_v-C_i$$
$$=-\la_1\ord_{E_1}D\log|x|_v-\la_1\ord_{E_2}D\log|y|_v-b_E\psi(v_{E_1})\log|x|_v-b_E\psi(v_{E_2})\log|y|_v-C_i$$$$\geq -\la_1\ord_{E_1}D\log|x|_v-\la_1\ord_{E_2}D\log|y|_v+\la_1C_i-(1+\la_1)C_i$$$$\geq \la_1\phi_D(p)-(1+\la_1)C_i.$$

Set $C':=\max\{(1+\la_1)C_i\}_{1\leq i\leq l'}$, we conclude that for all $p\in U_C$, we have $$\phi_D(f(p))\geq \la_1\phi_D(p)-C'.$$

For any $B> \max\{C, C'/(\la_1-1)\}$, we have $U_B\subseteq U_C$. Moreover, for any $p\in U_B$, we have $\phi_D(f(p))\geq \la_1B-C'>B$. It follows that $f(p)\in U_B$ and then $f(U_B)\subseteq U_B$.
\endproof

At last, we apply this attracting set $U_B$ to prove the following proposition which is an analogue of Theorem \ref{thmlodml} in our case.
%

\begin{pro}\label{prodmlcywt}Let $C$ be a curve in $\mathbb{A}^2_{\overline{\mathbb{Q}}}$ and $p$ be a closed point in $\mathbb{A}^2_{\overline{\mathbb{Q}}}.$ If there exists a branch $v_{C_1}$ of $C$ at infinity satisfying $v_{C_1}\in B(\{w^1_2,\cdots,w^s_2\})^\circ\setminus W(\theta^*)$ and the set $\{n\in \mathbb{N}|\,\,f^n(p)\in C\}$ is infinite, then $p$ is preperiodic.
\end{pro}

 \begin{figure}
\centering
\includegraphics[width=10.5cm]{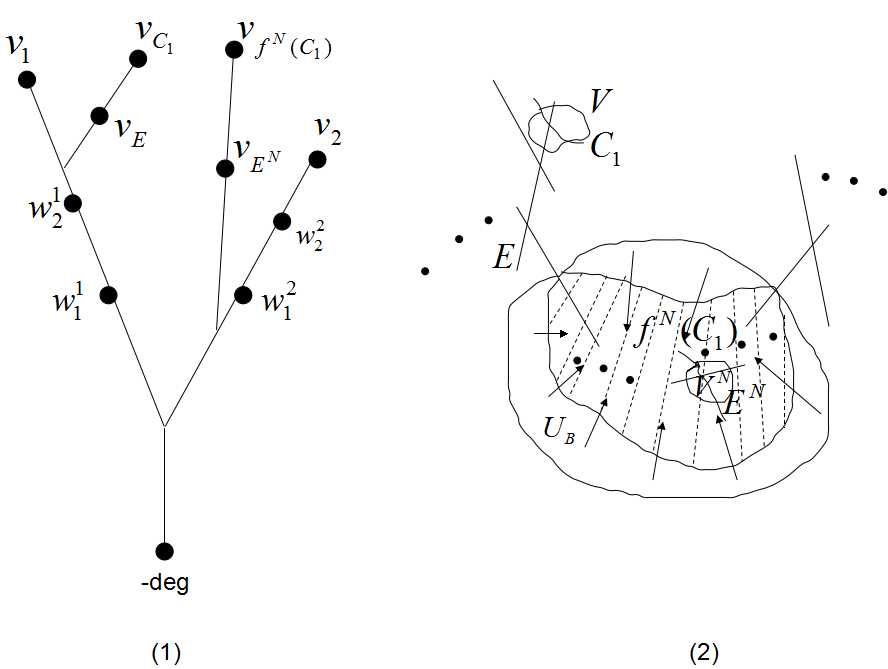}
\caption{}
\label{figc4}
\end{figure}

\proof[Proof of Proposition \ref{prodmlcywt}]Set $k:=\overline{\mathbb{Q}}$. We suppose that the set $\{n\in \mathbb{N}|\,\,f^n(p)\in C\}$ is infinite and $p$ is not preperiodic.
By Theorem \ref{thmsiegel}, we may suppose that $C$ is rational and has at most two places at infinity.  As in Equation (A), set $D:=(1+t)\sum_{i=1}^sr_iZ_{w_2}^i-\sum_{i=1}^sr_iZ_{w_1}^i$  where $t$ satisfies $(1+t)\sum_{i=1}^sr_iZ_{w_2}^i(w^j_2)-\sum_{i=1}^sr_iZ_{w_1}^i(w^j_2)<0$ for all $j=1,\cdots,s.$
There exists $N\geq 1$ such that $f^N_{\d}(v_{C_1})\not\in B(\{w^1_1,\cdots,w^s_1\})$ as in (1) of Figure \ref{figc4}.

Let $X$ be a compactification of $\mathbb{A}^2_k$, such that $D$ and $f^*D$ can be realized as a $\mathbb{R}$-divisor supposed by $X_{\infty}.$ Further, we may suppose that the center $c(v_{C_1})$ of $v_{C_1}$ is contained in a unique exceptional divisor $E$ and $c(v_{f^N(C_1)})$ is contained in a unique exceptional divisor $E^N.$ It follows that $\ord_E(D)<0$ and $\ord_{E^N}(D)>0.$ For convenience, we write $C$ for the Zariski closure of $C$ in $X.$

Let $\pi_1:\mathbb{P}^1_k\rightarrow C$ be a normalization. We may suppose that the branch $C_1$ is defined by the point $Q:=[1:0]\in \mathbb{P}^1_k .$ Set $\pi_2=f^N\circ \pi_1:\mathbb{P}^1_k\rightarrow f^N(C).$

Let $K$ be a number field such that $\pi_1$, $\pi_2$, $X$, $E$, $f$, $C$ and $p$ are all defined over $K.$

There exists a infinite sequence $\{n_1<n_2<\cdots\}$ of nonnegative integers such that $f^{n_i}(p)\in C.$ By contradiction, we suppose that $p$ is not preperiodic. We may suppose that for all $i\geq 1$, $C$ is smooth at $f^{n_i}(p)$. Set $p_i:=f^{n_i}(p)$ and $q_i:=\pi_1^{-1}(p_i).$ Write $q_i$ as form $[x_i:y_i].$ Observer that $q_i$ is a $K$ point in $\mathbb{P}^1_k$ for all $i\geq 0,$ then we may suppose that $x_i$'s and $y_i$'s and contained in $K.$  Let $S$ be a finite set of places $v\in \mathcal{M}_K$ containing $\mathcal{M}_K^{\infty}$ such that $f$ are defined in $O_S$ and $p$ is a $S$-integer. It follows that all $p_i$'s are $S$-integer. It follows that for all $v\in \mathcal{M}_K\setminus S$ there exists a number $C_v>0$ such that $|x_i/y_i|_v\leq C_v$ for all $i\geq 0$ and except a finite set of places, we have $C_v=1.$ By replacing $S$ by a bigger set, we may suppose that $C_v=1$ for all $v\in \mathcal{M}_K\setminus S$. By Northcott Property, we have $h_{\mathbb{P}^1_K}(q_i)\rightarrow \infty$ as $i\rightarrow \infty.$ Since $$h_{\mathbb{P}^1_K}(q_i)=\sum_{v\in \mathcal{M}_K}\log\max\{|x_i|_v,|y_i|_v\}$$$$=\sum_{v\in \mathcal{M}_K}\log\max\{|x_i/y_i|_v,1\}=\sum_{v\in S}\log\max\{|x_i/y_i|_v,1\},$$ there exists $v\in S$ such that by replacing $\{n_i\}_{i\geq 1}$ by a infinite subsequence, we have $|x_{i}/y_{i}|_v\rightarrow \infty$ as $i\rightarrow \infty.$ It follows that $q_i\rightarrow Q$ and $p_i\rightarrow c(v_{C_1})$ as $i\rightarrow \infty$ with respect to the topologies induced by $|\cdot|_v.$

Fix this place and by Proposition \ref{proexistenceofdgreenfunction}, there exists a $D$-Green function $\phi_D$ with respect to the topology induced by $|\cdot|_v.$ Since $E$ is the unique exceptional divisor containing the center $c(v_{C_1})$ of $v_{C_1}$ and $\ord_E(D)<0$, by the definition of $D$-Green function, there exists a neighborhood  $V$ of $c(v_{C_1})$ in $X$ with respect to the topology induced by $|\cdot|_v$ such that for all point $p'\in V\cap \mathbb{A}_k^2$, we have $\phi_D(p')<0.$

By Proposition \ref{proconattset}, there exists a real number $B>0$ such that $f(U_B)\subseteq U_B.$

%

Since $E^N$ is the unique exceptional divisor containing $c(v_{f^N(C_1)})$ and $\ord_{E^N}(D)>0$, by the definition of $D$-Green function, there exists a neighborhood  $V^N$ of $c(v_{f^N(C_1)})$ in $X$ with respect to the topology induced by $|\cdot|_v$ such that for all point $p'\in V^N\cap \mathbb{A}_k^2$, we have $\phi_D(p')>B.$ It follows that $V^N\cap \mathbb{A}^2_k\subseteq U_{B}.$ Since $q_i\rightarrow Q$ and $\pi_2(Q)=c(v_{f^N(C_1)})$, there exists $j>0$ such that $f^N(p_j)=\phi_2(q_j)\in V\cap \mathbb{A}^2_k\subseteq U_B.$ Then $f^r(p)\in U_B$ for all $r\geq n_j+N.$
Since $p_i\rightarrow c(v_{C_1})$ as $i\rightarrow \infty$, there exists $n_i\geq n_j+N$ such that $f^{n_i}(p)=p_i\in V.$ This contradicts the fact that $V\cap U_B=\emptyset$ and then we conclude our proposition.
\endproof

\subsection{Proof of Theorem \ref{thmsupponetwodml}}\label{subsectionproofend}
Set $k=\overline{\mathbb{Q}}$. We may suppose that $C$ can not be contracted to a point by $f^n$ for any $n\geq 0$ and $p$ is not preperiodic.
By Theorem \cite[Theorem 1.3]{Bell2010}, we may suppose that $Jf$ is not a constant. By Theorem \ref{thmsiegel}, we may suppose that $C$ has at most $2$ places at infinity.
Let $C_1,\cdots,C_t$, $t\in\{1,2\}$ be all branches of $C$ at infinite.

In the rest of this section we present our proof in the situation $t=2$ and we will give a remark for the situation $t=1$ in every case.

\medskip

\subsubsection*{\textbf{1) The case that $v_{C_i}\in W(\theta^*)$ for all branches $C_i$ of $C$ at infinity}}

\smallskip
\rem In the case $t=1$, we can use the same argument as in the case $t=2$.
\endrem

\subsubsection*{\textbf{1.1)}}  If $v_i$ is divisorial for all $i=1,\cdots,s$, by Theorem\ref{thmsuppfinitdivcarzfi} we have $R_{\{v_1,\cdots,v_s\}}=k[P]$ where $P$ is a polynomial in $k[x,y]\setminus k.$ Since $v_{C_i}\in W(\theta^*)=B(\{v_1,\cdots,v_s\})$ for all branches $C_i$ of $C$ at infinity, there exists $j_i\in \{1,\cdots,s\}$ such that $v_{j_i}<v_{C_i}.$

We have $v_{C_i}(P)\geq v_{j_i}(P)\geq0$ for all $i\in 1,2.$ Then the function $P|_{C}$ has no poles. It follows that $P|_{C}$ is a constant in $k.$ Then there exists an element $r\in \mathbb{A}^1(k)$ such that $C$ is contained in the fiber of $P:\mathbb{A}^2_k\rightarrow \mathbb{A}^1_k$ above $r.$ Pick a polynomial morphism $G:\mathbb{A}^1_k\rightarrow \mathbb{A}^1_k$ as in Theorem \ref{thmsuppfinitdivcarzfi}. Since $P\circ f^n=G^n\circ P$ for all $n\geq 0$ and the set $\{n\in \mathbb{N}|\,\,f^n(p)\in C\}$ is infinite, we have that $r$ is periodic under $G.$ Since $\{P-r=0\}$ has only finitely many irreducible component, then $C$ is periodic.

\smallskip

\subsubsection*{\textbf{1.2)}}Then we suppose that $v_1$ is nondivisorial.

\subsubsection*{\textbf{ 1.2.1)}}
\begin{figure}
\centering
\includegraphics[width=6.5cm]{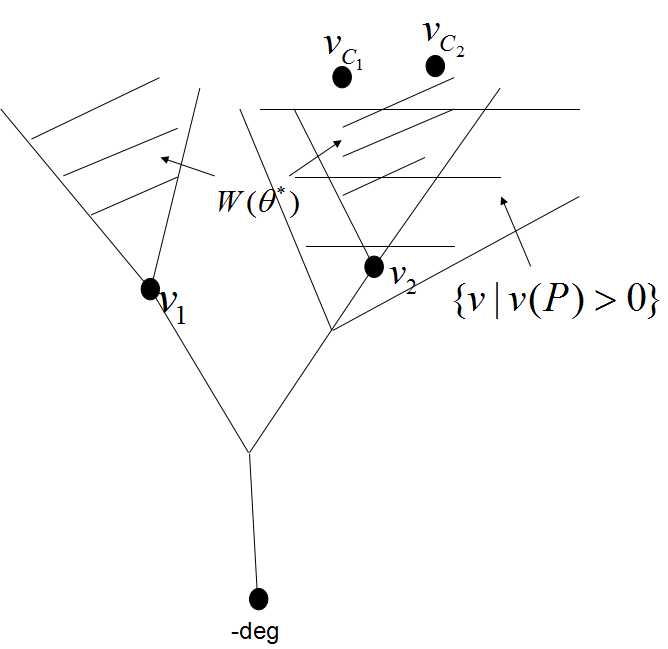}
\caption{}
\label{figc5}
\end{figure}
If for all $i=1,2$, we have $v_{C_i}\not\in B(\{v_1\})$, then we have $s=2$ and $v_{C_i}>v_2$ for all $i=1,2.$
See Figure \ref{figc5}.

Set $\psi:=R_{[-\deg,v_2]}\theta^*=r_1Z_{v_1\wedge v_2}+r_2Z_{v_2}\in \SH^+(V_{\infty}).$ Then we have $\psi(v)=0$ for all $v\geq v_2$ and $\langle\psi,\psi\rangle>0.$ By Theorem \ref{prodimtwothengeqz}, there exists a polynomial $P\in k[x,y]\setminus k$ satisfying $v_2(P)>0.$ Then we have $v_{C_i}(P)\geq v_{2}(P)> 0$ for all $i\in 1,2.$ It follows that $C$ is an irreducible component of $\{P=0\}$.
Apply the same argument for $f^n(C)$, $n\geq 0$, we have that $f^n(C)$ is an irreducible component of $\{P=0\}$. It follows that $C$ is periodic.

\smallskip

\subsubsection*{\textbf{ 1.2.2)}}
\begin{figure}
\centering
\includegraphics[width=6.5cm]{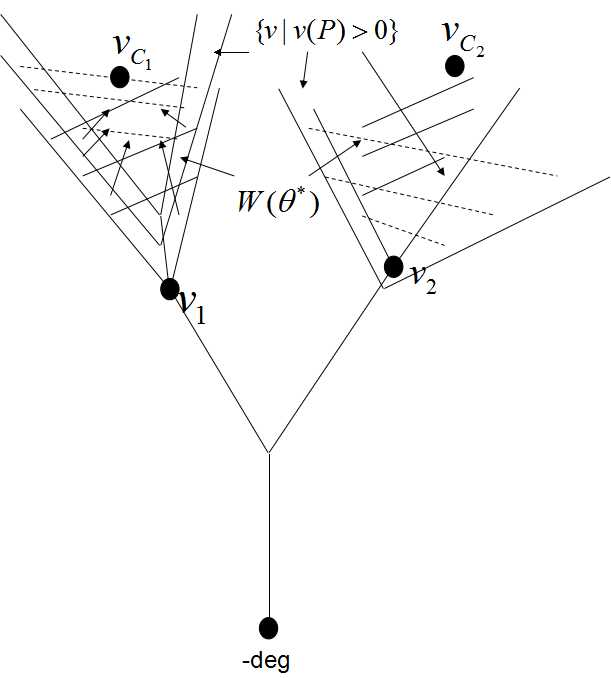}
\caption{}
\label{figc6}
\end{figure}

Otherwise, we may suppose that $v_{C_1}>v_1.$ It follows that $v_1$ is an irrational valuation. See Figure \ref{figc6}.

By Proposition \ref{prosegmentrep}, there exists $w>v_1$ such that for all $v\in \{v\in V_{\infty}|\,\,v_1<v\wedge w<w\}$ there exist $N_v\geq 0$ such that either $d(f^n,v)=0$ or $f^n_{\d}(v)\geq w$ for all $n\geq N_v.$
Pick a valuation $w_1\in (v_1,w)$ and apply \cite[Proposition 3.22]{Xieb} for $\{w_1\}\cup (\{v_1,\cdots,v_s\}\setminus \{v_1\})$. There exists a function $\psi\in \SH^+(V_{\infty})$ such that
$\psi(v)=0$ for all $v\in B(\{w_1\}\cup (\{v_1,\cdots,v_s\}\setminus \{v_1\}))$ and $\langle\psi,\psi\rangle>0.$
There exists $N\geq 0$ such that for all $n\geq N$, we have either $d(f^n,v_{C_i})=0$ or $f^N_{\d}(v_{C_i})\in W(\theta^*)\setminus\{v\in V_{\infty}|\,\,v_1\leq v\wedge w<w_1\}= B(\{w_1\}\cup (\{v_1,\cdots,v_s\}\setminus \{v_1\})).$ By replacing $C$ by some positive iterate, we may suppose that $d(f^n,v_{C_i})\neq 0$  and $f^n_{\d}(v_{C_i})\in B(\{v_2,w_1\})$ for all $n\geq 0$.
Rename $w_i:=v_i$ for $i\in \{1,\cdots,s\}\setminus \{1\}.$
Then $W(\theta^*)\setminus\{v\in V_{\infty}|\,\,v_1\leq v\wedge w<v_4\}=B(\{w_1,\cdots,w_s\}).$
By Theorem \ref{prodimtwothengeqz}, there exists a polynomial $P\in k[x,y]\setminus k$ satisfying $w_j(P)> 0$ for all $j=1,\cdots,s$.
For all $i=1,2$, there exists $j_i\in \{1,\cdots,s\}$ such that $v_{C_i}>v_{C_{j_i}}.$
Then we have $w_{C_i}(P)\geq w_{j_i}(P)> 0$ for all $i\in 1,2$. It follows that $C$ is an irreducible component of $\{\prod_{i=1}^sP_i=0\}$.
Apply the same argument for $f^n(C)$, $n\geq 0$, we have that $f^n(C)$ is an irreducible component of $\{\prod_{i=1}^sP_i=0\}$. It follows that $C$ is periodic.

\medskip

\subsubsection*{\textbf{2) The case that $s=1,t=2$, $v_{C_1}>v_1$ and $v_{C_2}\not\in W(\theta^*)$}}
\begin{figure}
\centering
\includegraphics[width=6.5cm]{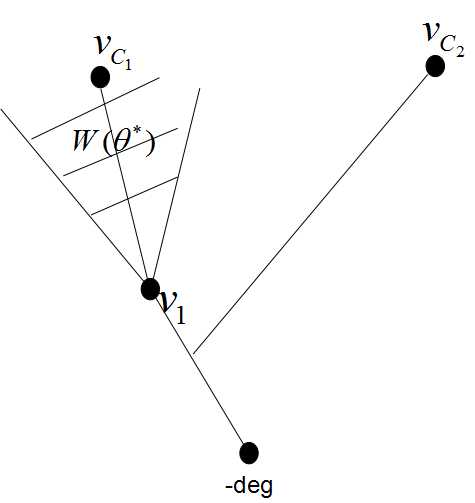}
\caption{}
\label{figc7}
\end{figure}
See Figure \ref{figc7}.

It follows that $v_1$ is divisorial with $\alpha(v_1)=0$. It follows that $\la_2/\la_1\geq 1$, and $\la_2/\la_1A(v_1)=A(v_1)+v_1(Jf)\leq A(v_1).$ Then  we have $A(v_1)\leq 0.$ By line embedding theorem, $f$ takes form $f=(F(x),G(x,y))$ where $\deg F=\la_1$ and $\deg_{y}G=\la_2/\la_1.$ Set $d_1:=\la_1$ and $d_2:=\la_2/\la_1$. Write $G$ in form $G=\sum_{i=0}^{d_2}A_i(x)y^i$ where $A_i\in k[x]$ and $A_{d_2}\neq 0$ in $k[x].$

For any $m\geq 0$, we may embed $\mathbb{A}^2_k$ in $\mathbb{F}_m$.
Let $L_{\infty}$ be the exceptional curve in $\mathbb{F}_m$ such that $v_{L_{\infty}}=v_*$ and $F_{\infty}$ the fiber of $\pi_m$ at infinity. Set $O:=L_{\infty}\cap F_{\infty}.$ By requiring $m$ large enough, we may suppose that $O\not\in C.$ The center of $C_1$ is at $L_{\infty}\setminus \{O\}$ and the center of $C_2$ is at $F_{\infty}\setminus \{O\}.$

Let $K$ be a number field such that $f$, $p$ and $C$ are all defined over $K$ and let $S$ be a finite subset of $\mathcal{M}_K$ containing $\mathcal{M}_K^{\infty}$ such that $f$ and $p$ are defined over $O_S.$

Let $h_1:C(K)\rightarrow \mathbb{R}$ be the function defined by $$h_1:(x,y)\mapsto \sum_{v\in \mathcal{M}_K}\log\max\{|x|_v,1\}$$ and  $h_2:C(K)\rightarrow \mathbb{R}$ be the function defined by $$h_2:(x,y)\mapsto \sum_{v\in \mathcal{M}_K}(\log\max\{|y|_v,1,|x|_v^m\}-\log\max\{1,|x|_v^m\}).$$ It follows that $h_1$ is a Weil height function with respect to the divisor $C\cdot {F_{\infty}}$ and $h_2$ is a Weil height function with respect to the divisor $C\cdot{L_{\infty}}.$

For all $v\in \mathcal{M}_K\setminus S$, we have $|x_n|_v\leq 1$ and $|y_n|_v\leq 1$. It follows that $$\log\max\{|y_n|_v,1,|x_n|_v^m\}-\log\max\{1,|x_n|_v^m\}=0$$ for all $v\in \mathcal{M}_K\setminus S$.

Since $\infty$ is a supperattracting point of $F,$ for all place $v\in K$, there exists $r_v>0$ such that $|f^n(x)|_v\rightarrow \infty$ for all $x\in K$ satisfying $|x|_v>r_v$ and further we may suppose that $r_v=1$ for $v\in \mathcal{M}_K\setminus S'$, where $S'$ is a finite subset of $\mathcal{M}_K.$

There exists an infinite sequence $\{n_1<n_2<\cdots\}$ of nonnegative integers such that $f^{n_i}(p)\in C$ for all $i\geq 0.$
Write $f^n(p)=(x_n,y_n)$ for all $n\geq 0$. Set $S_1\subseteq S$ consisting of places $v\in M_K$ such that $|x_n|_v\leq r_v$ for all $n\geq 0.$ Since $c(v_{C_2})\not\in L_{\infty}$, for all $v\in S$ there exists an neighborhood $U_v$ of $c(v_{C_2})$ with respect to the topology induced by $|\cdot|_v$ and $B_v\geq 0$, such that for all $(x,y)\in U_v\cap \mathbb{A}^2(K)$ we have  $|y|_v\leq B_v|x|_v^m$. For all $v\in S$, there exists $R_v>r_v+1$ such that $C\cap \{(x,y)\in \mathbb{A}^2(K)|\,\, |x|_v>R_v\}\subseteq U_v.$ By replacing $p$ by $f^n(p)$ for $n$ large enough, we may suppose that for all $v\in S\setminus S_1$, we have $|x_0|_v>R_v$.  If follows that $$\log\max\{|y_n|_v,1,|x_n|_v^m\}-\log\max\{1,|x_n|_v^m\}=\log\max\{|y_n|_v,|x_n|_v^m\}-\log(|x_n|_v^m)$$$$
\leq \log\max\{B_v|x|_v^m,|x_n|_v^m\}-\log\max(|x_n|_v^m)\leq \log\max\{B_v,1\}.$$

For all $v\in S_1$, we have $|x_n|_v\leq r_v$ for all $n\geq 0.$  There exists $D_v\geq 1$ such that $|A_i(x)|_v\leq D_v$ for all for all $i=1,\cdots,d_2$. It follows that $$|y_{n+1}|_v=|\sum_{i=0}^{d_2}A_i(x_n)y_n^i|_v$$$$\leq \sum_{i=0}^{d_2}|A_i(x_n)||y_n|^i\leq D_v\sum_{i=0}^{d_2}|y_n|^i\leq D_v(d_2+1)\max\{|y_n|,1\}^{d_2}. $$ It follows that $\max\{|y_{n+1}|_v,1\}\leq D_v(d_2+1)\max\{|y_n|,1\}^{d_2}$. It follows that there exits $D'_v\geq 0$ such that $\log\max\{|y_{n}|_v,1\}\leq (d_2+1/2)^nD'_v$ for all $v\in S_1$ and $n\geq 0.$ It follows that
$$\log\max\{|y|_v,1,|x|_v^m\}-\log\max\{1,|x|_v^m\}\leq \max\{(d_2+1/2)^nD'_v,1\}$$ for all $v\in S_1.$

Then we have $$h_2(f^{n_i}(p))\leq \sum_{v\in \mathcal{M}_K}(\log\max\{|y_{n_i}|_v,1,|x_{n_1}|_v^m\}-\log\max\{1,|x_{n_1}|_v^m\})$$
$$\leq \sum_{v\in S_1}\max\{(d_2+1/2)^{n_i}D'_v,1\}+\sum_{v\in S\setminus S_1}\log\max\{B_v,1\}$$
$$\leq (\#S_1\sum_{v\in S_1}\max\{D'_v,1\})(d_2+1/2)^{n_i}+\sum_{v\in S\setminus S_1}\log\max\{B_v,1\}$$ for all $i\geq 1.$ Since $p$ is not preperiodic, we have $S_1\neq\emptyset.$

On the other hand, there exists $B_2>0$, such that $\log |F^n(x)|_v>B_2d_1^n$ for all $v\in S\setminus S_1$ $n\geq 0$ and $x\geq r_v.$ Then we have $$h_1(f^{n_i}(p))=\sum_{v\in \mathcal{M}_K}\log\max\{|x_{n_i}|_v,1\}$$
$$\geq\sum_{v\in S\setminus S_1}\log\max\{|x_{n_i}|_v,1\}\geq \sum_{v\in S\setminus S_1}B_2d_1^n=\#(S\setminus S_1)B_2d_1^n.$$
If $\#(S\setminus S_1)=0$, $h_{1}(x_n)$ is bounded, then $x_0$ is preperiodic. Since $C$ is not a fiber, we have $C\cap \cup_{n=0}^{\infty}\{x=x_n\}$ is finite and then $p$ is preperiodic which contradics to our assumption. Then we have $\#(S\setminus S_1)>0$. It follows that $h_2(f^{n_i}(p))/h_1(f^{n_i}(p))\rightarrow 0$ as $i\rightarrow \infty$ which contradicts to Lemma \ref{lemhoverh}.

\medskip
\subsubsection*{\textbf{3) The case that $t=2$, $s=2$, $v_{C_1}\in W(\theta^*)$ and $v_{C_2}\not\in W(\theta^*)$}}
%

\begin{figure}
\centering
\includegraphics[width=6.5cm]{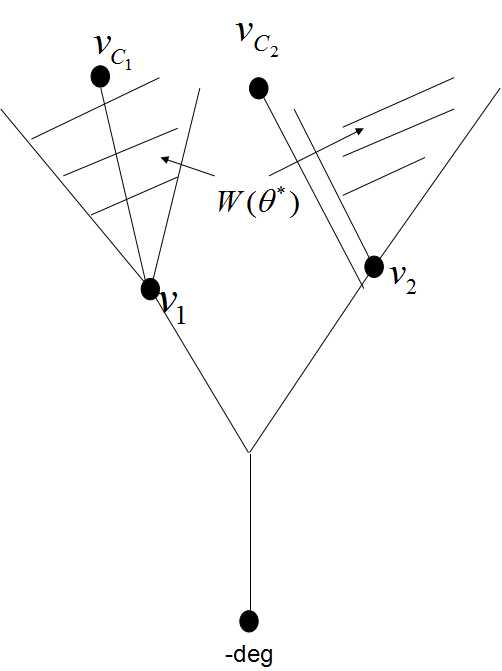}
\caption{}
\label{figc8}
\end{figure}

We may suppose that $v_{C_1}> v_1.$ See Figure \ref{figc8}.

By Theorem \ref{thmasequencecurves}, we may suppose that there exists a sequence of curves $\{C^i\}_{i\leq 0}$ with $2$ places at infinity such that
\begin{points}
\item[$\d$] $C^0=C$;
\item[$\d$] $f(C^i)=C^{i+1}$;
\item[$\d$] for all $i\in \mathbb{Z}$, the set $\{n\geq 0|f^n(p)\in C^i\}$ is infinite.
\end{points}
Let $C^j_i$'s be branches of $C^j$, we may suppose that $f(C^j_i)=C_i^{j+1}$ for $j\leq -1$ and $1\leq i\leq 2.$
Observe that $v_{C^j_1}>v_1$ for all $j\leq 0.$

The following lemma is a key ingredient of our proof in this case. It can be viewed as a modified version of Lemma \ref{lemexwtowcninw} to adapt this case.
\begin{lem}\label{lemseqtwoteqtwoneverinopen}If there exists an open set $W$ of $V_{\infty}$ containing $v_*$ such that for infinitely many $j\leq 0$ we have
$v_{C^j_i}\not\in W$ for all $i=1,\cdots,t$,
then Theorem \ref{thmsupponetwodml} holds.
\end{lem}

\proof[Proof of Lemma \ref{lemseqtwoteqtwoneverinopen}]
Set $\psi:=R_{[-\deg,v_1]}\theta^*\in \SH^+(V_{\infty})$, we have $\psi(v)=0$ for all $v\geq v_1,$ and $\langle\psi,\psi\rangle>0.$
By Lemma \ref{lemgeqm} there exists $M\leq 1$ such that for any set $B$ of valuations satisfying
\begin{points}\item[(1)]$B\setminus B(\{v_1\})$ has at most $1$ elements;
\item[(2)] $B\subseteq B(\{v_1\})\cup \{v\in V_{\infty}|\,\, \alpha(v)\leq M\}$;
\end{points}
there exists a function $\phi\in \mathbb{L}^2(V_{\infty})$
satisfying $\phi(v)=0$ for all $v\in B(B)$ and $\langle\phi,\phi\rangle>0.$

By Proposition \ref{prodmlcywt}, we may suppose that $v_{C^j_2}\not\in B(\{w^1_2,\cdots,w^s_2\})^{\circ}$ for all $j\leq 0$.
By Proposition \ref{prothetageztendstovs}, there exists $N\geq 0$ such that $\{v\in V_{\infty}|\,\,\alpha(v)\geq M\}\subseteq f^{-N}_{\d}(W)\cup B(\{w^1_2,\cdots,w^s_2\})^{\circ}.$ Set $W_1:=V_{\infty}\setminus (f^{-N}_{\d}(W)\cup B(\{w^1_2,\cdots,w^s_2\})^{\circ}).$ For all pair $w=(w_1,w_2)\in B(\{v_1\})\times W_1$, there exist $w_1'<w_1,w_2'<w_2$ and  a function $\phi_w\in \mathbb{L}^2(V_{\infty})$
satisfying $\phi_w(v)=0$ for all $v\in B(\{w_1',w_2'\})$ and $\langle\phi,\phi\rangle>0.$ Set $U_{w}:=B(\{w_1'\})^{\circ}\times B(\{w_2'\})^{\circ}$. By Theorem \ref{prodimtwothengeqz}, there exists a polynomial $P_w\in k[x,y]\setminus k$, such that $w_1'(P_w)> 0$ and $w_2'(P_w)>0$. Since $B(\{v_1\})\times W_1$ is compact, there exist $w^1,\cdots,w^m\in B(\{v_1\})\times W_1$, $m\in \mathbb{Z}^+$ such that $B(\{v_1\})\times W_1\subseteq \cup _{i=1}^mU_{w^i}.$ It follows that for all $(w_1,w_2)\in B(\{v_1\})\times W_1$, there exists $i\in \{1,\cdots, m\}$, such that $w_1(P_{w^i})> 0$ and $w_2(P_{w^i})> 0$.

There exists a infinite sequence of negative integers $\{j_1>j_2>\cdots\}$ such that $v_{C^{j_i}_2}\not\in W$ for all $i\geq 0.$ Then we have $v_{C^{j_i-N}_2}\not\in f^{-N}_{\d}(W)$ for all $i\geq 0.$  Then we have $v_{C^{j_i-N}_2}\in W_1$ for all $i\geq 0$. There exists $l_i\in \{1,\cdots,m\}$ such that $v_{C^{j_i}_1}(P_{l_i})> 0$ and $v_{C^{j_i}_2}(P_{l_i})> 0$. It follows that $P_{l_i}|_{C^{j_i}}=0$ and then $C^{j_i}$ is an irreducible component of $\{P_{l_i}=0\}.$ Since there are only finitely many irreducible components of $\{\prod_{i=1}^mP_i=0\}$, we conclude that $C$ is periodic.
\endproof

\subsubsection*{\textbf{3.1)}}
If $v_*$ is nondivisorial,
by \cite[Theorem 3.1]{Favre2011}, there exists an open set $W$ of $V_{\infty}$ containing $v_*$ such that
\begin{points}
\item[$\d$]
$v_{C^0_2}\not\in W$;
\item[$\d$]
$f_{\d}(W)\subseteq W$.
\end{points}
Then we have $W\subseteq f^j(W)$ for all $j\leq 0$. It follows that $v_{C^j_2}\not\in W$ for all $j\leq 0.$ Apply Lemma \ref{lemseqtwoteqtwoneverinopen}, we conclude our proposition in this situation.

\subsubsection*{\textbf{3.2)}} If $v_*$ is divisorial. There exists a smooth projective compactificaition $X$ of $\mathbb{A}^2$ containing a divisor $E$ satisfying $v_E=v_*$. By \cite[Lemma 4.6]{Favre2011}, we may suppose that for any point $t$ in $I(f)\cap E$, $t$ is not a periodic point of $f|_E$.

%
%
There exists a neighborhood $W$ of $v_*$ in $V_{\infty}$ such that

\begin{points}
\item for all valuation $v\in W$, $d(f,v)>0$ and the center of $v$ is contained in $E$;
\item for any point $t\in E$, we have $f_{\d}(U(t)\cap W)\subseteq U(f|_{E}(t))$.
\end{points}
For any valuation $v\in W$, denote by $c(v)$ the center of $v$ in $E.$
By Lemma \ref{lemseqtwoteqtwoneverinopen}, there exists $j_0\leq 0$ such that $v_{C^j_2}\in W$ for all $j\leq j_0.$ By replacing $C$ by $C^{j_0}$, we may suppose that $j_0=0.$

\smallskip

 \subsubsection*{ \textbf{3.2.1)}} We first treat the case that $\deg(f|_E)=1$. By replacing $f$ by a positive iterate, we may suppose that all periodic points of $f$ are fixed.
 By Lemma \ref{lemdegeeonecjcjpolo}, we may suppose that $W$ is a nice neighborhood of $v_*$. Recall that $W$ satisfies the following properties:
\begin{points}
\item for all valuation $v\in W$, $d(f,v)>0$ and the center of $v$ is contained in $E$;
\item for any point $t\in E$, we have $f_{\d}(U(t)\cap W)\subseteq U(f|_{E}(t))$;
\item for all $j\leq 0$ there exists $i\in \{1,\cdots,s\}$ satisfying $v_{C^j_i}\in W$,
we have $\deg f|_{C^j}\leq \la_1$ for all $j\leq -1$;
\item its boundary $\partial W$ is finite;
\item for any fixed point $x\in E$, $f_{\d}(U(x)\cap W)\subseteq U(x)\cap W$.
\end{points}

  %

For all $j\leq 0$ and $i=1,2$, set $m^j_i:=(C^j_i\cdot l_{\infty})$.
\begin{lem}\label{lemcjthetastabound}There exists $B>0$ such that $\left(\theta^*\cdot (\sum_{i=1}^2m_i^jZ_{v_{C^j_i}})\right)\leq B$ for all $j\leq 0.$
\end{lem}

\rem This Lemma holds also in the case that all branches $C_i$ of $C$ are not contained in $W(\theta^*)$ by the same proof.
\endrem

Then we have $m_1^j+m_2^j=\deg C^j$ and $B\geq \sum_{i=1}^2m_i^j\theta^*(v_{C^j_i})=m_2^j\theta^*(v_{C^j_i}).$
By Proposition \ref{prodmlcywt}, we may suppose that $v_{C^j_2}\not\in B(\{w^1_2,\cdots,w^s_2\})^{\circ}$ for all $j\leq 0$.
Since $V_{\infty}\setminus V$ is compact and $\theta^*$ is continuous, there exists $\delta>0$ such that $\theta^*\geq \delta$ on $V_{\infty}\setminus B(\{w^1_2,\cdots,w^s_2\})^{\circ}$. It follows that $B\geq m_2^j\theta^*(v_{C^j_i})\geq m_2^j\delta$. It follows that $m_2^j\leq \delta^{-1}B$ for all $j\leq 0.$

Since $\langle R_{[-\deg,v_1]}\theta^*,R_{[-\deg,v_1]}\theta^*\rangle>0,$
by Proposition \ref{prodimtwothengeqz}, there exists a polynomial $P$ satisfying $v_1(P)>0.$ Set $r:=v_1(P)$, then $v_{C^j_1}(P)>r$ for all $j\leq 0.$
Then $P|_{C^j}$ has at least $m_1^jr$ zeros but $P|_{C^j}$ has $\max\{0,-m_2^jv_{C^j_2}(P)\}$ poles. Observe that $-m_2^jv_{C^j_2}(P)\leq m_2^j\deg(P).$ If $m_1^jr>m_2^j\deg(P)$, then $P|_{C^j}=0$ and then $C^j$ is an irreducible component of $\{P=0\}$. Suppose that $C$ is not periodic. By replacing replacing $C$ be some $C^j$ for $j$ negative enough, we may suppose that $m_1^jr\leq m_2^j\deg(P)$ for all $j\leq 0$. Then we have $m_1^j\leq r^{-1}\deg(P)m_2^j\leq r^{-1}\deg(P)\delta^{-1}B$ and then $\deg C^j=m^j_1+m^j_2\leq (1+r^{-1}\deg(P))\delta^{-1}B$ for all $j\leq 0.$ We conclude our theorem by Proposition \ref{procurveboundeddegfibration}.

\proof[Proof of Lemma \ref{lemcjthetastabound}]By Lemma \ref{lembrafact}, we have
$m^{j-1}_id(f,v_{C^{j-1}_i})=\deg(f|_{C^{j-1}})m_{C^{j}_i}$ for all $j\leq 0$ and $i=1,2$. It follows that
$$\left(\theta^*\cdot (\sum_{i=1}^2m_i^{j-1}Z_{v_{C^{j-1}_i}})\right)=\la_1^{-1}\left(f^*\theta^*\cdot (\sum_{i=1}^2m_i^{j-1}Z_{v_{C^{j-1}_i}})\right)$$
$$=\la_1^{-1}\left(\theta^*\cdot f_*(\sum_{i=1}^2m_i^{j-1}Z_{v_{C^{j-1}_i}})\right)=\la_1^{-1}\left(\theta^*\cdot (\sum_{i=1}^2m_i^jd(f,v_{C^{j-1}_i})Z_{v_{C^{j}_i}})\right)$$
$$=\la_1^{-1}\deg(f|_{C^{j-1}})\left(\theta^*\cdot (\sum_{i=1}^2m_{C^{j}_i}Z_{v_{C^{j}_i}})\right)\leq \left(\theta^*\cdot (\sum_{i=1}^2m_{C^{j}_i}Z_{v_{C^{j}_i}})\right).$$ Then $B:=\left(\sum_{i=1}^2m_{C^{0}_i}Z_{v_{C^{0}_i}})\right)$ is what we require.
\endproof

\smallskip

\subsubsection*{ \textbf{3.2.2)}}  Then we suppose that $\deg f|_E\geq 2.$
There exists a number field $K$ such that $X,E,f,C,p$ are all defined in $K.$ For all $j\leq 0$, since $C^j$ is rational and contains infinitely many $K$-points, we have that $C^j$ is defined over $K.$ Then we have that $c(v(C^j_2))\in f^{j}(c(v(C^0_2)))$ is defined over $K$.  By Northcott property, we have that the set $\{c(v(C^j_2))\}_{j\leq 0}$ is finite. By replacing $f$ by a suitable iterate, we may suppose that $c(v(C^j_2))=c(v(C^0_2))$ for all $j\leq 0.$
Set $q:=c(v(C^j_2))$.
Let $W'$ be a neighborhood of $v_*$ in $W$ satisfying
\begin{points}
\item[$\d$]$v_{C^0_1}\not\in W'$;
\item[$\d$]$f_{\d}(U(q)\cap W')\subseteq U(q)\cap W'$.
\end{points}
It follows that $v_{C^j_2}\not\in W'$ for all $j\leq 0.$ By Lemma \ref{lemseqtwoteqtwoneverinopen}, we conclude our theorem in this case.

\medskip

\subsubsection*{\textbf{4})} Finally we treat the case that $v_{C_i}\not\in W(\theta^*)$ for all branches $C_i$ of $C$ at infinity.

By Theorem \ref{thmasequencecurves}, we may suppose that there exists a sequence of curves $\{C^j\}_{j\in \Z}$ with at most $2$ branches at infinity such that
\begin{points}
\item[$\d$] $C^0=C$;
\item[$\d$] $f(C^i)=C^{i+1}$;
\item[$\d$] $v_{C^i_j}\not\in W(\theta^*)$ for $j=1,\cdots,s$;
\item[$\d$] for all $i\in \mathbb{Z}$, the set $\{n\geq 0|f^n(p)\in C^i\}$ is infinite.
\end{points}
Let $C^j_i$'s be branches of $C^j$, we may suppose that $f(C^j_i)=C_i^{j+1}$ for $j\leq -1$ and $i=1,2.$ Since for branches $C^j_i$, we have $v_{C^j_i}\not\in  W(\theta^*)$, we have $d(f,v_{C^j_i})>0$. It follows that the number of branches of $C^j$ at infinity are the same for all $j\leq 0$.

\rem When $t=1$, it is possible that there exists $j_0\leq -1$ such that the number of branches of $C^{j}$ at infinity equals to 2 for all $j\leq j_0$. In this case, we may replace $C$ by $C^{j_0}$ and then we may also suppose that the number of branches of $C^j$ at infinity are the same for all $j\leq 0$.
\endrem

The following lemma is a key ingredient of our proof in this case. It plays the same role as Lemma \ref{lemexwtowcninw} does in the case $\#\Supp\Delta \theta^*\geq 3$.

\begin{lem}\label{lemsupplesstwotwonotinwnotinw}Let $L$ be a nonnegative integer. If there exists an open set $W$ of $V_{\infty}$ containing $v_*$ such that for infinitely many $j\leq 0$ we have
$v_{C^j_i}\not\in \cap_{l=0}^L f^{-l}_{\d}(W)$,
then the pair $(\A^2_k,f)$ satisfies the DML property for the curve $C$.
\end{lem}

\rem This Lemma holds also when $t=1$ by the same argument as in the case $t=2$.
\endrem

\proof[Proof of Lemma \ref{lemsupplesstwotwonotinwnotinw}]
By Proposition \ref{prodmlcywt}, we may suppose that for all $j\leq 0$ and all branches $C^j_i$ of $C^j$ at infinity, $v_{C^j_i}\not\in B(\{w^1_2,\cdots,w^s_2\})^{\circ}$.

By Proposition \ref{prothetageztendstovs}, there exists $N\geq 0$ such that $\{v\in V_{\infty}|\,\,\alpha(v)\geq M\}\subseteq (\cap_{l=0}^L f^{-N-l}_{\d}(W))\cup B(\{w^1_2,\cdots,w^s_2\})^{\circ}.$

Set $W_1:=V_{\infty}\setminus ((\cap_{l=0}^L f^{-N-l}_{\d}(W))\cup B(\{w^1_2,\cdots,w^s_2\}))^{\circ}.$ For all pair $w=(w_1,w_2)\in W_1^{2}$, there exist $w_i'<w_i$ for all $i=1,2$ and  a function $\phi_w\in \mathbb{L}^2(V_{\infty})$
satisfying $\phi_w(v)=0$ for all $v\in B(\{w_1',w_2'\})$ and $\langle\phi,\phi\rangle>0.$ Set $U_{w}:=\prod_{i=1}^2B(\{w_i'\})^{\circ}$. By Theorem \ref{prodimtwothengeqz}, there exists a polynomial $P_w\in k[x,y]\setminus k$, such that $w_i'(P_w)> 0$  for all $i=1,\cdots, s$. Since $W_1^2$ is compact, there exist $w^1,\cdots,w^m\in W_1^2$, $m\in \mathbb{Z}^+$, such that $W_1^2\subseteq \cup _{i=1}^mU_{w^i}.$ It follows that for all $(w_1,w_2)\in W_1^2$, there exists $i\in \{1,\cdots, m\}$, such that $w_1(P_{w^i})> 0$ and $w_j(P_{w^i})> 0$ for all $j=1,2$.

There exists a infinite sequence of negative integers $\{j_1>j_2>\cdots\}$ such that $v_{C^{j_i}_2}\not\in \cap_{l=0}^L f^{-l}_{\d}W$ for all $i\geq 0.$ Then we have $v_{C^{j_i-N}_2}\not\in \cap_{l=0}^L f^{-l-N}_{\d}W$ for all $i\geq 0.$  Then we have $v_{C^{j_i-N}_2}\in W_1$ for all $i\geq 0$. There exists $l_i\in \{1,\cdots,m\}$ such that $v_{C^{j_i}_r}(P_{l_i})\geq 0$ for all $r=1,\cdots,t$ and $v_{C^{j_i}_1}(P_{l_i})> 0$. It follows that $P_{l_i}|_{C^{j_i}}=0$ and then $C^{j_i}$ is an irreducible component of $\{P_{l_i}=0\}.$ Since there are only finitely many irreducible components of $\{\prod_{i=1}^mP_i=0\}$, we conclude that $C$ is periodic.
\endproof

\subsubsection*{\textbf{4.1)}}  If $v_*$ is nondivisorial,
by \cite[Theorem 3.1]{Favre2011}, there exists an open set $W$ of $V_{\infty}$ containing $v_*$ such that
\begin{points}
\item[$\d$]
$v_{C^0_i}\not\in W$ for all $i=1,2$;
\item[$\d$]
$f_{\d}(W)\subseteq W$.
\end{points}
Then we have $W\subseteq f^j(W)$ for all $j\leq 0$. It follows that $v_{C^j_i}\not\in W$ for all $i=1,2$ and $j\leq 0.$ Apply Lemma \ref{lemsupplesstwotwonotinwnotinw}, we conclude our proposition in this situation.

\rem In the case $t=1$, we can use the same argument as in the case $t=2$.
\endrem

\smallskip

\subsubsection*{\textbf{4.2)}} If $v_*$ is divisorial. There exists a smooth projective compactificaition $X$ of $\mathbb{A}^2$ containing a divisor $E$ satisfying $v_E=v_*$. We may suppose that for any point $x$ in $I(f)\cap E$, $x$ is not a periodic point of $f|_E$.

%
%
%

There exists a neighborhood $W$ of $v_*$ in $V_{\infty}$ such that
\begin{points}
\item for all valuation $v\in W$, $d(f,v)>0$ and the center of $v$ is contained in $E$;
\item for any point $x\in E$, we have $f_{\d}(U(x)\cap W)\subseteq U(f|_{E}(x))$;
\end{points}
For any valuation $v\in W$, denote by $c(v)$ the center of $v$ in $E.$
By Lemma \ref{lemseqtwoteqtwoneverinopen}, there exists $j_0\leq 0$ such that for all $j\leq j_0$, there exists a branch $C^j_i$ of $C^j$ at infinity such that $v_{C^j_i}\in W$. By replacing $C$ by $C^{j_0}$, we may suppose that $j_0=0.$

\smallskip

\subsubsection*{\textbf{4.2.1)}} We first treat the case that $\deg(f|_E)=1$. By Lemma \ref{lemdegeeonecjcjpolo}, we may suppose that $W$ is a nice neighborhood. By Lemma \ref{lemseqtwoteqtwoneverinopen} and by replacing $C$ by a suitable $C^{j_0}$, $j_0\leq 0$, we may suppose that
 for all $j\leq 0$, there exists a branch $C^j_i$ of $C^j$ at infinity such that $v_{C^j_i}\in W$.
Then we have $\deg f|_{C^j}\leq \la_1$ for all $j\leq -1$.

For all $j\leq 0$ and $i=1,2$, set $m^j_i:=(C^j_i\cdot l_{\infty})$. Then we have $m_1^j+m_2^j=\deg C^j$ for all $j\leq 0.$

By Lemma \ref{lemcjthetastabound}, there exists $B>0$ such that $B\geq \sum_{i=1}^2m_i^j\theta^*(v_{C^j_i}).$
By Proposition \ref{prodmlcywt}, we may suppose that $v_{C^j_i}\not\in B(\{w^1_2,\cdots,w^s_2\})^{\circ}$ for all $j\leq 0$ and $i=1,2$.
Since $V_{\infty}\setminus V$ is compact and $\theta^*$ is continuous, there exists $\delta>0$ such that $\theta^*\geq \delta$ on $V_{\infty}\setminus B(\{w^1_2,\cdots,w^s_2\})^{\circ}$. It follows that $B\geq \sum_{i=1}^2m_i^j\theta^*(v_{C^j_i})\geq \delta\sum_{i=1}^2m_i^j=\delta\deg(C^j)$. It follows that $\deg(C^j)\leq \delta^{-1}B$ for all $j\leq 0.$
Then we conclude our theorem by Proposition \ref{procurveboundeddegfibration}.

\rem In the case $t=1$, we can use the same argument as in the case $t=2$.
\endrem

\smallskip

\subsubsection*{\textbf{ 4.2.2)}} Then we may suppose that $\deg(f|_E)\geq 2.$
There exists a number field $K$ such that $X,E,f,C,p$ are all defined in $K.$ For all $j\leq 0$, since $C^j$ is rational and contains infinitely many $K$-points, we have that $C^j$ is defined over $K.$ Then if there exists a point $x\in C^j\cap E$, we have $[K(x):K]\leq 2.$ Let $P$
be the set of points $x\in E$ such that $f|_E^n(x)\in I(f)$ for some $n\geq 0$ and satisfying $[K(x):K]\leq 2.$ Observe that for all $x\in P$, $x$ is not periodic. By Northcott property, we have that $P$ is a finite set. Set $L:=\# P$.

By Lemma \ref{lemsupplesstwotwonotinwnotinw}, we may suppose that there exists $j_0\leq 0$ such that for all $j\leq j_0,$ there exists a branch $C^j_i$ of $C^j$ at infinity satisfying $v_{C^j_i}\in W.$ By replacing $C$ by $C^{j_0}$, we may suppose that $j_0=0.$

\rem When $t=1$, we are always in the following case 4.2.2.1) and the argument is the same as the case $t=2$.
\endrem

\subsubsection*{\textbf{ 4.2.2.1)}} If there exists $j_0\leq 0$ for which $v_{C^j_i}\in W$ for all branches $C^j_i$ of $C^j$ at infinity and $j\leq j_0,$ by replacing $C$ by $C^{j_0}$, we may suppose that $j_0=0.$

For $i=1,2$ and all $j\leq 0$, we have $[K(c(v_{C_i^j})):K]\leq 2$. By Northcott property, the set $\{c(v_{C_i^j}),j\leq 0\}$ is finite. It follows that $c(v_{C^0_i})$ is periodic for $i=1,2.$ By replacing $f$ by some positive iterate, we may suppose that there exists $x_i\in E$ which is fixed by $f|_E$ and satisfying $x_i=c(v_{C^j_i})$ for all $j\leq 0.$ Let $W'$ be a neighborhood of $v_*$ in $W$ satisfying
\begin{points}
\item[$\d$]$v_{C^0_i}\not\in W'$ for $i=1,\cdots,t$;
\item[$\d$]$f_{\d}(U(x_i)\cap W')\subseteq U(x_i)\cap W'$ for $i=1,\cdots,t$.
\end{points}

It follows that $v_{C^j_i}\not\in W'$ for all $i=1,2$ and $j\leq 0.$ By Lemma \ref{lemsupplesstwotwonotinwnotinw}, we conclude our theorem in this case.

\subsubsection*{\textbf{ 4.2.2.2)}}

If there exists $i_0\in \{1,2\}$ and $j_0\leq 0$ such that $v_{C^j_{i_0}}\in W$ for all $j\leq j_0,$ we may suppose that $i_0=1$ and by replacing $C$ by $C^{j_0}$, we may suppose that $j_0=0.$

By the argument in the previous paragraph, we may suppose that there exists an infinite set $S$ of index $j\leq 0$ such that $v_{C^j_2}\not\in W.$ By the same argument in the previous paragraph, we may suppose that there exists $x\in E$ which is fixed by $f|_E$ and satisfying $x=c(v_{C^j_1})$ for all $j\leq 0.$ Let $W'$ be a neighborhood of $v_*$ in $W$ satisfying
\begin{points}
\item[$\d$]$v_{C^0_1}\not\in W'$;
\item[$\d$]$f_{\d}(U(x)\cap W')\subseteq U(x)\cap W'$.
\end{points}
It follows that $v_{C^j_1}\not\in W'$ for all $j\leq 0.$ By Lemma \ref{lemsupplesstwotwonotinwnotinw}, we conclude our theorem in this case.

\subsubsection*{\textbf{ 4.2.2.3)}}
Otherwise, there exists $j_i\leq 0$ such that $v_{C^{j_i}_i}\not\in W$ for all $i=1,2$ Pick $j_0=\min\{j_1,j_2\}-1$. It follows that for all $i=1,2$ and $j\leq j_0$, if $v_{C^j_i}\in W$, we have $c(v_{C^j_i})\in P.$ By replacing $C$ by $C^{j_0}$, we may suppose that $j_0=0.$

If there exists $i\in\{1,2\}$ and $j\leq -L$, such that $v_{C^j_i}\in \cap_{l=0}^Lf_{\d}^{-l}(W)$, then we have $\{c(v_{C^{j}_i}),\cdots,c(v_{C^{j+L}_i})\}\subseteq P.$ Since there are not periodic points in $P$, we get a contradiction.

It follows that $v_{C^j_i}\not\in \cap_{k=0}^Lf_{\d}^{-k}(W)$ for all $i=1,2$ and $j\leq -L.$
Then we conclude our theorem by Lemma \ref{lemsupplesstwotwonotinwnotinw}.
\endproof

\addcontentsline{toc}{section}{}

\printindex

\bibliography{dd}

\end{document}